\documentclass[11pt]{amsart}

\usepackage{epsfig, amsthm, amsfonts, color}

\newcommand{\alg}{\ensuremath{\mathcal{A}}}

\newcommand{\df}{\ensuremath{\partial}}
\newcommand{\dfi}{\ensuremath{\partial_{int}}}
\newcommand{\dfe}{\ensuremath{\partial_{ext}}}

\newcommand{\reeb}{\ensuremath{X_\alpha}}
\newcommand{\ps}[1]{\ensuremath{\mathbf{#1}}}
\newcommand{\nov}{\ensuremath{T}}
\newcommand{\ca}[1]{\ensuremath{\mathcal{C}^0(\mathcal{A}_{#1}, \df_{#1})}}
\DeclareMathOperator{\image}{Im}

\newcommand{\cc}{\ensuremath{\mathbb{C}}}
\newcommand{\rr}{\ensuremath{\mathbb{R}}}
\newcommand{\zz}{\ensuremath{\mathbb{Z}}}

\newcommand{\nn}{\ensuremath{\mathbb{N}}}

\theoremstyle{plain}
\newtheorem{thm}{Theorem}[section]
\newtheorem{cor}[thm]{Corollary}
\newtheorem{lem}[thm]{Lemma}
\newtheorem{claim}[thm]{Claim}
\newtheorem{prop}[thm]{Proposition}

\theoremstyle{definition}
\newtheorem{defn}[thm]{Definition}
\newtheorem*{ex}{Example}

\theoremstyle{remark}
\newtheorem*{rem}{Remark}

\begin{document}

\title{Invariants of Legendrian Knots in Circle Bundles}

\author[J. Sabloff]{Joshua M. Sabloff} \address{University of
  Pennsylvania, Philadelphia, PA  19104} \email{sabloff@math.upenn.edu}
\urladdr{http://math.upenn.edu/\textasciitilde sabloff} 



\begin{abstract}
  Let $M$ be a circle bundle over a Riemann surface that supports a
  contact structure transverse to the fibers.  This paper presents a
  combinatorial definition of a differential graded algebra (DGA) that
  is an invariant of Legendrian knots in $M$.  The invariant
  generalizes Chekanov's combinatorial DGA invariant of Legendrian
  knots in the standard contact $3$-space using ideas from Eliashberg,
  Givental, and Hofer's contact homology.  The main difficulty lies in
  dealing with what are ostensibly $1$-parameter families of
  generators for the DGA; these are solved using ``Morse-Bott''
  techniques.  As an application, the invariant is used to distinguish
  two Legendrian knots that are smoothly isotopic, realize a
  non-trivial homology class, but are not Legendrian isotopic.
\end{abstract}

\maketitle

\setcounter{tocdepth}{2}
\tableofcontents

\newpage
\section{Introduction}
\label{sec:intro}

\subsection{Legendrian Knots}
\label{sec:problems}

A Legendrian knot is a smooth embedding of the circle into a contact
manifold that is everywhere tangent to the contact
planes.\footnote{See Section~\ref{sec:basics} for more of the basic
  notions of contact geometry.}  One fundamental problem in Legendrian
knot theory is to determine when two Legendrian knots are not (or are)
isotopic through Legendrian knots, even if they are isotopic as smooth
knots.  In other words, the problem is to find effective invariants of
Legendrian isotopy.  For null-homologous knots, the first step is to
define the the Thurston-Bennequin and rotation numbers (the so-called
``classical invariants'').  Eliashberg and Fraser proved that these
two invariants classify unknots in tight contact $3$-manifolds
\cite{yasha-fraser}, and Etnyre and Honda proved that the classical
invariants also classify the Legendrian figure eight knots and torus
knots \cite{etnyre-honda:knots}.

The classical invariants are not the whole story, however.  Chekanov
\cite{chv} built an invariant of Legendrian isotopy that distinguishes
two Legendrian $5_2$ knots that have the same classical invariants.
Chekanov's invariant (commonly referred to as the Chekanov-Eliashberg
algebra) is a non-commutative differential graded algebra (DGA) over
$\zz_2$, freely generated by the double points of the $xy$ diagram of
a Legendrian knot.  The differential comes from counting immersed
polygons with corners at the double points and whose edges lie in the
diagram.  The ``stable tame isomorphism'' type --- and consequently
the homology --- of the DGA is invariant under Legendrian isotopy.  It
is difficult to extract information from this invariant, but Chekanov
defined a ``linearized'' theory of the DGA that, when it exists, has
proven to be useful in distinguishing Legendrian knots that have the
same classical invariants; also see Fuchs et al.
\cite{efm,fuchs:augmentations}.  In \cite{lenny:computable}, Ng
generated more powerful methods for mining the DGA for information and
used them to distinguish Legendrian knots that the linearized DGA
cannot.

At the same time that Chekanov was developing his invariant,
Eliashberg, Givental, and Hofer \cite{egh} constructed a geometric
invariant for Legendrian knots.\footnote{Givental's contributions to
  the theory came after the name ``Chekanov-Eliashberg algebra'' came
  into common usage, but the reasons that Hofer's name was omitted are
  lost in the sands of time.}  Their \textbf{relative contact
  homology} is a non-commutative DGA generated by the Reeb chords of
the knot (in a sufficiently generic setup in which the chords are
isolated).  For a contact manifold $(M, \alpha)$, let $(M \times \rr,
d(e^t \alpha))$ be its symplectization; if $L$ is a Legendrian knot in
$M$, then $L \times \rr$ is a Lagrangian submanifold of the
symplectization.  The differential for the DGA is defined by counting
rigid $J$-holomorphic disks in $M \times \rr$ whose boundaries lie in
$L \times \rr$.  In \cite{ens}, it was proven that Chekanov's
invariant is a combinatorial translation of Eliashberg, Givental, and
Hofer's contact homology for the standard contact structure on
$\rr^3$.

Whereas the classification results or Eliashberg-Fraser and
Etnyre-Honda hold in any tight contact manifold, thus far,
non-classical invariants have been rendered computable only for the
standard structure on $\rr^3$ and for the space of contact elements to
the plane.  The results in this paper use the ideas of relative
contact homology to extend the range of tight manifolds for which
there is a combinatorially computable non-classical invariant. In
particular, the goal of this paper is to define a combinatorial theory
for Legendrian knots in circle bundles that have contact structures
transverse to their fibers.

\subsection{Combinatorial Contact Homology}
\label{sec:comb-intro}

The key feature used in defining a combinatorial translation of
relative contact homology is that the projection along the Reeb flow
is a fibration over a complex base. For example, in the $\rr^3$ case,
the Reeb field points along the positive $z$ direction.  Thus, the
projection along the Reeb flow is equivalent to projection to the $xy$
plane.  Consequently, the Reeb chords of a knot correspond to the
double points of its $xy$ diagram, which shows that the generators of
the Chekanov-Eliashberg algebra are the same as those of the relative
contact homology algebra.  The fact that the projection from the
almost complex manifold $(M \times \rr, J)$ to the base of the
fibration is holomorphic leads to a correspondence between spaces of
immersed polygons used to define the combinatorial differential and
the moduli spaces of rigid holomorphic disks in the geometric theory.
See Section $7$ of \cite{ens} for more details.

Circle bundles with a contact structure transverse to their fibers
constitutes another class of contact manifolds for which the Reeb flow
induces a fibration.  For these \textbf{contact circle bundles}, a
contact form may be chosen so that the Reeb field points along the
fibers.\footnote{See Section~\ref{sec:basics} for a more detailed
  construction of these contact structures.}  Thus, the projection
along the Reeb flow is the same as the bundle projection to the base.
Eliashberg, Givental, and Hofer's theory says that the DGA of a
Legendrian knot $L$ should be generated by the Reeb chords.  There are
two types of Reeb chords for a Legendrian knot in a contact circle
bundle.  The first type come from double points of the projection of
$L$ to the base. Over a double point, there are chords that start on
one strand of $L$, possibly wrap around the fiber a few times, and
finish on the other strand.  These are isolated, and to each double
point there corresponds two sets of chords, depending on the
starting strand.  Each set is indexed by the winding number of the
chord around the fiber.  The second type live over every point of $L$:
they start and end at the same point, traversing the fiber at least
once.  For each winding number around the fiber, there is a
$1$-parameter family of these chords, parameterized by $L$.  This would
seem to indicate that the DGA would have to be uncountably generated.

In order to deal with this complication, ``Morse-Bott'' methods must
be applied to relative contact homology.  In particular, each
$1$-parameter family is replaced by the critical points of a Morse
function on that family, and the differential is adjusted to reflect
this perturbation.  These techniques lead to the following theorem:

\begin{thm} \label{thm:overall}
  Let $(E, \alpha)$ be a contact circle bundle.  Let $L$ be a
  Legendrian knot in $(E, \alpha)$.  Then there is a
  combinatorially-defined filtered DGA whose ``stable tame
  isomorphism'' type is invariant under Legendrian isotopy of $L$.
\end{thm}

The combinatorial definition of the invariant occupies
Section~\ref{sec:defn}, and Theorem~\ref{thm:dga} states that the
invariant is a genuine DGA.  The definitions necessary to understand
the invariance properties of the DGA are set out in
Section~\ref{sec:isomorphisms}, culminating in
Theorem~\ref{thm:invariance}.

As an application, the invariant of Theorem~\ref{thm:overall} may be
used to prove:

\begin{prop} \label{prop:overall-ex}
  There exist examples of Legendrian knots in circle bundles that are
  smoothly isotopic but not Legendrian isotopic that, additionally,
  satisfy one of the following two topological properties:
  \begin{enumerate}
  \item The knot has nontrivial degree in the fiber.
  \item The projection of the knot to the base $F$ represents a
  nontrivial homology class.
  \end{enumerate}
\end{prop}

The remainder of the paper is divided into five sections.
Section~\ref{sec:geom} surveys the basic notions of contact geometry
and, in particular, the geometry of contact structures on circle
bundles.  The section ends with a combinatorial description of the
projection of a Legendrian knot to the base, paying particular
attention to the case of knots in lens spaces.

In Section~\ref{sec:defn}, the definition of the DGA appears in three
parts: first, the underlying algebra is defined in
Section~\ref{sec:algebra}.  Second, the differential is described in
Section~\ref{sec:diffl}, with a simple example to illustrate how to
compute it in Section~\ref{sec:simple-ex}.  Finally, the algebraic
notions necessary for the statement of invariance are detailed in
Section~\ref{sec:isomorphisms}.  As mentioned above, the two most
important theorems are Theorems~\ref{thm:dga} and
\ref{thm:invariance}, which combine to make Theorem~\ref{thm:overall}
precise.  The proofs of these theorems are delayed until
Sections~\ref{sec:d2} and \ref{sec:invariance}, respectively.
Section~\ref{sec:apps} contains the computations necessary to prove
Proposition~\ref{prop:overall-ex}.

\section{Geometry of Contact Circle Bundles}
\label{sec:geom}

\subsection{Basic Notions of Contact Geometry}
\label{sec:basics}

Let $M$ be a closed, oriented $3$-manifold.  A \textbf{contact
  structure} $\xi$ on $M$ is a completely non-integrable tangent
$2$-plane field. If $\xi$ is the kernel of a $1$-form $\alpha$, then
the non-integrability condition is equivalent to $\alpha \wedge
d\alpha$ being nowhere vanishing.  Such a contact structure is called
\textbf{co-oriented} and $\alpha$ is called a \textbf{contact form}. A
co-oriented contact structure is called \textbf{positive} if $\alpha
\wedge d\alpha$ gives the correct orientation on $M$.  A contact form
picks out a special vector field transverse to $\xi$ called the
\textbf{Reeb field} \reeb.  This vector field is is defined by the
equations
\begin{equation*}
  \begin{split}
    d\alpha(X_\alpha, \cdot) &= 0, \\
    \alpha(X_\alpha) &= 1.
  \end{split}
\end{equation*}

Darboux's Theorem says that every contact form is locally isomorphic
to the standard contact form on $\rr^3$:
\begin{equation*}
  \alpha_0 = dz+xdy.
\end{equation*}
Thus, only the global properties of contact manifolds are of interest.
To organize the study of global properties, Eliashberg divided contact
structures on $3$-manifolds into two classes: \textbf{overtwisted},
which contain an embedded disk $D$ that is tangent to $\xi$ along
$\partial D$, and \textbf{tight}, which do not.  Eliashberg proved
that overtwisted structures are classified up to contact isotopy by
the homotopy class of their underlying $2$-plane fields
\cite{yasha:overtwisted}.  Tight structures are more rigid, as
evidenced, for example, by the existence of a unique tight contact
structure on $S^3$ \cite{yasha:20yrs}, by the existence of a manifold
that does not admit a tight structure
\cite{etnyre-honda:non-existence}, and by the Bennequin inequality
(see below).

One effective way to study global properties of contact manifolds is
to look at the Legendrian knots they support.  A \textbf{Legendrian
  knot} is an embedded circle $L \subset M$ that is always tangent to
$\xi$.  An ambient isotopy of $L$ through other Legendrian knots is a
\textbf{Legendrian isotopy}.  Note that Legendrian knots are
plentiful; for example, any smooth knot can be continuously
approximated by a Legendrian knot.

As mentioned in the introduction, there are two ``classical''
invariants for null-homologous Legendrian knots up to Legendrian
isotopy.  The first classical invariant is the
\textbf{Thurston-Bennequin number} $tb(L)$, which measures the
twisting of the contact planes around the knot $L$.  More precisely,
let $\hat{L}$ result from pushing $L$ out a small distance along a
vector field that is transverse to $\xi$ along $L$.  Define $tb(L)$ to
be the linking number of $L$ and $\hat{L}$.  The second classical
invariant, the \textbf{rotation number}, is defined for
\emph{oriented} Legendrian knots.  It measures the twisting of the
tangent direction to $L$ inside $\xi$.  More precisely, let $\Sigma$
be a Seifert surface for $L$.  This means that $\partial \Sigma = L$,
so $[\Sigma]$ is a class in $H_2 (M,L;\zz)$.  Trivialize $\xi$ over
$\Sigma$; then $r(L; [\Sigma])$ is the winding number of the oriented
tangent direction to $L$ with respect to this trivialization.  In a
tight contact manifold, the classical invariants are restricted by the
Bennequin inequality, which was originally proved by Bennequin
\cite{bennequin} and generalized by Eliashberg \cite{yasha:knots} to:
\begin{equation*}
  tb(L) + \left| r(L; [\Sigma]) \right| \leq -\chi(\Sigma).
\end{equation*}

The contact manifolds of interest in this paper are circle bundles $E$
over closed Riemann surfaces $F$ with contact structures transverse to
the fibers.  So long as the Euler number $e(E)$ is negative, Giroux
\cite{giroux:bundles} and Honda \cite{honda:classification-2} proved
that such contact structures exist, are tight, and are unique up to
contact isotopy.\footnote{The existence result is really a consequence
  of the ``Milnor-Wood'' inequalities for contact structures; see
  \cite{confol, giroux:bundles, sato-tsuboi}.}  If, in addition, the
contact structure is invariant along the fibers, the Lutz \cite{lutz}
proved that it is unique up to equivariant contactomorphism.

The following construction of the unique invariant contact structure
transverse to the fibers of $E \overset{\pi}{\to} F$ will be used
throughout this paper.  Let $\mathcal{E} \overset{\pi}{\to} F$ be a
Hermitian line bundle and let $E$ be its unit circle bundle.  Let
$\alpha$ be the restriction of a unitary connection form on
$\mathcal{E}$ to $E$.  Thus, the curvature $\Omega \in \Omega^2(F)$ is
given by:
\begin{equation}
  \label{eqn:curvature}
  \pi^* \Omega = i d\alpha.
\end{equation}
The connection form is a positive contact form if and only if its
curvature is strictly negative.  It follows that $e(E) < 0$.  The Reeb
field for $\alpha$ points along the fibers.  Note that $E$ is
holomorphically filled by the unit disk bundle of $\mathcal{E}$, so it
is tight \cite{yasha:filling}.

\begin{defn} \label{defn:ccb}
  A \textbf{contact circle bundle} is a circle bundle $E
  \overset{\pi}{\to} F$ together with a contact form $\alpha$ as
  described above.
\end{defn}

The standard contact $S^3$ is the simplest example of a contact circle
bundle.  Consider $S^3$ to be the unit sphere inside $\cc^2$.  Let
$\xi_0$ be the $2$-plane field defined by the complex tangencies to
$S^3$.  It is straightforward to check that $\xi_0$ is the kernel of
the form
\begin{equation*}
  \alpha_0 = \frac{1}{2} \sum_{j=1,2} x_j dy_j - y_j dx_j.
\end{equation*}
The Reeb field of $\alpha_0$ generates the Hopf fibration.

Note that the standard structure on $S^3$ is an example of this
construction.  For a more detailed introduction to contact geometry,
see \cite{aeb, confol, etnyre:intro}.

\subsection{Diagrams of Legendrian Knots}
\label{sec:diagrams}

This section explains how a Legendrian knot $L$ in a contact circle
bundle can be described combinatorially by three pieces of data: its
projection $\pi(L)$ to the base $F$, a choice of Reeb chord above
every double point of $\pi(L)$, and an integer vector with one
component for each region of $F \setminus \pi(L)$.  The last two
characterize the interaction of $L$ with the topology of $E$ in a
similar manner to Turaev's shadow link representation for topological
knots in circle bundles \cite{turaev:shadow}.

\begin{figure}[tbp]
  \centering{\input{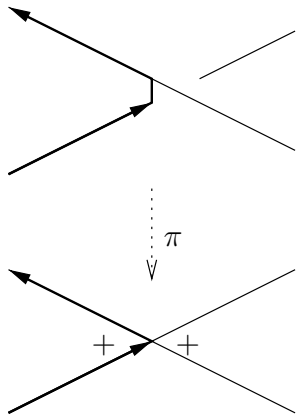}}
  \caption{The representation of a choice of chord above a double point.}
  \label{fig:decorated-diagram}
\end{figure}

Recall that the Reeb flow goes around the fibers of $E$.  Thus, as
mentioned in Section~\ref{sec:comb-intro}, above a double point of
$\pi(L)$, there exist countably many isolated Reeb chords that start
and end on different strands.  These come from concatenating a short
chord between the strands with a chord that winds several times around
the fiber.  Above each double point, choose one of the two Reeb chords
with distinct endpoints and with minimum length among all chords with
the same endpoints. Display this choice on $\pi(L)$ as follows: in
each quadrant near a double point, the orientation of $F$ gives a
counter-clockwise orientation to the edges of $\pi(L)$ that bound the
quadrant; see Figure~\ref{fig:decorated-diagram}. If the starting
point of the chosen Reeb chord lies on the incoming strand, then
decorate the quadrant with a $+$.  Otherwise, leave it blank.  After
following this algorithm, each double point has two opposing quadrants
decorated with a $+$.  Refer to this decorated diagram as $\pi(L)^+$.

The next step is to quantize the holonomy of $\alpha$ over the
piecewise smooth boundary of a surface immersed in $F$.  Of particular
interest in this section will be the regions of the graph $\pi(L)
\subset F$; later on, the definition of the invariant will require the
consideration of immersed disks.

\begin{figure}[tbp]
  \centering{\input{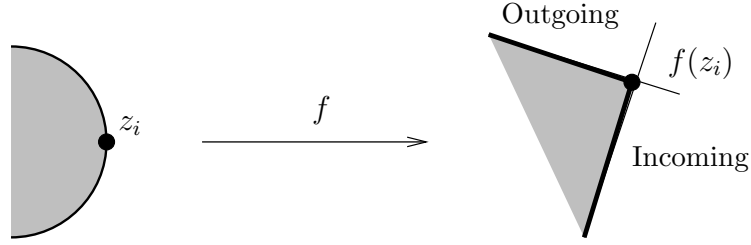}}
  \caption{The map $f$ sends a marked point $z_i$ to a double point.}
  \label{fig:corner-defect}
\end{figure}

Let $\Sigma$ be an oriented surface with non-empty boundary and let
$\{z_1, \ldots, z_m\}$ be marked points on $\partial \Sigma$.  Let $f:
\Sigma \to F$ be an orientation-preserving immersion on the interior
of $\Sigma$.  Assume that $f$ extends smoothly to $\partial \Sigma$
away from the marked points and that $f$ maps each marked point to a
double point of $\pi(L)$.  At a marked point, $f|_{\partial \Sigma}$
can turn from one strand above the double point to another, as in
Figure~\ref{fig:corner-defect}(a).  With respect to the orientation on
$\partial \Sigma$, say that $f|_{\partial \Sigma}$ passes from the
\textbf{incoming strand} to the \textbf{outgoing strand}.  

Above each $f(z_i)$, let $x_i$ Reeb chord with distinct endpoints.
Let $l(x_i)$ be the length of $x_i$; assume that $0 < l(x_i) < 2
\pi$.\footnote{This length may be defined either in terms of the
  Hermitian metric on $\mathcal{E}$ or, more topologically, by
  $\int_{x_i} \alpha$.} Define the number $\epsilon_i$ as follows:
\begin{equation}
  \epsilon_i = 
  \begin{cases}
    +1 & \quad x_i \text{ flows from the incoming to the outgoing
      strand,} \\
    -1 & \quad \text{otherwise.} 
  \end{cases}
\end{equation}

\begin{defn}
  \label{defn:defect} Let the immersion $f$ and the chords $x_1, \ldots,
  x_m$ be defined as above.  The \textbf{defect} of $f$ with respect
  to $x_1, \ldots, x_m$ is the integer defined by:
  \begin{equation}
    \label{eqn:defect}
    n(f; x_1, \ldots, x_m) = \frac{1}{2\pi} \left(\int_\Sigma f^*\Omega 
      + \sum_{j=1}^m \epsilon_j l(x_j) \right).
  \end{equation}
  Extend the defect linearly to formal chains of immersions.
\end{defn}

Here is another description of the defect.  Suppose that $\partial
\Sigma$ is connected.  Choose the component of $\partial \Sigma
\setminus \{z_1, \ldots, z_m\}$ that lies between $z_1$ and $z_2$ and
lift it to a Legendrian curve in $f^* E$.  Start lifting the next
component at a length $\epsilon_2 l(x_2)$ along the fiber away from
the end of the first lift.  Proceed this way until all components of
$\partial \Sigma \setminus \{z_1, \ldots, z_m\}$ have been lifted.
The defect is the winding number around the fiber of the curve defined
by the lifted path together with the Reeb chords; in other words, the
defect is a coarse measure of the holonomy of the connection $\alpha$
above the curve $f(\partial \Sigma)$.  This construction justifies the
assertion that the defect is an integer.

As a final note on the defect, suppose that $\epsilon_i = +1$ for $n$
of the $m$ chords.  Since the curvature of $\alpha$ is strictly
negative and the Reeb chords have length less than $2\pi$, the
following holds:

\begin{lem} \label{lem:defect-bound}
  The defect is bounded above by $n-1$.
\end{lem}

Now it is time to put all of this information together in a diagram.
The information about the defects of the regions of $F \setminus
\pi(L)$ is recorded in an integer vector $\vec{n}(L)$ that has one
component $n_i$ for each region $R_i$ of $\pi(L)$.  If the boundary of
$R_i$ contains the double points under the previously chosen Reeb
chords $x_{1}, \ldots, x_{m}$, then the components of $\vec{n}(L)$ are
defined by:
\begin{equation}
  n_i = n(R_i; x_{1}, \ldots, x_{m}).
\end{equation}
Note that the signs $\epsilon_j$ are positive if $R_i$ covers a
quadrant with a $+$ at $\pi(x_j)$ and negative otherwise.

\begin{defn}
  \label{defn:knot-diagram}
  Let $L$ be a Legendrian knot in a contact circle bundle $(E,
  \alpha)$ with chosen Reeb chords over each double point of $\pi(L)$.
  A pair $(\Gamma^+, \vec{n})$, where $\Gamma^+$ is a circle immersed
  in $F$ together with $+$ decorations at its double points and
  $\vec{n}$ is an integer vector with one component for each region of
  $F \setminus \Gamma$, is a \textbf{diagram} for $L$ if there exists
  an orientation-preserving diffeomorphism $\phi$ of $F$ that
  satisfies:
  \begin{enumerate}
  \item $\phi(\Gamma) = \pi(L)$,
  \item $\phi$ sends quadrants near double points of $\Gamma$ labeled
    with a $+$ to such quadrants in $\pi(L)$, and
  \item $\vec{n} = \vec{n}(L)$.
  \end{enumerate}
\end{defn}

The immersed circle $\Gamma$ will frequently be viewed as a $4$-valent
graph in $F$.  The components of the vector $\vec{n}$ are restricted
by the following easy consequence of the definition of the defect and
of the Chern-Weil theorem:

\begin{prop}
  \label{prop:defect-euler}
  Let $(\Gamma^+, \vec{n})$ be the diagram of a Legendrian knot in $E
  \to F$.  Then:
  \begin{equation*}
    \label{eqn:defect-euler}
    \sum_i n_i = e(E).
  \end{equation*}
\end{prop}

\begin{figure}[tbp]
  \centering{\input{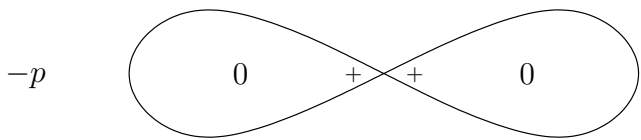}}
  \caption{The diagram of a local Legendrian unknot in $L(p,1)$.}
  \label{fig:trivial-diagram}
\end{figure}

\begin{ex}
  The projection of a small, null-homologous unknot in the lens space
  $L(p,1)$ with $tb = -1$ to $S^2$ is given by the Whitney immersion;
  see Figure~\ref{fig:trivial-diagram}.  Choose the shorter of the two
  Reeb chords at the crossing. The components of $\vec{n}$ associated
  with the two lobes of the immersion are both $0$.  By
  Proposition~\ref{prop:defect-euler}, the defect for the remaining
  region must be $-p$.
\end{ex}

\begin{figure}[tbp]
  \centering{\input{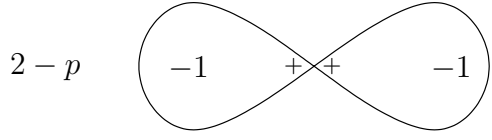}}
  \caption{The diagram of the Whitehead double of a fiber of $L(p,1)$.}
  \label{fig:simple-example}
\end{figure}

\begin{ex}
  The knot $K$ pictured in Figure~\ref{fig:simple-example} is
  topologically equivalent to the Whitehead double of a fiber of
  $L(p,1)$ over a point in $S^2$.
\end{ex}

As in Section 8 of Chekanov's paper \cite{chv}, it is possible to
combinatorially characterize Legendrian knot diagrams when $F \simeq
S^2$.  In order to specify which pairs $(\Gamma^+, \vec{n})$ are
diagrams for Legendrian knots, some notation is necessary.  Think of
$\Gamma$ as a $4$-valent graph with vertices $p_1, \ldots, p_r$ and
regions $U_1, \ldots, U_{r+2}$ in $F \setminus \Gamma$. Following
Chekanov, define two vector spaces over \rr: $Q_\Gamma$ is generated
by the vertices $p_j$, and $V_\Gamma$ is generated by the regions
$U_j$.  The Euler number $e$ of $E \to S^2$ generates a hyperplane in
$V_\Gamma$:
\begin{equation}
  V_{\Gamma, e} = \left\{ \sum_i c_i U_i \in V_\Gamma \; | \; \sum_i c_i =
  2 \pi e(E) \right\}.
\end{equation}
Let $U_+^1(p_j)$ and $U_+^2(p_j)$ be the two regions that are
decorated with a $+$ near $p_j$.  Let $U_-^1(p_j)$ and
$U_-^2(p_j)$ be the other two regions.  Define a linear map $\Psi:
Q_\Gamma \to V_\Gamma$ by:
\begin{equation}
  \Psi(p_j) = U_+^1(p_j) + U_+^2(p_j) - U_-^1(p_j) - U_-^2(p_j).
\end{equation}

Properly combined, these spaces and maps characterize diagrams of
Legendrian knots in $E$:

\begin{prop}
  \label{prop:lens-diagrams} Let $Q_\Gamma^{>0}$ be the vectors in
  $Q_\Gamma$ with all positive components, and let $V_{\Gamma,e}^{<0}$
  be the vectors in $V_{\Gamma,e}$ with all negative components.
  Then $(\Gamma^+, \vec{n})$ is the diagram of a Legendrian knot in $E
  \to S^2$ if and only if
  \begin{equation}
    \Psi(Q_\Gamma^{>0}) \cap (2\pi \vec{n} - V_{\Gamma,e}^{<0}) \neq
    \emptyset.
  \end{equation}
\end{prop}

\begin{proof}
  First suppose that $(\Gamma^+, \vec{n})$ is the diagram of a
  Legendrian knot $L$.  Let $\phi \in \text{Diff}^+(S^2)$ have the
  property that $\phi(\Gamma^+) = \pi(L)^+$.  Let $x_j$ be the chosen
  Reeb chord at the vertex $p_j$.  Define $q_L \in Q_\Gamma$ by:
  \begin{equation*}
    q_L = \sum_j l(x_j) p_j.
  \end{equation*}
  Clearly, $q_L \in Q_\Gamma^{>0}$.  Similarly, define an element
  $v_\phi \in V_\Gamma$ by:
  \begin{equation*}
    v_\phi = \sum_j k_j U_j,
  \end{equation*}
  where $k_j = \int_{U_j} \Omega$, the total curvature of $\alpha$
  over $U_j$.  By Proposition~\ref{prop:defect-euler} and the fact
  that the curvature is strictly negative, $v_\phi \in
  V_{\Gamma,e}^{<0}$.  Thus, by the definition of the defect,
  \begin{equation*}
    \Psi(q_L) = 2\pi \vec{n} - v_\phi.
  \end{equation*}
  The first half of the proposition follows.
  
  Conversely, suppose that there exist elements $q \in Q_\Gamma^{>0}$
  and $v \in V_{\Gamma,e}^{<0}$ that satisfy 
  \begin{equation} \label{eqn:q-v}
    \Psi(q) = 2\pi \vec{n} - v.
  \end{equation}
  On one hand, it is not hard to construct $\phi \in
  \text{Diff}^+(S^2)$ so that $v = v_\phi$.  It remains to show that
  Legendrian curve $L$ that projects to $\phi(\Gamma^+)$ is actually a
  knot.  To do this, it suffices to show that the total holonomy
  around $\phi(\Gamma)$ is an integral multiple of $2\pi$.
  
  To measure the holonomy, orient $\Gamma$ and let $\{C_1, \ldots, C_m
  \}$ be the Seifert circles of $\Gamma$.\footnote{This is where the
    assumption that $F \simeq S^2$ is necessary.  More generally, this
    proof works for knots with null-homologous projections in any base
    $F$.}  Each $C_i$ bounds a disk $D_i$ with corners; let $k_i$ be
  the total curvature over $D_i$ and let $n_i$ be the defect of $D_i$
  for the Reeb chords chosen by the ``$+$'' decorations on $\Gamma^+$.
  Note that if $D_i$ covers several regions, the defect and the total
  curvature of $D_i$ are the sum of the defects or curvatures,
  respectively, of those regions. The holonomy around a lift of
  $\phi(\Gamma)$ is the sum of the $k_i$, counted with sign: if
  $\Gamma$ goes around $\partial D_i$ counter-clockwise, then $D_i$
  contributes $k_i$ to the holonomy; otherwise, $D_i$ contributes
  $-k_i$.
  
  Equation (\ref{eqn:q-v}) implies that $k_i = 2\pi n_i - \sum_j
  \epsilon_j q_j$, where the $q_j$ are summed over all corners of
  $D_i$.  It is not hard to see that the $\epsilon_j q_j$ cancel out
  in the signed sum of the curvatures $k_i$.
  
  It follows that $\sum_i \pm k_i = \sum_i \pm 2\pi n_i$.  Since the
  $n_i$ are integers this proves the proposition.
\end{proof}

\section{A Combinatorial Definition of the Invariant}
\label{sec:defn}

There are three building blocks of the invariant promised in
Theorem~\ref{thm:overall}.  The first is a graded filtered algebra \alg\
associated to a knot diagram.  The second is a differential on the
algebra that, roughly speaking, counts immersed disks with boundary in
the knot diagram.  The last building block is an algebraic equivalence
relation on the DGAs whose equivalence classes are invariant under
Legendrian isotopy.  This section considers each building block in
turn.

\subsection{The Underlying Algebra \alg}
\label{sec:algebra}

\subsubsection{Generators for the Algebra}
\label{sec:generators}

\begin{figure}[tbp]
  \centerline{\input{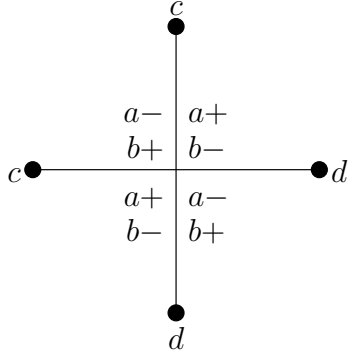}}
  \caption{Labels at a double point of $\Gamma$.}
  \label{fig:ab-labels}
\end{figure}

Let $(\Gamma^+, \vec{n})$ be a diagram for a Legendrian knot $L$.  The
definition of the generators of \alg\ requires a few additional
decorations on the diagram.  Number the double points from $1$ to $m$.
To the double point labeled $i$, associate two countable sets of
generators $\{a_i^k\}_{k=0,1,\ldots}$ and $\{b_i^k\}_{k=0,1,\ldots}$.
In addition, label the quadrants around each double point as in
Figure~\ref{fig:ab-labels}.\footnote{These labels play the same role
  as the $\pm$ labels in the Chekanov-Eliashberg differential.}  

\begin{figure}[tbp]
  \centerline{\input{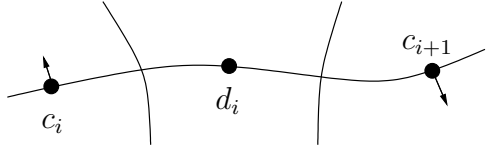}}
  \caption{Decorations at the edge points of $\Gamma$.}
  \label{fig:generators}
\end{figure}

On the interior of each edge $e$ of $\Gamma$, choose a point
$\hat{e}$.  Number these edge points from $1$ to $2m$ by traversing
$\Gamma$, starting from an arbitrary edge.  To the edge point
$\hat{e}_{2i-1}$, associate a countable set of generators
$\{c^k_i\}_{k=1, 2, \ldots}$.  Similarly, to the edge point
$\hat{e}_{2i}$, associate a countable set of generators
$\{d^k_i\}_{k=1, 2, \ldots}$.  Hereafter, the edge points
$\hat{e}_{2i-1}$ will be referred to as $c_i$, and the edge points
$\hat{e}_{2i}$ will be referred to as $d_i$.  Further, choose a
direction transverse to $\Gamma$ at each $c_i$; see
Figure~\ref{fig:generators}.

\begin{defn} \label{defn:special-pt}
  The \textbf{special points} of $\Gamma$ consist of the edge points
  $\hat{e}$ and the double points.
\end{defn}

\begin{defn} \label{defn:algebra}
  Let $\alg_{(\Gamma^+, \vec{n})}$ be the filtered unital
  based\footnote{A based algebra is an algebra together with a choice
    of generating set.}  associative algebra with coefficients in
  $\zz_2$ freely generated by the letters $a_i^k$, $b_i^k$,
  $c_i^{k'}$, and $d_i^{k'}$ for $i = 1, \ldots, m$, $k = 0, 1, 2,
  \ldots$, and $k' = 1, 2, 3, \ldots$.  The subscript $(\Gamma^+,
  \vec{n})$ will be dropped if there is no danger of confusion.
  
  The $k^{th}$ level of the filtration $F^0\alg \subset F^1\alg
  \subset \cdots$ is generated as a vector space by the monomials
  \begin{equation} \label{eqn:filtration}
    \left\{x_1^{j_1} \cdots x_m^{j_m}  \in \alg \vert \sum_i
      j_i \leq k
    \right\}.
  \end{equation}
\end{defn}

\begin{rem}
  Recall that there are two types of Reeb chords for a Legendrian knot
  $L \subset E$:
  \begin{enumerate}
  \item Those that start and end on different strands.  As mentioned
    in Section~\ref{sec:diagrams}, over a double point of $\pi(L)$,
    there are countably many such isolated Reeb chords, indexed by
    their ``winding number'' around the fiber.
  \item Those that begin and end at the same point on $L$.  There are
    countably many of them over any given point of $L$, indexed by
    their ``winding number'' around the fiber.
  \end{enumerate}
  
  Assume that $\Gamma = \pi(L)$.  The generators $a^k_i$ and $b^k_i$
  of $\alg_{(\Gamma^+, \vec{n})}$ correspond to the first type of Reeb
  chord that project to the double point $i$, wind around the fiber
  $k$ times, and whose endpoints lie on different strands of $L$ above
  the double point.  The chord corresponding to the $a_i^0$
  (respectively, $b_i^0$) generator is the one that passes from the
  incoming to the outgoing strand when the boundary of the quadrant
  labeled $a_i+$ (resp. $b_i+$)is given a counter-clockwise
  orientation.
  
  By definition, the \textbf{length} $l$ of these generators is
  \begin{equation} \begin{split}
      l(a_i^k) &= l(a_i^0) + 2 \pi k, \\
      l(b_i^k) &= l(b_i^0) + 2 \pi k.
    \end{split}
  \end{equation}
  
  The generators $c^k_i$ and $d^k_i$ represent the Reeb chords that
  start and end at the same point and project to the points labeled
  $c_i$ and $d_i$, respectively.  Note that 
  \begin{equation}
    l(c^k_i) = l(d^k_i) = 2 \pi k.
  \end{equation}
\end{rem}

Instead of working directly with the algebra $\alg$, it will sometimes
be useful to organize the generators into power series in
$\alg[[\nov]]$.  Of course, all equations involving these power series
should be understood as a sequence of equations relating the
coefficients of $\nov^k$ on the left hand side to those on the right
hand side for each $k$.  To formalize the power series idea, define:

\begin{defn}
  \label{defn:alg-pwr-series}
  A filtered unital based associative algebra \alg\ is
  \textbf{generated by power series} if the chosen generators of \alg\ 
  may be grouped into sequences $\{x_i^j\}_{j \geq k_i}$, where $x_i^j
  \in F^j \alg$.  For short, say that \alg\ is generated by the power
  series
  \begin{equation*}
    \ps{x}_i = \sum_{j=k_i}^\infty x_i^j \nov^j, \quad i=1, \ldots, m.
  \end{equation*}
\end{defn}

\begin{rem}
  So long as there is a universal bound on the number of generators in
  each level of the filtration, the algebra may be generated by power
  series.
\end{rem}

The generating power series of $\alg_{(\Gamma^+, \vec{n})}$ are:
\begin{align*}
  \ps{a}_i &= \sum_{k=0}^\infty a^k_i \nov^k &
  \ps{b}_i &= \sum_{k=0}^\infty b^k_i \nov^k \\
  \ps{c}_i &= \sum_{k=1}^\infty c^k_i \nov^k &
  \ps{d}_i &= \sum_{k=1}^\infty d^k_i \nov^k 
\end{align*}
Note that, since the power series $\ps{c}_i$ starts in degree one, the
element $1 + \ps{c}_i$ is invertible in $\alg[[\nov]]$.

\subsubsection{Grading}
\label{sec:grading}

\begin{figure}[tbp]
  \centerline{\input{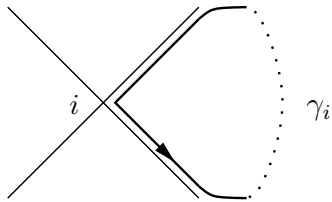}}
  \caption{Combinatorial construction of the capping path $\gamma_i$.}
  \label{fig:capping-path}
\end{figure}

A generator of $\alg_{(\Gamma^+, \vec{n})}$ has a grading if it is
associated to an edge point or to a double point that satisfies an
additional geometric condition.  Define a \textbf{capping path}
$\gamma_i$ for the double point $i$ to be one of the two oriented
paths in $\Gamma$ that begin at $i$, run along $\Gamma$ until they
first return to $i$, and induce a counter-clockwise orientation on the
quadrant they bound near $i$; see Figure~\ref{fig:capping-path}.

There are two quantities associated to a \emph{contractible} capping
path $\gamma_i$: its holonomy and its rotation.  To define the
holonomy, assume that the quadrant that $\gamma_i$ bounds has an $a+$
label (see Figure~\ref{fig:ab-labels}); the following constructions
are analogous if there is a $b+$ label.  In the former case, say that
$\gamma_i$ is a \textbf{capping path for $a_i$}; in the latter, say
that $\gamma_i$ is a \textbf{capping path for $b_i$}.  Using the
standard Seifert Circle algorithm, construct a chain of embedded
surfaces $\Sigma_i$ with $\partial \Sigma_i = \gamma_i$.  With this in
hand, define the holonomy $k_i$ of the capping path by:
\begin{equation}
  k_i = - n (\Sigma_i; a_i).
\end{equation}
Note that $k_i$ may be computed by summing up (with sign) the defects
listed in $\vec{n}$ of the embedded surfaces that comprise $\Sigma_i$.

To define the rotation, note that the contractibility assumption
implies that $\gamma_i \subset D^2 \subset F$ for some embedded disk
$D^2$ in $F$.  The rotation number of $\gamma_i$ with respect to a
trivialization of $TD^2$ may be computed as the fractional number of
counter-clockwise rotations made by the tangent vector to
$\gamma_i$.\footnote{In practice, this trivialization will be given by
  the presentation of the diagram.}  Call this number $r_D(\gamma_i)$.
If $\Sigma_i \subset D^2$, then the trivialization of $TD^2$ induces
one on $T\Sigma_i$.  In general, however, $\Sigma_i$ may not lie
inside the disk, and corrections to the trivialization must be made.
If $p$ is some point in $F$ outside the disk, then the
\textbf{rotation number} of $\gamma_i$ is:
\begin{equation*}
  r(\gamma_i) = r_D(\gamma_i) + \chi(F) (p \cdot \Sigma_i).
\end{equation*}

\begin{defn}
  \label{defn:grading} Without loss of generality, assume
  that the strands of $\Gamma$ meet orthogonally at $i$.  Suppose that
  $a_i$ has a contractible capping path.  The grading of
  $a_i^{-n(\gamma_i; a_i)}$ is a number in $\frac{1}{e(E)} \zz$
  defined by:
  \begin{equation}
    |a_i^{-n(\gamma_i; a_i)}| = 2r(\gamma_i) - \frac{1}{2}.
  \end{equation}
  Let $\mu_E = -\frac{\chi(F)}{e(E)}$. For arbitrary $k = 0, 1, 2,
  \ldots$, define:
  \begin{equation} \label{eqn:shifted-grading}
    \begin{split}
      |a_i^k| &= |a_i^{-n(\gamma_{a_i}; a_i)}| + \left(
        k+n(\gamma_{a_i}; a_i) \right) 2 \mu_E,\\
      |b_i^k| &= 2\mu_E (2k+1) - 1 - |a_i^k|.
    \end{split}
  \end{equation}
  Furthermore, the generators $c_j^k$ and $d_j^k$ are graded as
  follows:
  \begin{equation}
    \begin{split}
      |c_j^k| &= 2k\mu_E, \\
      |d_j^k| &= 2k\mu_E -1.
    \end{split}
  \end{equation}
\end{defn}

The grading is well-defined up to the choice of capping path.  If
$\Gamma$ is contractible in $F$, then there are many possible choices.
Any two capping paths will differ by a path that traverses the entire
knot an integral number of times.  If $n(L)$ is the total holonomy
around the knot, traversing the knot will add $2r(L)$ to the rotation
of $\gamma_i$ and $n(L)$ to $n(\gamma_i; a_i)$.  Thus:

\begin{prop}
  \label{prop:grading-ambig}
  For a knot whose projection is contractible, the grading on \alg\ in
  Definition~\ref{defn:grading} is well-defined modulo $2r(L) + 2\mu_E
  n(L)$.
\end{prop}

Note that this reduces to the ambiguity of Chekanov's grading in
\cite{chv} when $L$ is null-homotopic in $E$.

\begin{rem}
  This grading is a combinatorial translation of the Conley-Zehnder
  index of a Reeb chord.  Note that $\mu_E$ is the Maslov index of a
  fiber.  See \cite{egh} or \cite{ens} for a more geometric
  construction.
\end{rem}

\subsection{The Differential on \alg}
\label{sec:diffl}

The differential on \alg\ is the sum of two parts: an ``external''
differential and an ``internal'' differential.  This is analogous to
the differential in Morse-Bott theory, which splits into flowlines
\emph{between} critical submanifolds and flowlines \emph{inside} the
critical submanifolds.  Since the definitions of the internal and
external differentials are somewhat involved, a simple example is
calculated in Section~\ref{sec:simple-ex}.

\subsubsection{The External Differential}
\label{sec:ext-diffl}

The external differential \dfe\ is defined in much the same way as the
Chekanov-Eliashberg differential \cite{chv, ens}: it is a count of
certain immersed disks in $F$ with boundary in the diagram $\Gamma$.
There are, however, several new features that change the nature of the
boundary conditions for the disks.  First of all, there are now many
generators associated to each double point.  Secondly, there are new
generators associated to the edge points.  Lastly, the topology of the
bundle, as recorded by the defect, must be taken into account.

The first step is to define appropriate spaces of immersions of marked
disks.  In this context, a \textbf{marked disk} is a disk $D^2$
together with $m+1$ distinct marked points $z, w_1, \ldots, w_m$,
arranged counter-clockwise on $\partial D^2$.  Let $x$ be a label in
the set $\{a_i, b_i, d_i\}_{i=1 \ldots n}$ and let $y_1, \ldots, y_m$
be labels in the set $\{a_i, b_i, c_i\}_{i=1 \ldots n}$.

\begin{defn}
  \label{defn:ext-disks}
  The \textbf{space of immersed marked disks} $\Delta(x; y_1, \ldots,
  y_m)$ consists of orientation-preserving immersions of marked disks
  $f: D^2 \to F$ that satisfy:

  \begin{enumerate}
  \item The map $f$ sends $\partial D^2$ to $\Gamma$ and is smooth off
    of the marked points.  The restriction of $f$ to $\partial D^2$ is
    an immersion off of the marked points.

  \item The map $f$ sends $z$ to $x$ and $w_j$ to $y_j$.
    
  \item Near each marked point whose image is a double point, the
    image of $f$ in a neighborhood of the marked point covers but one
    quadrant of $\Gamma$.  If $x$ is a double point, then $f$ covers
    one of the two quadrants labeled $x+$; say that $f$ has a
    \textbf{positive corner} at $x$. If $y_j$ is a double point, then
    $f$ covers a quadrant labeled ${y_j}-$; say that $f$ has a
    \textbf{negative corner} at $y_j$.
    
  \item If $f|_{\partial D^2}$ crosses an edge point labeled $c_i$,
    then that $c_i$ must be the image of a marked point.
  \end{enumerate}
  
  Two immersions $f$ and $g$ are equivalent in $\Delta$ if there
  exists a smooth automorphism $\phi$ of the marked disk so that $f =
  g \circ \phi$.
\end{defn}

Using these spaces, the exterior differential of each generator of
\alg\ may be defined via power series formulae.

\begin{defn}
  \label{defn:ext-diffl}
  Let $x, y_1, \ldots, y_m$ be as above.  Let $n(f; x, y_1, \ldots,
  y_m)$ be the defect of $f$.  The
  \textbf{exterior differential} \dfe\ is defined by the formula:
  \begin{equation} \label{eqn:ext-diffl}
    \dfe \ps{x} = \sum_{(y_1, \ldots, y_m)}\sum_{f \in 
      \Delta(x; y_1, \ldots, y_m)} \tilde{\ps{y}}_1 \cdots
    \tilde{\ps{y}}_m \nov^{-n(f; x, y_1, \ldots, y_m)}, 
  \end{equation}
  where
  \begin{equation*}
    \tilde{\ps{y}} = \begin{cases}
      1 + \ps{y} & y = c_j \text{ and the transversal
        at } c_j \text{ points into } \image f; \\
      (1 + \ps{y})^{-1} & y = c_j \text{ and the transversal
        at } c_j \text{ points out of } \image f; \\
      \ps{y} & \text{otherwise.}
    \end{cases}
  \end{equation*}      
  If $m=0$, then the contribution of $\Delta(x)$ is $1$ for every
  distinct element of $\Delta(x)$.
  
  Extend \dfe\ to all of \alg\ via the Leibniz rule.
\end{defn}

\begin{rem}
  Strictly speaking, the defect of an element of $\Delta(x; y_1,
  \ldots, y_m)$ is a slight abuse of notation since $x$ or some of the
  $y_i$ may lie on an edge rather than at a double point.  Such chords
  are ignored in the calculation of the defect.
  
  The defect of an element $f$ of $\Delta(x; y_1, \ldots, y_m)$ may be
  easily computed combinatorially. Let $\tilde{n}$ be the defect of
  $f$ computed with respect to the choice of chords given by the
  diagram $\Gamma^+$, i.e. add up the defects in $\vec{n}$ of the
  regions of $\Gamma$ that lie in the image of $f$.  The defect of $f$
  with respect to the chords $x, y_1, \ldots, y_n$ comes from making
  the following adjustments to $\tilde{n}$:
  \begin{enumerate}
  \item If $x$ is a double point label and the $+$ decoration from
    $\Gamma^+$ does not lie in the same quadrant as the label $x+$,
    then add $1$ to $\tilde{n}$.
  \item For $i=1, \ldots, m$, if $y_i$ is a double point label and the
    $+$ decoration lies in the same quadrant as the label $y_i-$, then
    subtract $1$ from $\tilde{n}$.
  \end{enumerate}
  The result is $n(f; x, y_1, \ldots, y_m)$.
\end{rem}

\begin{prop}
  \label{prop:ext-diffl-good}
  The exterior differential is well-defined.  In particular, each
  coefficient in equation (\ref{eqn:ext-diffl}) is given by a finite
  sum.
\end{prop}

\begin{proof}
  Fix $k \geq 0$ and let $y_1^{k_1} \cdots y_m^{k_m}$ be a term in the
  differential of $x^k$ coming from a map $f \in \Delta(x;y_1, \ldots,
  y_m)$.  Equation (\ref{eqn:ext-diffl}) implies:
  \begin{equation} \label{eqn:boundary-condition}
    k+ n(f; x, y_1, \ldots, y_m) = \sum_{j=1}^m k_j,
  \end{equation}
  and hence there are finitely many choices for the $k_j$.
  
  It now suffices to prove that, up to equivalence, the sets
  $\Delta(x; y_1, \ldots, y_m)$ are finite.  Suppose that $\Gamma$
  coincides with $\pi(L)$. In this situation, the defining equation
  for the defect (see Definitions~\ref{defn:defect} and
  \ref{defn:knot-diagram}) is:
  \begin{equation} \label{eqn:defect-new}
    n(f; x, y_1, \ldots, y_m) = \frac{1}{2\pi} \left(\int_D f^* \Omega +
    l(x^0) - \sum_{j=1}^m l(y_j^0) \right).
  \end{equation}
  Here, it is understood that $l(c^0) = l(d^0) = 0$.  The terms
  $l(y_j^0)$ are understood to be $0$ if $y_j$ is a $c_j$ generator,
  and similarly for $l(x^0)$.  Using the fact that the curvature
  $\Omega$ is strictly negative and that $l(x^k) = l(x^0) + 2\pi k$,
  equations (\ref{eqn:boundary-condition}) and (\ref{eqn:defect-new})
  lead to:

  \begin{lem} \label{lem:bounded-curvature}
    \begin{equation} \label{eqn:bounded-curvature}
      0 > \int_D f^*\Omega \geq -l(x^k) + \sum_{j=1}^m l(y_j^k).
    \end{equation}
  \end{lem}
  
  In particular, the total curvature of the regions of $\pi(L)$
  covered by $f$ is bounded below by $-l(x^k)$.  This is a uniform
  bound for any map $f$ that contributes to $\dfe x^k$.  Thus, there
  exists a bound on the number of regions that $f$ covers, and hence
  there exist finitely many possible maps $f$.
\end{proof}

\begin{rem}
  Proposition~\ref{prop:ext-diffl-good} does \emph{not} imply that the
  sum in (\ref{eqn:ext-diffl}) is finite; it merely shows that it is
  well-defined as a power series.
\end{rem}

\subsubsection{The Internal Differential}
\label{sec:int-diffl}

\begin{figure}[tbp]
  \centerline{\input{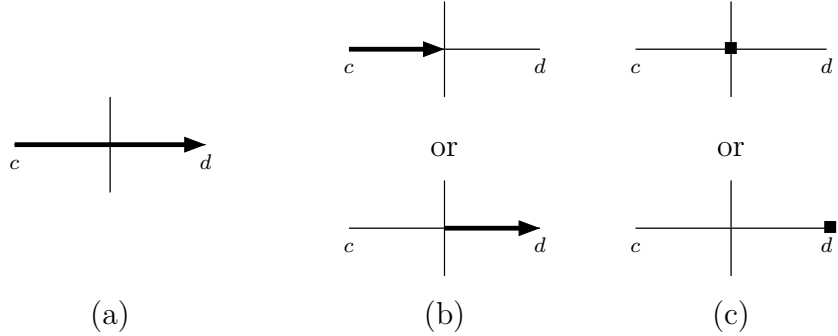}}
  \caption{(a) A flowline.  (b) Half flowlines.  (c)  Vertex flowlines.}
  \label{fig:internal-disks}
\end{figure}

To define the internal differential \dfi, it helps to think of
$\Gamma$ as having a ``Morse function'' with its maxima at the $c_i$
and its minima at the $d_i$.  Schematically, the internal differential
is defined by counting the flowlines, the ``half-flowlines'' between
the extrema and the crossings, and ``vertex flowlines'' that stay
fixed at the crossings and the minima.  See
Figure~\ref{fig:internal-disks}.

\begin{figure}[tbp]
  \centerline{\input{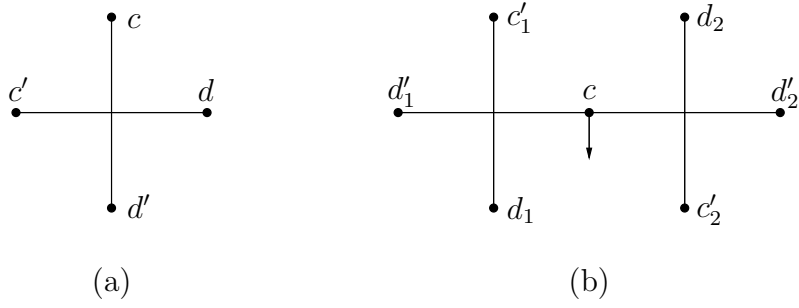}}
  \caption{(a) The configuration for $a$ and $b$ in 
    Definition~\ref{defn:internal-diffl}.  (b) The configuration for
    the differential of $c$ in Definition~\ref{defn:internal-diffl}.}
  \label{fig:internal-diffl}
\end{figure}

\begin{defn}
  \label{defn:internal-diffl}
  Let the labels near a double point be as in
  Figure~\ref{fig:internal-diffl}(a) and near a point labeled with a
  $c$ as in Figure~\ref{fig:internal-diffl}(b).  The \textbf{internal
    differential} \dfi\ is given by:
  \begin{description}
  \item[For \ps{a}] The internal differential comes from the two
    half-flowlines that start at the double point and end at the
    adjacent $d$ points:
    \begin{equation*}
      \dfi \ps{a} = \ps{a} \ps{d}' + \ps{d} \ps{a}.
    \end{equation*}
  \item[For \ps{b}] The first two terms come from the two
    half-flowlines that start at the double point and end at the
    adjacent $d$ points; the last comes from a vertex flowline at the
    double point:
    \begin{equation*}
      \dfi \ps{b} = \ps{b} \ps{d} + \ps{d}' \ps{b} + \ps{b} \ps{a}
      \ps{b} \nov.
    \end{equation*}
  \item[For \ps{c}] The first set of terms comes from the flowline and
    the half-flowline that depart $c$ to the left in
    Figure~\ref{fig:internal-diffl}(b).  The second set comes from a
    flowline and the half-flowline that go to the right:
    \begin{equation*}
      \dfi \ps{c} = (1+\ps{c}) (\ps{d}_1' + \ps{b}_1 \ps{a}_1 \nov)
      + (\ps{d}_2' + \ps{b}_2 \ps{a}_2 \nov) (1 + \ps{c}).
    \end{equation*}
    If the relative positions of the $c_1'$ and $d_1$ on the vertical
    strand at the left of Figure~\ref{fig:internal-diffl}(b) are
    reversed, then $\ps{b}_1 \ps{a}_1 \nov$ becomes $\ps{a}_1 \ps{b}_1
    \nov$, and similarly for the $c_2'$ and $d_2$.
  \item[For \ps{d}] The internal differential comes from a vertex
    flowline at $d$:
    \begin{equation*}
      \dfi \ps{d} = \ps{d} \ps{d}.
    \end{equation*}
  \end{description}
  Extend \dfi\ to all of \alg\ via the Leibniz rule.
\end{defn}

\subsubsection{The Full Differential}
\label{sec:full-diffl}

The internal and external differentials combine to form the full
differential on \alg:

\begin{thm}
  \label{thm:dga}
  The full differential $\df = \dfi + \dfe$ on \alg\ satisfies $\df
  \circ \df = 0$.  Further, \df\ preserves the filtration on \alg\ 
  and, when the grading is defined as in Section~\ref{sec:grading},
  has degree $-1$.
\end{thm}

The proof that $\df \circ \df = 0$ will be delayed until
Section~\ref{sec:d2}.

\begin{proof}[Proof that \df\ preserves filtration]
  Suppose that $\ps{y}_1 \cdots \ps{y_m} \nov^l$ is a term in $\df
  \ps{x}$. Expanding, this gives the terms $y_1^{j_1} \cdots
  y_m^{j_m}$ in $\df x^k$ for $k=l+\sum k_j$. So long as $l \geq 0$,
  \df\ preserves filtration.
  
  For \dfi, $l \geq 0$ by inspection.  For \dfe,
  Lemma~\ref{lem:defect-bound} and the fact that the disks that
  contribute to \dfe\ have but one positive corner imply that $l \geq
  0$.
\end{proof}

\begin{proof}[Proof that \df\ has degree $-1$]
  For the internal differential, the theorem follows from direct
  calculations.  For example, consider the internal differential of
  $a^k$:
  \begin{equation*}
    \dfi a^k = \sum_{j=0}^{k-1}a^j d^{k-j} + d^{k-j} a^j.
  \end{equation*}
  On the left hand side, the grading is:
  \begin{equation*} \label{eqn:internal-deg-calc}
    \begin{split}
      |a^k| &= |a^j| + 2 \mu_E (k-j) \\
      &= |a^j| + |d^{k-j}| + 1.
    \end{split}
  \end{equation*}
  This is one more than the grading of the terms on the right hand
  side.  The other computations using Definition~\ref{defn:grading}
  are similar.
  
  For the external differential, suppose that $y_1^{k_1} \cdots
  y_m^{k_m}$ is an arbitrary term in $\dfe x^k$ that comes from a disk
  $f$.  Let $\gamma_x$ be the capping path for $x$, and $\gamma_i$ be
  the capping path for $y_i$, as defined in Section~\ref{sec:grading}.
  Note that if the algorithm produces a capping path for $a_i$, then
  the same gradings are produced if $-\gamma{a_i}$ is used as a
  capping path for $b_i$.  For $c_i$ and $d_i$, define the capping
  path to be the constant path.  Define $r(\gamma_{c_i}) = 0$ and
  $n(\gamma_{c_i}; c_i) = 0$, and similarly for $d_i$.  
  
  Combine the capping paths of the generators $x, y_1, \ldots, y_m$
  with $f|_{\partial D}$ to form a loop:
  \begin{equation*}
    \gamma = \image f|_{\partial D} \sqcup -\gamma_x \sqcup \gamma_1 
    \sqcup \cdots \sqcup \gamma_m.
  \end{equation*}
  Since $\gamma$ is a smooth closed curve in $\Gamma$, either
  \begin{equation} \label{eqn:total-maslov}
    2r(\gamma) + 2 \mu_E n(\gamma) \equiv 0 \mod \mu(L)
  \end{equation}
  if $\Gamma$ is contractible in $F$, or $\gamma$ itself must be
  contractible in $F$, so:
  \begin{equation} \label{eqn:total-maslov-2}
    r(\gamma) = 0.
  \end{equation}  

  On the other hand, it is possible to calculate $2r(\gamma) + 2 \mu_E
  n(\gamma)$ piece by piece:
  \begin{equation} \label{eqn:r+n-gamma}
    \begin{split}
      r(\gamma) & = r(\image f|_{\partial D}) - r(\gamma_x) + 
      \sum_{i=1}^m r(\gamma_i),\\
      n(\gamma) & = n(f; x,y_1, \ldots, y_m) - n(\gamma_x; x) +
      \sum_{i=1}^m n(\gamma_i; y_i).
    \end{split}
  \end{equation}    
  Let $l$ be the number of convex corners of $f$.  It is not hard to
  see that $r(\partial f) = 1 - \frac{l}{4}$.  In the case that $x$ is
  not a $d$, combining (\ref{eqn:boundary-condition}) with the
  equations (\ref{eqn:r+n-gamma}) yields:
  \begin{equation*}
    \begin{split}
      2r(\gamma) + 2 \mu_E n(\gamma) &= 1 - \left( 2r(\gamma_x)
        - \frac{1}{2} + 2\mu_E (n(\gamma_x;x)+k) \right) \\
      &\quad +\sum_{y_i \neq c_j} \left(2r(\gamma_i) - \frac{1}{2} +
        2\mu_E (n(\gamma_i;y_i) +
        k_i) \right) + \sum_{y_i = c_j} 2\mu_E k_i \\
    &= 1 - |x^k| + \sum_{i=1}^m |y_i^{k_i}|.
    \end{split}
  \end{equation*}
  The case where $x=d_j$ is entirely similar.  The result now follows
  from (\ref{eqn:total-maslov}) or (\ref{eqn:total-maslov-2}).
\end{proof}

\subsection{A Simple Example}
\label{sec:simple-ex}

\begin{figure}[tbp]
  \centering{\input{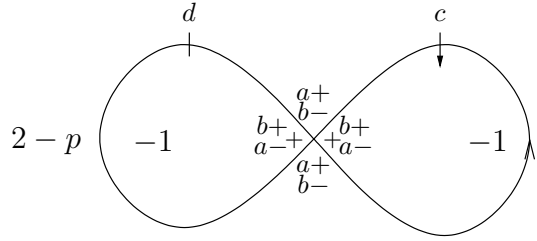}}
  \caption{The labeling of the diagram of the knot $L \subset L(p,1)$.}
  \label{fig:simple-ex}
\end{figure}

Let $L \subset L(p,1)$ be the knot introduced in the second example of
Section~\ref{sec:diagrams}.  The labeled diagram $(\Gamma^+,
\vec{n})$ for $L$ is pictured in Figure~\ref{fig:simple-ex}.  The
algebra for $(\Gamma^+, \vec{n})$ is generated by the power series
\ps{a}, \ps{b}, \ps{c}, and \ps{d}.

\subsubsection{The Grading}
\label{sec:simple-ex-grading}

The curve that runs counter-clockwise around the right-hand lobe of
$\Gamma$ is a capping path for $b$; call it $\gamma_b$.  The rotation
number of $\gamma_b$ is $\frac{3}{4}$ and $n(\gamma_b; b) =-1$, so
\begin{equation*}
  |b^1| = 2 \cdot \frac{3}{4} - \frac{1}{2} = 1.
\end{equation*}
By (\ref{eqn:shifted-grading}), 
\begin{align*}
  |a^k| &= \frac{4(k+2)}{p}-2 & |c^k| &= \frac{4k}{p} \\
  |b^k| &= 1 + \frac{4(k-1)}{p} & |d^k| &= \frac{4k}{p} -1.
\end{align*}

Since $L$ is null-homologous and its rotation number is $0$, the
grading is well-defined in $\frac{1}{p}\zz$.

\subsubsection{The Differential}
\label{sec:simple-ex-diffl}

Begin with the external differential of \ps{a}.  The disk $f$ shown in
Figure~\ref{fig:simple-disks}(a) contributes one of the two terms of
$\dfe \ps{a}$; exchanging the positions of the $a+$ and $b-$ gives the
other.  By the remark after Definition~\ref{defn:ext-diffl}, the
defect associated to the disk in Figure~\ref{fig:simple-disks}(a) is
$3-p$.  Thus,
\begin{equation} \label{simple-dfe}
  \dfe \ps{a} = \ps{b}(1+\ps{c})^{-1} \nov^{p-3} +
  (1+\ps{c})^{-1} \ps{b} \nov^{p-3}.
\end{equation}
Note that the $(1+\ps{c})^{-1}$ factor appears because the transverse
direction at $c$ points out of the disk.

\begin{figure}[tbp]
  \centering{\input{figs/simple-disks.pstex_t}}
  \caption{(a) The disk giving the term $\ps{b}(1+\ps{c})^{-1} 
    \nov^{p-3}$ in $\dfe \ps{a}$. (b) The disk giving the term
    $(1+\ps{c})\nov$ in $\dfe{b}$.}
  \label{fig:simple-disks}
\end{figure}

The disk shown in Figure~\ref{fig:simple-disks}(b) contributes one of
the two terms in $\dfe \ps{b}$.  The defect in the diagram is the same
as that of the disk, so this gives:
\begin{equation}
  \dfe \ps{b} = \nov + (1+\ps{c})\nov.
\end{equation}

The external differential for \ps{d} comes from disks whose boundaries
pass through $d$.  These are the left-hand lobe of the knot and the
disk in Figure~\ref{fig:simple-disks}(a).  The defect $n(\cdot;a)$ of
the left-hand lobe is one less than the defect in the diagram. The
defect of the outer disk agrees with the defect of the diagram.  As a
result,
\begin{equation}
  \dfe \ps{d} = \ps{a}\nov^{2} + \ps{b}(1+\ps{c})^{-1} \ps{b} \nov^{p-2}.
\end{equation}

The internal differential is calculated using
Definition~\ref{defn:internal-diffl} and
Figure~\ref{fig:internal-diffl}:
\begin{equation} \label{simple-dfi}
  \begin{split}
    \dfi \ps{a} &= \ps{a} \ps{d} + \ps{d} \ps{a}, \\
    \dfi \ps{b} &= \ps{b} \ps{d} + \ps{d} \ps{b} + \ps{b} \ps{a}
    \ps{b} \nov, \\
    \dfi \ps{c} &= (1+\ps{c})(\ps{d} + \ps{a} \ps{b} \nov) + (\ps{d} +
    \ps{b} \ps{a} \nov)(1+\ps{c}), \\
    \dfi \ps{d} &= \ps{d} \ps{d}. 
  \end{split} 
\end{equation}

Equations (\ref{simple-dfe}) through
(\ref{simple-dfi}), together with the fact that $\dfe \ps{c} = 0$,
determine the differential \df\ on \alg.

\subsection{Algebraic Notions}
\label{sec:isomorphisms}

The definition of the DGA in Sections~\ref{sec:algebra} and
\ref{sec:diffl} depends on the diagram $(\Gamma^+, \vec{n})$, the
choice of transverse direction at each $c$ point, and the position of
the $c$ and $d$ generators.  This section contains the definition of a
suitable equivalence relation on the algebras $(\alg, \df)$ which
captures how they change with the choices mentioned above.  

Before beginning the definitions, some notation is necessary.
Throughout this section, ``algebra'' means ``free, unital, graded,
based, associative, and filtered algebra''.  Any homomorphism of
algebras is assumed to respect these properties.  The filtration $F^k
\alg$ is always exhaustive and ascending:
\begin{equation}
  \cdots \subset F^k \alg \subset F^{k+1} \alg \subset \cdots
\end{equation}
Lastly, generators with filtration $k$ are denoted by $x^k$ when it is
convenient to emphasize the filtration.

The building block of the equivalence relation is the ``elementary
isomorphism'', which are strung together to form ``tame
isomorphisms''. Let \alg\ and $\bar{\alg}$ be algebras with a given
correspondence $x \leftrightarrow \bar{x}$ between their generators.
Choose a generator $y \in F^k\alg$ and an element $u \in
F^k\bar{\alg}$ so that $\bar{y}$ does not appear in $u$.  

\begin{defn} \label{defn:tame}
  The \textbf{elementary isomorphism} $\phi_y: \alg \to \bar{\alg}$ is
  defined on the generators by:
  \begin{equation}
    \phi_y(x) = \begin{cases}
      \bar{y} + u & \quad x = y, \\
      \bar{x} & \quad \text{otherwise}. 
    \end{cases}
  \end{equation}

  A \textbf{tame isomorphism} is the composition of a (possibly
  infinite) sequence of elementary isomorphisms $\cdots \circ
  \phi_{x_2} \circ \phi_{x_1}$ with the property that, for each level
  $k$ of the filtration, there exists $N_k \in \nn$ so that $x_j \not
  \in F^k\alg$ for all $j \geq N_k$.
\end{defn}

One convenient way to define a tame isomorphism is to use power
series.  Let \alg\ and $\bar{\alg}$ be algebras generated by power
series with a given correspondence $\ps{x} \leftrightarrow
\bar{\ps{x}}$ between their generating series.  Given a generating
series $\ps{y} \in \alg[[\nov]]$, let \ps{u} be a power series in
$\bar{\alg}[[\nov]]$ so that:
\begin{enumerate}
\item The coefficient $u^k$ of $\nov^k$ in \ps{u} lies in
  $F^k\bar{\alg}$ for all $k$;
\item $y^k$ never appears in $u^k$.
\end{enumerate} 

\begin{lem}
  \label{lem:tame-isom}
  The isomorphism $\boldsymbol{\psi}: \alg[[\nov]] \to
  \bar{\alg}[[\nov]]$ defined on the generating series by:
  \begin{equation}
    \boldsymbol{\psi}(\ps{x}) = \begin{cases}
      \bar{\ps{y}} + \ps{u} & \quad \ps{x} = \ps{y}, \\
      \bar{\ps{x}} & \quad \text{otherwise}
    \end{cases}
  \end{equation}
  gives a tame isomorphism $\psi: \alg \to \bar{\alg}$.
\end{lem}

\begin{proof}
  Let $\psi_k$ be the map that sends $y^k$ to $\bar{y}^k + u^k$ and is
  the identity on all other generators.  The first condition above
  implies that $\psi_k$ is a filtered map.  The second condition
  guarantees that $y^k$ does not appear in $u^k$.  The map $\psi$ is
  the composition of elementary isomorphisms $\cdots \circ \psi_2
  \circ \psi_1 \circ \psi_0$.
\end{proof}

\begin{rem}
  Suppose that there is an ordering on the generators of \alg\ that is
  compatible with the filtration $F^k\alg$; for instance, ordering the
  generators by length in the algebra of a knot.  Then, by the same
  argument as in the proof above, any map $\psi$ defined on the
  generators by the formula:
  \begin{equation*}
    \psi(x) = x + u_x,
  \end{equation*}
  where $u_x$ is a term containing only generators that lie strictly
  below $x$ in the ordering, is a tame isomorphism.
\end{rem}

Let $\mathcal{E}$ be an algebra generated by the two power series
$\boldsymbol{\alpha}$ and $\boldsymbol{\beta}$.  Suppose that
$|\alpha^k| +1 = |\beta^k|$, and define
\begin{equation}
  \begin{split}
    \df_{\mathcal{E}} \boldsymbol{\beta} &= \boldsymbol{\alpha}, \\
    \df_{\mathcal{E}} \boldsymbol{\alpha} &= 0.
  \end{split}
\end{equation}

\begin{defn} \label{defn:stable-tame}
  The \textbf{stabilization} $S(\alg, \df)$ of a DGA $(\alg, \df)$ is
  the DGA generated by the union of the generators of $\alg$ and those
  of $\mathcal{E}$.  $S(\alg, \df)$ inherits its grading, filtration,
  and differential from $(\alg, \df)$ and $(\mathcal{E},
  \df_{\mathcal{E}})$. 

  Two algebras \alg\ and $\bar{\alg}$ are \textbf{stable tame
  isomorphic} if there exist stabilizations $S_1, \ldots, S_m$ and
  $\bar{S_1}, \ldots, \bar{S_n}$ and a tame isomorphism:
  \begin{equation*}
    \psi: S_1 \left( \cdots S_m(\alg) \cdots \right) \to
    \bar{S_1} \left( \cdots \bar{S_n}(\bar{\alg}) \cdots \right).
  \end{equation*}
  
  Two DGAs $(\alg, \df)$ and $(\bar{\alg}, \bar{\df})$ are
  \textbf{stable tame isomorphic} if \alg\ and $\bar{\alg}$ are stable
  tame isomorphic via a chain map intertwining \df\ and $\bar{\df}$.
\end{defn}

Stable tame isomorphism is the appropriate equivalence relation on the
DGAs, as evidenced by the following theorem:

\begin{thm} \label{thm:invariance}
  Let $L$ be a Legendrian knot in a contact circle bundle $E \to F$.
  Let $(\alg, \df)$ be the DGA for a decorated and labeled diagram of
  $L$.  The stable tame isomorphism type of \alg\ is:
  \begin{enumerate}
  \item Independent of the choice of transverse direction at each
    $c_j$;
  \item Independent of whether a given edge point has a $c_i$ or $d_i$ label.
  \item Invariant under Legendrian isotopy of $L$.
  \end{enumerate}
\end{thm}

Notice that $(F^0 \alg, \df)$ is a sub-DGA of the DGA for a knot
diagram.  Since all of the stable tame isomorphisms involved in
Theorem~\ref{thm:invariance} preserve the filtration, the following
holds:

\begin{cor}
  The stable tame isomorphism type of the sub-DGA  $(F^0 \alg, \df)$
  is an invariant of $L$ as in Theorem~\ref{thm:invariance}.
\end{cor}

\begin{rem}
  The equivalence relation defined by stable tame isomorphism is
  stronger than the one defined by quasi-isomorphism; in other words,
  the homology of $(\alg, \df)$ is also an invariant.  To prove this,
  it suffices to show that the inclusion $i: \alg \to S(\alg)$ is a
  chain equivalence on the level of vector spaces.  Let $\tau: S(\alg)
  \to \alg$ be the natural projection.  Clearly, $\tau \circ i =
  id_\alg$.  On the other hand, $i \circ \tau$ is chain homotopic to
  the identity on $S(\alg)$.  To see this, define a vector space map
  $H: \alg \to S(\alg)$ by:
  \begin{equation} \label{eqn:H-defn}
    H(w) = \begin{cases}
      x \beta^k y &\quad w=x \alpha^k y \text{ and } x \in \alg, \\
      0 &\quad \text{otherwise.}
    \end{cases}
  \end{equation}
  It is straightforward to check that $H$ is the required chain
  homotopy, i.e. that it satisfies:
  \begin{equation} \label{eqn:H-chain}
    \tau \circ i + Id_\alg = H \circ \df + \df \circ H.
  \end{equation}
\end{rem}

\section{Applications}
\label{sec:apps}

\subsection{An Example with Nontrivial Degree in the Fiber}
\label{sec:fiber-ex}

\subsubsection{The Knots $L_1$ and $L_2$}
\label{sec:fiber-knots}

\begin{figure}[tbp]
  \centerline{\input{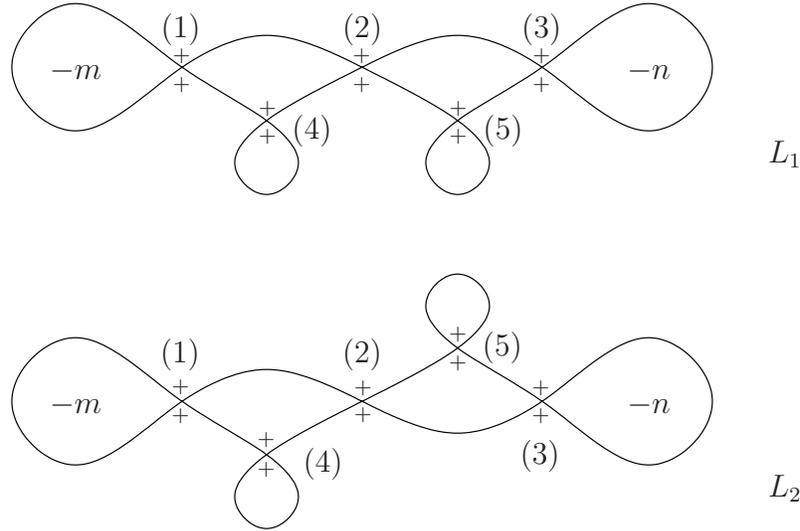}}
  \caption{The knots $L_1$ and $L_2$ inside $L(p,1)$ for $p \geq
    m+n+2$.  The outer regions have defects of $m+n-p$ and the
    unlabeled inner regions all have defect equal to $0$.  The corners
    of the bounded regions of the diagram are labeled by $b+/a-$.}
  \label{fig:fiber-ex}
\end{figure}

The two knots $L_1$ and $L_2$ in $L(p,1)$ pictured in
Figure~\ref{fig:fiber-ex} furnish the examples that prove the first
part of Proposition~\ref{prop:overall-ex}.  Assume that $p \geq n+m+2$
and that $n,m \geq 2$.  A simple application of
Proposition~\ref{prop:lens-diagrams} shows that the diagrams in
Figure~\ref{fig:fiber-ex} are indeed projections of Legendrian knots.

\begin{prop}
  \label{prop:fiber-ex}
  The Legendrian knots $L_1$ and $L_2$ are smoothly isotopic, but not
  Legendrian isotopic. If $n \neq m$, then $L_1$ represents a
  nontrivial class in $H_1(L(p,1))$.
\end{prop}

By undoing the small loops at crossings $4$ and $5$, it is clear that
$L_1$ and $L_2$ are topologically isotopic.  After introducing some
algebraic machinery in Section~\ref{sec:char-alg}, the second part of
the proposition will be proven in Section~\ref{sec:fiber-proof}.

\subsubsection{The Characteristic Algebra}
\label{sec:char-alg}

Working directly with the DGAs $(\alg_1, \df_1)$ and $(\alg_2, \df_2)$
of the knots $L_1$ and $L_2$ is a difficult task.  To facilitate
computations for Legendrian knots in $\rr^3$, Ng introduced the
\textbf{characteristic algebra} \cite{lenny:computable}.  A direct
adaptation of his definition to the circle bundle situation would
still yield a complicated object.  Since $F^0 \alg$ is an invariant
sub-DGA of $(\alg, \df)$, however, it is possible to restrict the
definition of the characteristic algebra to the lowest-energy
generators.

\begin{defn}
  \label{defn:char-alg} 
  Let $I$ be the two-sided ideal in $F^0 \alg$ generated by the set
  \begin{equation*}
    \{\df x_i^0 \; \vert \; x_i^0 \text{ a generator in } F^0 \alg\}.
  \end{equation*}
  The low-energy \textbf{characteristic algebra} \ca{} of the filtered
  DGA $(\alg, \df)$ is $(F^0 \alg) / I$.
\end{defn}

Suppose that $(\alg_1, \df_1)$ and $(\alg_2, \df_2)$ are related by an
elementary isomorphism $\phi$ that sends $x$ to $x'+v$.  Since
\begin{equation}
  \phi(\df_1 x) = \df_2 \phi (x) = \df_2 (x'+v),
\end{equation}
the Leibniz rule shows that $\phi$ identifies $I_1$ and $I_2$.  Thus,
say that \ca{1} and \ca{2} are \textbf{tame isomorphic} if $\alg_1$
and $\alg_2$ are tame isomorphic as algebras (after possibly adding
``trivial'' generators to both $\alg_i$ and $I_i$) so that the tame
isomorphism identifies $I_1$ and $I_2$.

Adding a stabilization $\mathcal{E}$ to \alg\ adds two generators
$e_1$ and $e_2$ and the relation $e_2=0$ to \ca{}.  Thus, under
stabilization, \ca{} changes by the addition of one generator and no
further relations.

Put together, the equivalence relation of stable tame isomorphism on
\alg\ translates to:

\begin{defn}[Ng \cite{lenny:computable}]
  \label{defn:equiv-char-alg}
  Two characteristic algebras \ca{1} and \ca{2} are
  \textbf{equivalent} if they are tame isomorphic after adding a
  finite number of generators (with no additional relations) to each.
\end{defn}

Theorem~\ref{thm:dga} and the discussion above prove:

\begin{prop}[Ng \cite{lenny:computable}]
  \label{prop:char-alg}
  If two Legendrian knots are Legendrian isotopic, then their
  low-energy characteristic algebras are equivalent.
\end{prop}

\begin{rem}
  It is also possible to adapt Chekanov's graded linearization method
  (see \cite{chv, efm, fuchs:augmentations}) to this setting, but it
  cannot distinguish $L_1$ and $L_2$.  There are other situations in
  which it is useful, however.
\end{rem}

\subsubsection{Proof of Proposition~\ref{prop:fiber-ex}}
\label{sec:fiber-proof}

In order to calculate \ca{}, it is necessary to know the differentials
only of $a_i^0$ and $b_i^0$; the other generators of \alg\ never come
into play. With this in mind, drop the superscript $0$ from the
notation for the remainder of this section.  The internal differential
does not contribute, and the only disks that contribute to the
external differential are those whose defect is $0$ and do not contain
a $c^k$ term.  In the case at hand, the only such disks come from the
four regions in the center of the diagram.  In particular, note that
the regions at the ends have defect $-1$ for the $b_1+$ or $b_3+$
chords, so they do not contribute to $F^0 \alg$.

It follows that $\df_i a_j = 0$ for $i=1,2$ and $j=1, \ldots, 5$.
Further,
\begin{align*}
  \df_1 b_1 &= a_4 a_2 & \df_2 b_1 &= a_4 a_2 \\
  \df_1 b_2 &= a_1 a_4 + a_5 a_3 & \df_2 b_2 &= a_1 a_4 + a_3 a_5 \\
  \df_1 b_3 &= a_2 a_5 & \df_2 b_3 &= a_5 a_2 \\
  \df_1 b_4 &= 1+ a_2 a_1 & \df_2 b_4 &= 1+ a_2 a_1 \\
  \df_1 b_5 &= 1+ a_3 a_2 & \df_2 b_5 &= 1+ a_2 a_3 \\
\end{align*}

Some simple manipulations of both characteristic algebras show that
$a_4$ and $a_5$ are trivial in both \ca{1} and \ca{2}.  Further,
$a_1=a_3$ in \ca{1}.  Thus, up to equivalence,
\begin{align*}
  \ca{1} &= \zz_2 \langle a_1, a_2 \rangle /(1 = a_1 a_2 = a_2 a_1), \\
  \ca{2} &= \zz_2 \langle a_1, a_2, a_3 \rangle /(1 = a_2 a_1 = a_2
  a_3).
\end{align*}
These algebras cannot be equivalent, since all elements in \ca{1} are
invertible from both sides, while the element $a_2 \in \ca{2}$ is only
invertible from the right.  This finishes the proof of
Proposition~\ref{prop:fiber-ex}.

\subsection{An Example with Nontrivial Homology in the Base}
\label{sec:base1-ex}

The second example presents two knots $L_3$ and $L_4$ whose
projections to the base $F$ represent nontrivial homology classes in
$H_1(F)$, and hence provides the necessary examples to prove the
second half of Proposition~\ref{prop:overall-ex}.  Let $\gamma$ be a
simple closed curve that represents a nontrivial homology class in
$H_1(F)$.  Let $\alpha$ be an arbitrary contact form on $E$.  Lift
$\gamma$ to a Legendrian curve $\tilde{\gamma}$ in $E$.  Choose a
closed form $\beta \in \Omega^1(F)$ whose period around $\gamma$ is
equal to the negation of the holonomy of $\tilde{\gamma}$.  For
example, take $\beta$ to be an appropriately scaled Poincar\'e dual to
$\gamma$.  As a result, $\gamma$ is the projection of a closed
Legendrian curve for the contact form $\alpha + \pi^*\beta$.

\begin{figure}[tbp]
  \centerline{\input{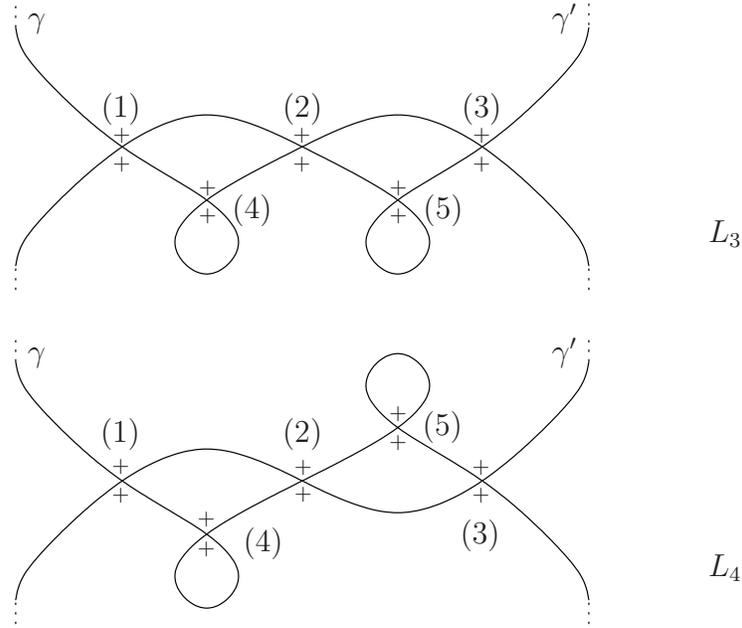}}
  \caption{The knots $L_3$ and $L_4$ in $E \to F$. The region between
    $\gamma$ and $\gamma'$ and outside the the four central regions
    has defect $0$, and the outer region has defect $p \leq -1$.}
  \label{fig:base1-ex}
\end{figure}

To form $L_3$ and $L_4$, let $\gamma'$ be a small translation of
$\gamma$ in $F$.  Remove small segments of $\gamma$ and $\gamma'$ and
glue in the middle portions of $L_1$ and $L_2$.  See
Figure~\ref{fig:base1-ex}.

\begin{prop}
  \label{prop:base1-ex}
  The Legendrian knots $L_3$ and $L_4$ are topologically isotopic, but
  not Legendrian isotopic.
\end{prop}

As with $L_1$ and $L_2$, the only disks that contribute to the
exterior differential on $F^0 \alg$ lie in the four central regions.
Thus, the proof of Proposition~\ref{prop:base1-ex} is the same as that
of Proposition~\ref{prop:fiber-ex}.

\section{Proof that \alg\ is a DGA}
\label{sec:d2}

The goal of this chapter is to prove that $\df \circ \df = 0$, which
will complete the proof of Theorem~\ref{thm:dga}.  In outline, the
proof is similar to those in \cite{chv, ens}, which, in turn, are
combinatorial realizations of the standard proofs in
Morse-Witten-Floer theory (see \cite{austin-braam, schwarz}, for
example).

\subsection{Outline of the Proof}
\label{sec:outline}

It is convenient to prove the theorem on the level of power series.
Suppose that $\ps{x}$ is a generating power series in
$\alg_{(\Gamma^+, \vec{n})}$.  Let $\ps{y}_1 \cdots \ps{y}_n$
represent a disk or flowline $f \in \Delta^U(x;y_1, \ldots y_n)$, and
let $\ps{v}_1 \cdots \ps{v}_m$ represent a disk or flowline $g \in
\Delta^U(y_i; v_1, \ldots, v_m)$.  Thus,
\begin{equation} \label{eqn:d2-term}
  \ps{y}_1 \cdots \ps{y}_{i-1} \ps{v}_1 \cdots \ps{v}_m \ps{y}_{i+1}
  \cdots \ps{y}_n \nov^{-n(f)-n(g)}
\end{equation}
is a term in $\df (\df \ps{x})$.  The pair $(f,g)$ is called a
\textbf{broken disk at \ps{x}}.  Two broken disks are equivalent if
their component maps are.

The proof that $\df(\df \ps{x}) =0$ proceeds by organizing the terms
in $\df (\df \ps{x})$ --- in other words, the broken disks at \ps{x}
--- into canceling pairs.  The two maps that make up a broken disk
can be glued together at $y_i$ to form a new immersed disk.  The
resulting \textbf{obtuse disk} satisfies all of the requirements of
the definition of $\Delta^U(x;y_1, \ldots, v_1, \ldots, v_m, \ldots,
y_n)$ except that one corner covers three quadrants rather than one.
Conversely, Section~\ref{sec:degeneration} shows that every obtuse
disk splits into a broken disk in exactly two ways.  Since the two
broken disks comes from the same obtuse disk, they contribute
identical terms to $\df(\df \ps{x})$.  This will complete the proof.

\subsection{An Equivalent Definition of \df}
\label{sec:thick-flowlines}

A more unified treatment of broken and obtuse disks may be achieved by
altering the definition of the full differential.  Instead of
separating \df\ into internal and external pieces, all terms in the
unified definition come from immersed disks in a modification of the
diagram $\Gamma$. The motivation for this definition is the
perturbation of the ``Morse-Bott'' contact form as described in
\cite{bourgeois:mb}.

The unified definition requires some combinatorial modifications to
the diagram $\Gamma$.  First, fix the following terminology.  A
\textbf{full lattice} is the set of lines $\{x=n\}_{n \in \zz} \cup
\{y=n\}_{n \in \zz}$ in $\rr^2$.  The intersections of these lines in
the lattice are called \textbf{lattice points}.  A
\textbf{half-lattice} is the portion of a full lattice that lies below
the diagonal line $y = x + \frac{1}{2}$.  This line is called the
\textbf{edge} of the half-lattice.  There are two different concepts
of distance on a half-lattice.  Let $(m,n)$ be a lattice point.  The
\textbf{distance $\delta$ to the edge} is given by either:
\begin{align*}
  \delta_{ab} (m,n) &= m-n, \quad \text{or} \\
  \delta_{cd} (m,n) &= m-n+1.
\end{align*}
There is one last piece of terminology necessary: a
\textbf{translation of a full lattice} is a linear map $\tau$ that
sends lattice points to lattice points.  A \textbf{translation of a
  half-lattice} is the restriction of a translation of a full lattice
to the points in the half lattice whose images lie in the
half-lattice.

\begin{figure}[tbp]
  \centerline{\input{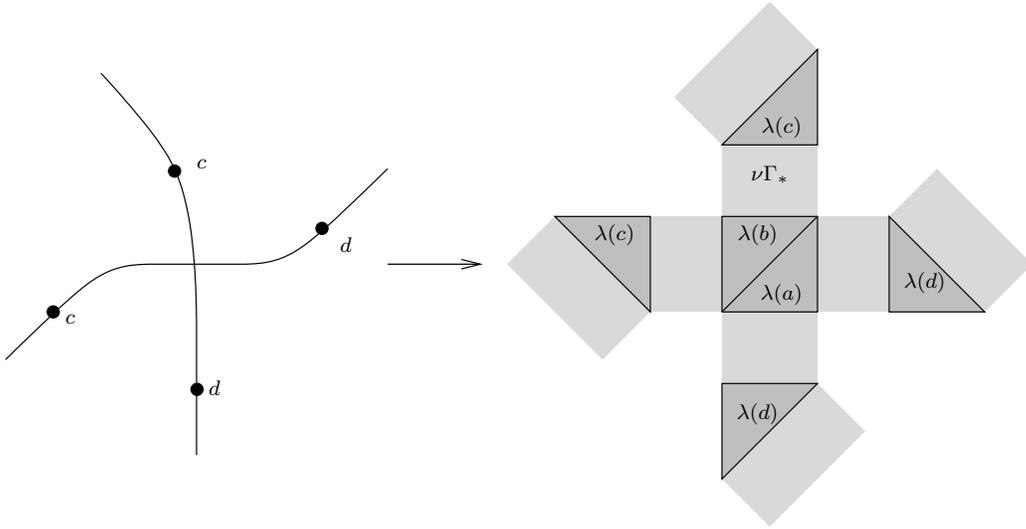}}
  \caption{A schematic representation of the modified diagram.}
  \label{fig:modified-diagram}
\end{figure}

\begin{figure}[tbp]
  \centerline{\input{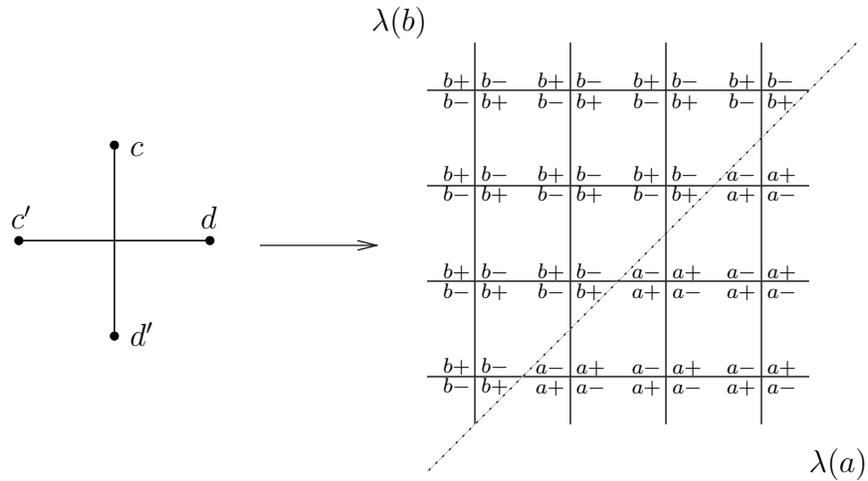}}
  \caption{A piece of the half-lattices $\lambda(a)$ and $\lambda(b)$ 
    at a double point.}
  \label{fig:ab-lattices}
\end{figure}

\begin{figure}[tbp]
  \centerline{\input{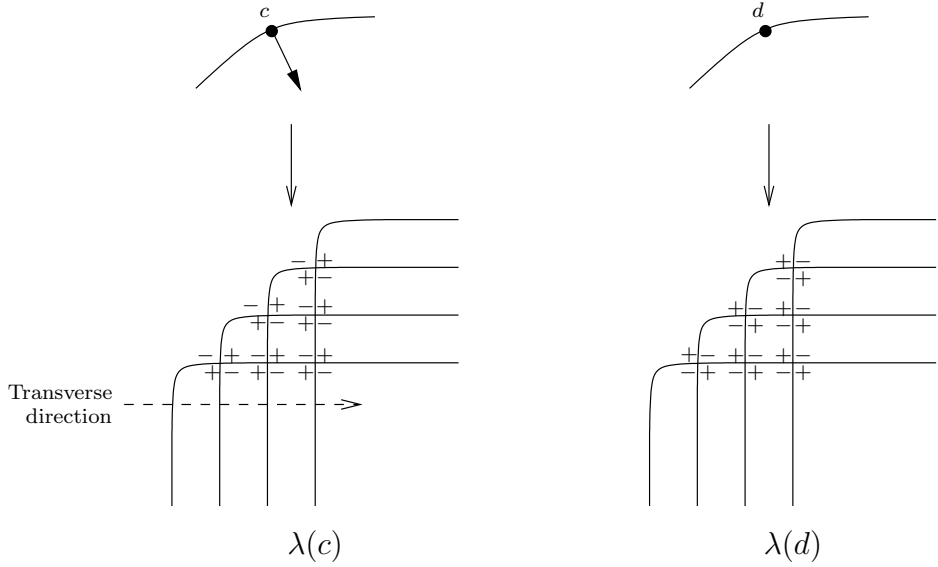}}
  \caption{(a) A piece of the half-lattice $\lambda(c)$; (b) A piece
    of the half-lattice $\lambda(d)$.}
  \label{fig:cd-lattices}
\end{figure}

The idea behind the modifications to $\Gamma$ is to replace every
special point by an embedding of a full or half lattice; see
Figure~\ref{fig:modified-diagram}.  Fix a small neighborhood $U$ of
$\Gamma$.  As pictured in Figure~\ref{fig:ab-lattices}, replace a
double point of $\Gamma$ by a full lattice, split into half-lattices
$\lambda(a)$ and $\lambda(b)$ that are matched up along their edges.
Equip $\lambda(a)$ and $\lambda(b)$ with the distance function
$\delta_{ab}$.  Label the quadrants around each double point in
$\lambda(a)$ (resp., $\lambda(b)$) with $a+$ and $a-$ (resp., $b+$ and
$b-$) as in Figure~\ref{fig:ab-lattices}.  Take note of the position
of the adjacent $c$ and $d$ special points in
Figure~\ref{fig:ab-lattices}.  The special points labeled $c$ and $d$
become embedded half-lattices $\lambda(c)$ and $\lambda(d)$,
respectively, that are capped off along the lattice edge as pictured
in Figure~\ref{fig:cd-lattices}.  Both are equipped with the distance
function $\delta_{cd}$.  Let $\Gamma_*$ denote $\Gamma \setminus
\{\text{special points}\}$.  Outside the lattices, $U$ has the
structure of the normal bundle $\nu \Gamma_*$.

Roughly speaking, the differential counts immersed disks whose
boundaries lie in $U$ and that have specified behavior in the
lattices.  More precisely, let $(D^2; z, w_1, \ldots, w_m)$ be a disk
with marked points ordered counter-clockwise on $\partial D^2$, as in
Section~\ref{sec:ext-diffl}.  Let $x, y_1, \ldots, y_m$ be special
points of $\Gamma$.  Define a smooth map $\rho: F \to F$ that retracts
$U$ onto $\Gamma$ and collapses the half-lattices $\lambda$ to their
associated special points.

\begin{defn}
  \label{defn:thick-flowlines}
  The \textbf{space of $U$-immersed marked disks} $\Delta^{U}(x; y_1,
  \ldots, y_m)$ consists of immersions $f: D^2 \to F$ of marked disks
  that satisfy the following:

  \begin{enumerate}
  \item The map $f$ sends $z$ to a lattice point of $\lambda(x)$ and
    $w_j$ to a lattice point of $\lambda(y_j)$.
    
  \item For every marked point $z$ (respectively, $w_i$), there exists
    a neighborhood $N_z$ (resp., $N_i$) such that the map $f$ sends
    $\partial D^2 \cap N_*$ to the lines of the corresponding
    half-lattice $\lambda$.  Elsewhere, $f$ sends $\partial D$ either
    to a line of a half-lattice or to a section of $\nu \Gamma_*$.
    Further, $f|_{\partial D \setminus \{z, w_1, \ldots, w_m\}}$ is
    smooth.  Notice that $f(\partial D) \subset U$.
    
  \item Let the neighborhoods $N_*$ be as above.  Near $z$ and $w_i$,
    the image of $f|_{N_*}$ covers but one quadrant of $\Gamma$.  At
    the lattice point $f(z)$, $f|_{N_z}$ covers one of the two
    quadrants labeled $x+$ and no other quadrants; say that $f$ has a
    \textbf{positive corner} at $x$. At the lattice point $f(w_j)$,
    $f|_{N_j}$ covers one of the two quadrants labeled $y_j-$ and no
    others; say that $f$ has a \textbf{negative corner} at $y_j$.
    
  \item The map satisfies:
    \begin{equation*}
      \delta(f(z)) = \sum_j \delta(f(w_j)) - n(\rho \circ f; x, y_1, 
      \ldots, y_m).
    \end{equation*}
    Here, $\delta$ is either $\delta_{ab}$ or $\delta_{cd}$, as appropriate.
    This is the \textbf{energy condition}.
  \end{enumerate}
  
  Finally, say that two immersions $f$ and $g$ are \textbf{equivalent}
  in $\Delta^{U}(x; y_1, \ldots, y_n)$ if there exists a smooth
  automorphism $\phi$ of the marked disk, neighborhoods $N_z, N_1,
  \ldots, N_m$ in $D^2$,\footnote{These neighborhoods may be different
    from those specified in the first condition.} and translations
  $\tau_x, \tau_{y_1}, \ldots, \tau_{y_m}$ of the half-lattices
  $\lambda(x), \lambda(y_1), \ldots, \lambda(y_m)$ (extended to
  translations of the full lattice at the double points) so that:
  \begin{equation}
    \label{eq:equiv}
    \begin{split}
      f &= \tau_x \circ g \circ \phi \quad \text{in } N_z, \\
      f &= \tau_{y_i} \circ g \circ \phi \quad \text{in } N_i,\\
      \rho \circ f &= \rho \circ g \circ \phi.
    \end{split}
  \end{equation}
\end{defn}

\begin{rem}
  If $\image f \subset U$, then $\rho \circ f$ is no longer an
  immersion.  Define the defect of such an $f$ by summing the lengths
  of the chords that lie over double points:
  \begin{equation} \label{eqn:flowline-defect}
    2\pi n(f; x, y_1, \ldots, y_m) = \begin{cases}
      l(x) - \sum_{y_i \neq c_j} l(y_i) & x \neq d_j, \\
      -\sum_{y_i \neq c_j} l(y_i) &  x = d_j.
    \end{cases}
  \end{equation}
  The defect of an ``external'' disk $f$ with $\image f \not \subset
  U$ is defined to be the defect of $\rho \circ f$.  Note that, by
  Lemma~\ref{lem:defect-bound} and the proof of
  Proposition~\ref{prop:equiv-d} below, the defect is always
  non-positive.  Together with the energy condition, this implies:
  \begin{equation} \label{eqn:energy}
      \delta(f(z)) \geq \sum_j \delta(f(w_j)).
  \end{equation}
\end{rem}

\begin{rem}
  Geometrically, a lattice point $p$ in $\lambda(x)$ corresponds with
  the generator $x^{\delta(p)}$ in \alg.  The equivalence relation for
  $\Delta^U$ allows a single equivalence class to range over all
  possible combinations of corners in each lattice, while fixing the
  macroscopic geometry of the disk.  Thinking of an equivalence class
  in $\Delta^U$ as a sum over all of these possibilities leads to a
  power series formula.  This idea is realized by the following
  proposition:
\end{rem}

\begin{prop}
  \label{prop:equiv-d}
  The total differential \df\ defined in Section~\ref{sec:full-diffl}
  is equivalent to:
  \begin{equation} \label{eqn:new-tot-diffl}
    \df \ps{x} = \sum_{f \in \Delta^U(x; y_1, \ldots,
      y_m)}
    \ps{y}_1 \cdots \ps{y}_m \nov^{-n(f; x, y_1, 
      \ldots, y_m)}.
  \end{equation}
\end{prop}

\begin{figure}[tbp]
  \centerline{\input{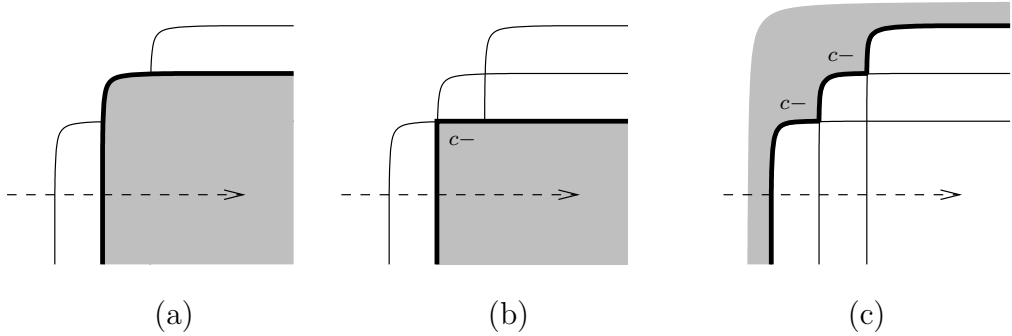}}
  \caption{Corners of $f$ in $\lambda(c)$ when the transverse direction 
    points into $f$ and $f$ (a) has no corners at a $c$ or (b) has a
    corner at a $c$; (c) the corners when the transverse direction
    points out of $f$.}
  \label{fig:multiple-c}
\end{figure}

\begin{proof}  There are two types of maps in $\Delta^U$:  those who
  image lies entirely inside $U$ and the rest.  As will be stated in
  Claims~\ref{clm:thick-int-diffl} and \ref{clm:thick-ext-diffl}, the
  former type of map gives \dfi, while the latter gives \dfe.

  \begin{claim} \label{clm:thick-int-diffl}
    \begin{equation}
      \dfi \ps{x} = \sum_{\substack{f \in \Delta^U(x;y_1, \ldots, y_m)
          \\
          \image f \subset U}}
        \ps{y}_1 \cdots \ps{y}_m \nov^{-n(f; x, y_1, 
      \ldots, y_m)}.
    \end{equation}
  \end{claim}
  
  Suppose that $f \in \Delta^U(x;y_1, \ldots, y_m)$ satisfies $\image
  f \subset U$.  Recall from Section~\ref{sec:int-diffl} that the
  terms in the internal differential correspond to flowlines,
  half-flowlines, and vertex flowlines in $\Gamma$.  The first step in
  proving the claim is to show that the image of $\rho \circ f$ (which
  is a line segment inside of $\Gamma$) coincides with one of these
  three types of flowlines.  This follows from three facts:
  \begin{enumerate}
  \item The condition that $f|_{\partial D}$ is a section of $\nu
    \Gamma_*$ away from the special points shows that $\image
    f|_{\partial D}$ must run straight from one lattice to another
    without turning around mid-edge.
  \item Since $f$ is an immersion, Figure~\ref{fig:cd-lattices} shows
    that $\image f$ cannot pass through the lattices $\lambda(c)$ or
    $\lambda(d)$. In particular, if $f$ has a corner at a $c$ or a
    $d$, then $\image (\rho \circ f)$ must have an end there.
  \item Since all of the corners of $\image f$ cover but one quadrant,
    $\image (\rho \circ f$) can only have an end at or pass straight
    through a double point.
  \end{enumerate}
  
  The next step is to show that ``thickened flowlines'' such as $f$
  give the terms in the definition of the internal differential.  This
  requires a careful examination of the structure of the ends of each
  type of thickened flowline.  Note that the description of the
  behavior of the ends of the flowlines in the lattices is unique up
  to the equivalence relation in
  Definition~\ref{defn:thick-flowlines}.
  
  Suppose that $\image (\rho \circ f)$ coincides with a full flowline.
  Such a thickened flowline has either one or two corners in
  $\lambda(c)$.  In either case, an examination of
  Figure~\ref{fig:cd-lattices}(a) shows that at least one of the
  corners must be positive. Thus, by the third condition of the
  definition of $\Delta^U$ and Figure~\ref{fig:cd-lattices}, a
  thickened flowline may only have a single negative corner in
  $\lambda(d)$. See Figure~\ref{fig:flowline}.
    
  Regardless of the transverse direction, there are thick flowlines
  with a single positive corner that leave the diagram of $\lambda(c)$
  in Figure~\ref{fig:cd-lattices} both down and to the right.  The
  latter is shown at the top of in Figure~\ref{fig:flowline}.  In the
  notation of Figure~\ref{fig:internal-diffl}, these give the
  $\ps{d}_1' + \ps{d}_2'$ terms in $\dfi \ps{c}$.
    
  If $f$ is a thick flowline that has two corners in $\lambda(c)$ and
  leaves the lattice to the right, then it has the form depicted at
  the bottom of Figure~\ref{fig:flowline}.  In counter-clockwise
  order, the flowline has a positive corner in $\lambda(c)$, a
  negative corner in $\lambda(d)$, and finally a negative corner back
  in $\lambda(c)$.  This gives the word $\ps{d}_2' \ps{c}$ in $\dfi
  \ps{c}$.  The order of the \ps{c} and \ps{d} is reversed for
  flowlines that leave downward, which gives the $\ps{c}\ps{d}_1'$
  term in $\dfi \ps{c}$.  Overall, the full flowlines contribute
  $(1+\ps{c}) \ps{d}_1' + \ps{d}_2' (1+\ps{c})$ to $\dfi \ps{c}$.
    
  \begin{figure}[tbp]
    \centerline{\input{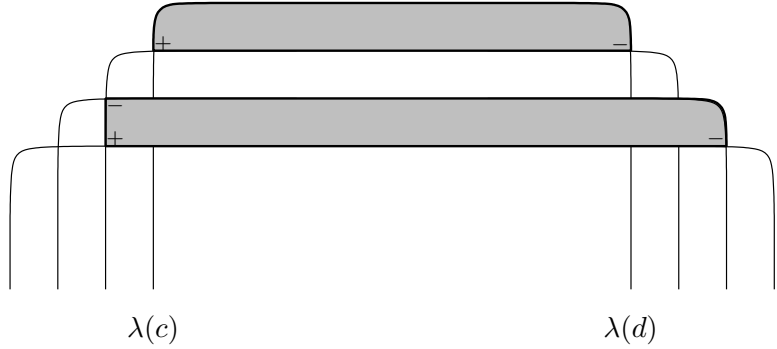}}
    \caption{The corners of full flowlines in $\Delta^U$.}
    \label{fig:flowline}
  \end{figure}
  
  The analysis of the other types of flowlines in
  Figure~\ref{fig:flowline} is similar.  Adding up the contributions
  yields \dfi\, and hence proves Claim~\ref{clm:thick-int-diffl}.
  
  \begin{claim} \label{clm:thick-ext-diffl}
    \begin{equation}
      \dfe \ps{x} = \sum_{\substack{f \in \Delta^U(x;y_1, \ldots, y_m)
          \\
          \image f \not\subset U}}
        \ps{y}_1 \cdots \ps{y}_m \nov^{-n(f; x, y_1, \ldots, y_m)}.
    \end{equation}
  \end{claim}
  
  Let $g \in \Delta(x; y_1', \ldots y_n')$. It suffices to show that
  the terms in $\dfe \ps{x}$ that come from $g$ are the same as those
  in
  \begin{equation} \label{eqn:ext-f}
    \sum_{\substack{f \in \Delta^U(x;y_1, \ldots, y_m)
        \\
        \rho \circ f = g}}
    \ps{y}_1 \cdots \ps{y}_m \nov^{-n(f; x, y_1, \ldots, y_m)}.
  \end{equation}
  Let $f$ be a map that contributes to the sum in (\ref{eqn:ext-f}).
  The maps $\rho \circ f$ and $g$ share their positive corner and
  their negative corners at double points.  Further, up to
  equivalence, $g$ uniquely determines $\image f$ in $\lambda(x)$ and
  in the double point lattices.  The structure of the half lattices
  $\lambda(c)$ and $\lambda(d)$ shows that $\rho \circ f$ can send
  neither a positive marked point to a $c$, nor a negative marked
  point to a $d$.  The proof of the claim now rests on the
  reconciliation of the marked points that $f$ and $g$ send to $c$
  points.
  
  The fourth and final condition in the definition of $\Delta(x; y_1',
  \ldots y_n')$ states that if $g|_{\partial D^2}$ crosses a special
  point labeled $c$, then that $c$ must be the image of a marked
  point.  That requirement lead to a $(1+\ps{c})$ or a
  $(1+\ps{c})^{-1}$ (depending on the transverse direction at $c$) in
  the word generated by $g$.  If the transverse direction at $c$
  points into $\image g$, then, as shown in
  Figures~\ref{fig:multiple-c}(a) and (b), the map $f$ may either pass
  through $\lambda(c)$ without a marked point or may have a single
  negative $c$ corner.  The former map gives the $1$, while the latter
  gives the \ps{c} in $(1+\ps{c})$.  On the other hand, if the
  transverse direction points out of $\image g$, then, as shown in
  Figure~\ref{fig:multiple-c}(c), the map $f$ may have any number of
  negative corners in $\lambda(c)$.  Up to equivalence, $f$ is
  determined by the number of $c$ corners in each lattice
  $\lambda(c)$.  Since
  \begin{equation*}
    (1+\ps{c})^{-1} = 1 + \ps{c} + \ps{cc} +  \ps{ccc} + \cdots,
  \end{equation*}
  summing the contributions of all possible maps $f$ gives the same
  terms as does $g$ in the original definition.  This completes the
  proof of Claim~\ref{clm:thick-ext-diffl}, and hence the proof of
  Proposition~\ref{prop:equiv-d}.
\end{proof}

\subsection{Gluing Broken Disks}
\label{sec:gluing}

In light of the unified definition of the differential, the definition
of a broken disk needs to be modified.  Suppose that $f \in
\Delta^U(x; y_1, \ldots, y_n)$ and $g \in \Delta^U(y_i; v_1, \ldots,
v_m)$.  Denote the marked points in the domain of $f$ by $\{z, w_1,
\ldots, w_n\}$ and those in the domain of $g$ by $\{z', w_1', \ldots,
w_m'\}$.  The pair $(f,g)$ is a \textbf{broken disk at \ps{x}} if
$f(w_i) = g(z')$.  Up to the equivalence relation in $\Delta^U$,
certain extra assumptions may be made about the form of $f$ and of
$g$.  Since the corners of $f$ and $g$ match at $y_i$, the images of
$f$ and $g$ share an edge $e$ in the lattice $\lambda(y_i)$.  Assume
that they share this edge until either $f$ or $g$ has another corner.
The edge $e$ has one end in $\lambda(y_i)$ and another in either
$\lambda(y_{i-1})$ or $\lambda(v_1)$, whichever comes
first.\footnote{If $i=1$, then $y_0$ is taken to be $x$.  Similarly,
  if $m=0$, then $v_1$ is taken to be $y_i$.  Clearly, it suffices to
  treat the case $n \geq 1$.}  The notion of which corner comes first
is unambiguous if $v_1 \neq y_{i-1}$.  If, on the other hand, $v_1 =
y_{i-1}$, then the geometry of the lattices and the energy condition
dictate the order of the corners.

\begin{figure}[tbp]
  \centerline{\input{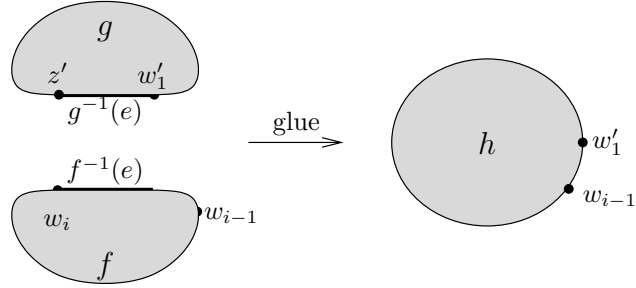}}
  \caption{Gluing the domains of $f$ and $g$ along $f^{-1}(e)$ 
    and $g^{-1}(e)$.}
  \label{fig:glued-domains}
\end{figure}

For notational convenience, suppose that the ends of $e$ lie in
$\lambda(y_i)$ and $\lambda(v_1)$.  Glue the domains of $f$ and $g$
along $f^{-1}(e)$ and $g^{-1}(e)$ as shown in
Figure~\ref{fig:glued-domains}.  Remove the marked points $w_i$ and
$z'$, but retain $w_1'$.  Define a new map $h: D^2 \to F$ by piecing
together $f$ and $g$ on the glued domain and smoothing as necessary.
By the discussion above, $h$ covers three quadrants at $h(w'_1)$; this
is the \textbf{obtuse corner} for $h$.  Say that the obtuse corner is
\textbf{positive} if two of the three quadrants covered are positive,
and \textbf{negative} otherwise.  The result is an \textbf{obtuse
  disk}, which satisfies all of the requirements of
Definition~\ref{defn:thick-flowlines} for $\Delta^U(x;y_1, \ldots,
v_1, \ldots, v_m, \ldots, y_n)$, except for the existence of one
obtuse corner. The map $h$ satisfies the energy condition because:
\begin{equation} \label{eqn:h-energy}
  \begin{split}
    \delta(h(z)) &= \delta(f(z)) \\
    &= \sum_j \delta(f(w_j)) - n(f) \\
    &= \sum_{j\neq i} \delta(f(w_j)) + \sum_k \delta(g(w_k')) - n(f) -
    n(g) \\
    &= \sum_j \delta(h(w_j)) - n(h).
  \end{split}
\end{equation}

The last line requires that the defects combine correctly:
\begin{multline}
  n(f; x, y_1, \ldots, y_n) + n(g;y_i, v_1, \ldots, v_m) =\\ n(h;x,
  y_1, \ldots, v_1, \ldots, v_m, \ldots, y_n).
\end{multline}
This follows directly from the definition of the defect (and equation
(\ref{eqn:flowline-defect}) if any of $f$, $g$, or $h$ are flowlines).
Thus, the $-n(h; \cdots)$ is the exponent of \nov\ in equation
(\ref{eqn:d2-term}), regardless of which broken disk was glued to form
it.

\subsection{Degeneration of Obtuse Disks}
\label{sec:degeneration}

\begin{figure}[tbp]
  \centerline{\input{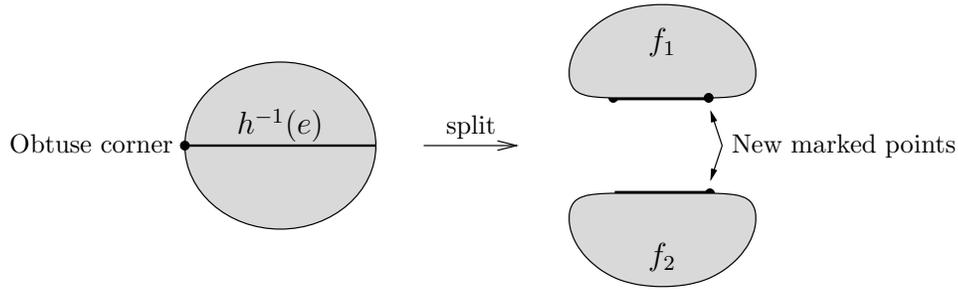}}
  \caption{Splitting the domains along $h^{-1}(e)$.}
  \label{fig:domain-split}
\end{figure}

It remains to show that any obtuse disk $h$ splits into exactly two
different broken disks.  Suppose that $h$ has a positive corner at $x$
and negative corners at $y_1, \ldots, y_n$.  Suppose further that the
obtuse corner lies in $\lambda(o)$, which may be either $\lambda(x)$
or some $\lambda(y_i)$. There are two line segments in $\lambda(o)$
that start at the obtuse corner and point into the interior of $\image
h$.  Each segment $e$ splits $h$ into a broken disk as follows: extend
$e$ along the lattice and, if necessary, by sections of $\nu \Gamma_*$
and along subsequent lattices until one of the following occurs:
\begin{enumerate}
\item The segment intersects the smooth part of the boundary of the
  obtuse disk,
\item The segment intersects itself, or
\item The segment returns to the obtuse corner.
\end{enumerate}
See Figure~\ref{fig:chv-splitting}. Notice that all three ending
conditions must occur within some lattice $\lambda(p)$, though there
is some ambiguity about the strand upon which $e$ enters $\lambda(p)$.
The preimage $h^{-1}(e)$ divides the domain $D^2$ into two sub-disks
$D^2_1$ and $D^2_2$.  Parameterize $e$ by $[0,1]$ so that $e(0) \in
\lambda(o)$ and $e(1) \in \lambda(p)$. Place an extra marked point at
$h^{-1}(e(1))$ in each of the new domains, as shown in
Figure~\ref{fig:domain-split}.  Let $f_j = h|_{D^2_j}$ for $j=1, 2$.
One of the $f_i$ has a positive corner at $e(1)$; the other a negative
corner at $e(1)$.\footnote{It is possible that the second disk has two
  negative corners in case (3).} Since there are exactly two possible
splitting segments $e$, the following lemma completes the proof of
Theorem~\ref{thm:dga}:

\begin{figure}[tbp]
  \centerline{\input{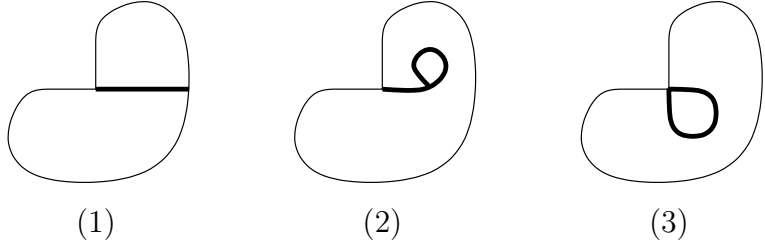}}
  \caption{The three ending conditions for the edge $e$ (shown in bold).}
  \label{fig:chv-splitting}
\end{figure}

\begin{lem} \label{lem:broken-disk}
  There exists a unique choice for the behavior of $e$ and $h$ (up to
  equivalence) in $\lambda(o)$ so that the pair $(f, g)$ is a
  broken disk.
\end{lem}

\begin{proof}
  Since $f(h^{-1}(e(1))) = g(h^{-1}(e(1)))$, it is sufficient to show
  that $f$ and $g$ are in the appropriate $\Delta^U$ space.  Consider
  each of the conditions in Definition~\ref{defn:thick-flowlines} in turn.

  \subsubsection{Conditions 1 and 2}
  \label{sec:conditions-12}
  Both $f$ and $g$ are immersions since $h$ is. They inherit the first
  two properties of Definition~\ref{defn:thick-flowlines} from $h$ and
  the construction of $e$ in $\lambda(o)$ and $\lambda(p)$.  

  \subsubsection{Condition 3}
  \label{sec:condition-3}
  The third property in Definition~\ref{defn:thick-flowlines} requires
  the existence of a unique (up to equivalence) choice of $e$ and $h$
  so that each $f_j$ has exactly one positive corner.  After the
  splitting construction, there are two positive corners to distribute
  between the $f_j$: one from $h$ and the other at $e(1)$.  Thus, for
  the third condition, the lemma reduces to showing that there exists
  a unique choice of behavior for $e$ and $h$ in $\lambda(o)$ such
  that the positive corner at $e(1)$ does not lie in the same disk as
  does the positive corner from $h$.  There are three cases:
  \begin{description}
    
  \item[Neither $\image f$ nor $\image g$ lies in $U$] In this case,
    $\lambda(p)$ must be a double-point lattice.  Since the splitting
    edge $e$ must leave the lattice $\lambda(o)$, there is complete
    freedom of choice about where $e$ enters $\lambda(p)$.  In
    particular, $e$ may end at either an $a$ corner or a $b$ corner,
    and only one of these choices will place the positive corner at
    $e(1)$ in the correct disk.
    
  \item[Both $\image f$ and $\image g$ lie in $U$] If $\image f_j
    \subset U$, then the argument at the beginning of the proof of
    Proposition~\ref{prop:equiv-d} shows that the image of $\rho \circ
    f_j$ coincides with that of a flowline. It is possible, however,
    for $f_j$ to have more or less than one positive corner.  For
    convenience, call such maps $f_j$ ``flowlines'' for the rest of
    the section, even though they might not contribute to \dfi.

    With this in mind, the second case follows from:

    \begin{claim} \label{clm:all-neg-flowline}
      No flowline can have all negative corners.
    \end{claim}

    \begin{proof}
      An inspection of the lattices in Figures~\ref{fig:ab-lattices}
      and \ref{fig:cd-lattices} shows that no vertex flowline can have
      all negative corners, and that every flowline that has an end in
      $\lambda(c)$ has at least one positive corner.  This leaves the
      case of a half-flowline from a double point to a $d$, but any
      such half-flowline must have at least one positive corner in
      either $\lambda(a)$ or $\lambda(b)$.
    \end{proof}
    
  \item[Only one of $\image f$ or $\image g$ lies in $U$] In other
    words, one of the $f$ or $g$ is a flowline (say that $f$ is) and
    the other is an external disk.  Claim~\ref{clm:all-neg-flowline}
    shows that $f$ has at least one positive corner.  The goal is to
    show that $f$ has \emph{exactly} one positive corner.  Suppose,
    instead, that $f$ has two positive corners for \emph{any} choice
    of the behavior of $e$ and $h$ in $\lambda(p)$, and proceed by
    contradiction.
  
    First, suppose that $f$ is a vertex flowline in some $\lambda(d)$.
    If this were the case, then $e(1)$ would also be in $\lambda(d)$,
    so $g$ would have a corner in $\lambda(d)$.  This corner must be
    positive, which contradicts the assumption that the flowline
    contains both positive corners.  Thus, $f$ cannot be a vertex
    flowline in $\lambda(d)$.
    
    Next, suppose that $f$ is a full flowline.  As before, $g$ cannot
    have any corners in $\lambda(d)$, so any corner of $f$ in
    $\lambda(d)$ is also a corner of $h$.  The energy condition and a
    quick examination of Figure~\ref{fig:cd-lattices} show that while
    $f$ must have a negative corner in $\lambda(d)$, it cannot have a
    positive corner there. Thus, the flowline can only have a single
    negative corner in $\lambda(d)$.  This gives a contradiction, for
    a full flowline can have at most one positive corner in
    $\lambda(c)$.
    
    If $f$ is a half-flowline to a $d$, the arguments above show that
    $f$ must have two positive corners in $\lambda(a) \cup
    \lambda(b)$.  Since the obtuse disk has but one positive corner,
    one of the positive corners of the flowline must be at $e(1)$.  As
    shown in Figure~\ref{fig:half-to-d}, either shifting $e$ or
    shifting the entire obtuse disk $h$ in $\lambda(a) \cup
    \lambda(b)$ gives a flowline with only one positive corner, which
    is a contradiction.

    \begin{figure}[tbp]
      \centerline{\input{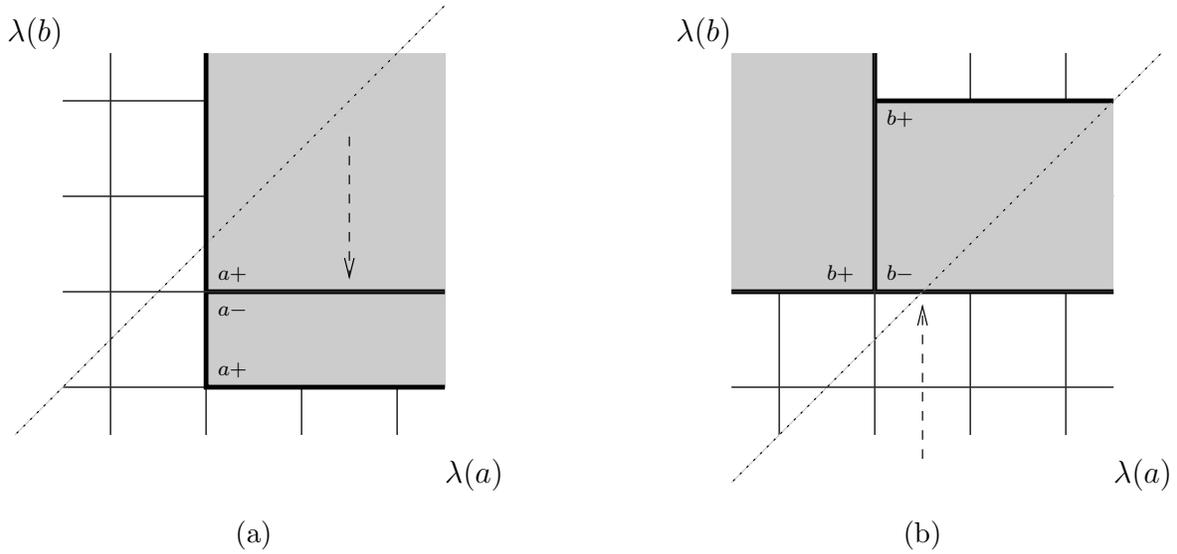}}
      \caption{(a) Shifting $e$ down until it meets the boundary of $h$ at an
        $a$ crossing yields a flowline with only one positive corner;
        (b) shifting the boundary of the obtuse disk up until $e$ meets
        the boundary at a $b$ corner yields a flowline with only one
        positive corner.  The other cases differ by a reflection.}
      \label{fig:half-to-d}
    \end{figure}
    
    Next, consider a half-flowline from a $c$.  In order to have two
    positive corners, the flowline must have corners that cover either
    both an $a-$ and an $a+$ or both a $b+$ and a $b-$ in the double
    point lattice.  If both of the corners in the double point lattice
    are corners of the obtuse disk, then an inspection of
    Figure~\ref{fig:ab-lattices} shows that the obtuse disk would not
    satisfy (\ref{eqn:energy}).  If $e(1)$ lies at the negative corner
    of the flowline in the double point lattice, then the external
    disk would have a positive corner, which is a contradiction.
    Finally, if the flowline and the external disk are glued together
    at the flowline's positive corner, then the edge $e$ may be
    shifted so that the external disk has a positive corner; see
    Figure~\ref{fig:corner-at-ab}.  In sum, the existence of any
    half-flowline from a $c$ with two positive corners leads to a
    contradiction.
    
    \begin{figure}[tbp]
      \centerline{\input{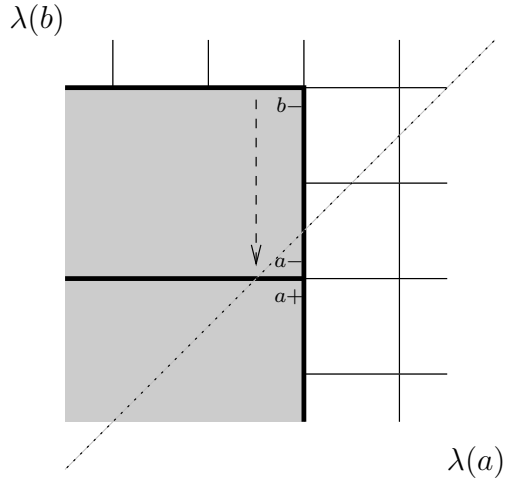}}
      \caption{The segment $e$ may be shifted down so that the external
        disk has a positive corner in $\lambda(a) \cup \lambda(b)$.}
      \label{fig:corner-at-ab}
    \end{figure}
    
    Suppose that the flowline is a vertex flowline in $\lambda(c)$;
    the possibilities are pictured in Figure~\ref{fig:bad-vertex-c}.
    The vertex flowline must have a positive corner at $e(1)$;
    otherwise, $g$ would have a positive corner.  Once the gluing is
    performed, however, there is always a negative corner that is
    further away from the lattice edge than the remaining positive
    corner.  But this is impossible by (\ref{eqn:energy}), giving the
    desired contradiction.

    \begin{figure}[tbp]
      \centerline{\input{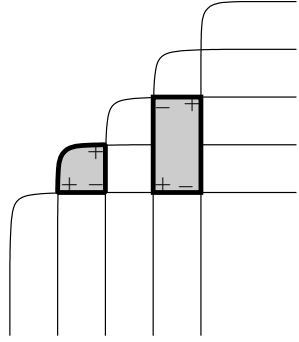}}
      \caption{Vertex flowlines in $\lambda(c)$.}
      \label{fig:bad-vertex-c}
    \end{figure}

    \begin{figure}[tbp]
      \centerline{\input{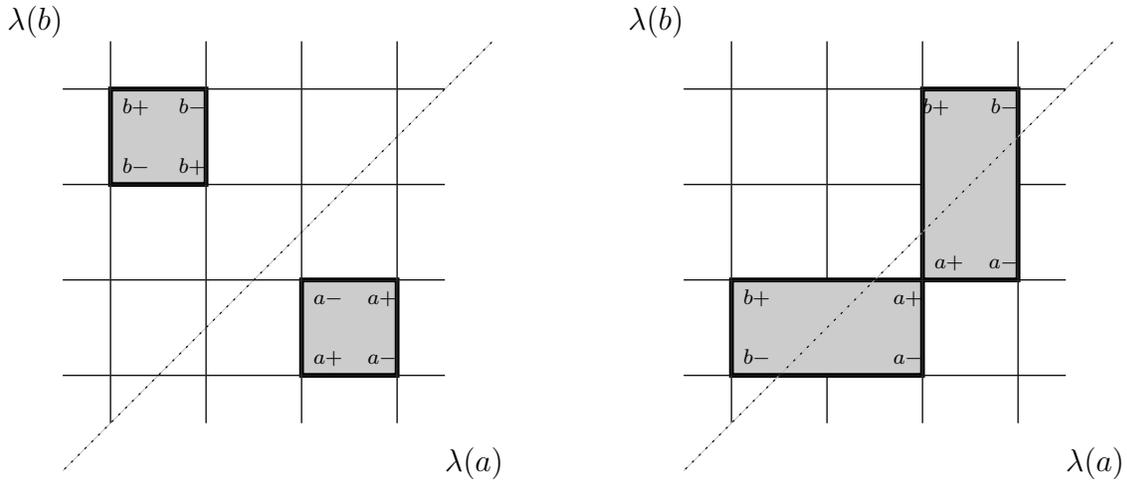}}
      \caption{Vertex flowlines in $\lambda(a) \cup \lambda(b)$.}
      \label{fig:bad-vertex-ab}
    \end{figure}
    
    Finally, Figure~\ref{fig:bad-vertex-ab} depicts the possible vertex
    flowlines in $\lambda(a) \cup \lambda(b)$.  All of these satisfy:

    \begin{claim} \label{clm:vertex-flowline}
      The sum of the distances of the to the lattice edge, counted
      with sign, is zero.
    \end{claim}
    
    \begin{proof}
      Put coordinates on the full lattice $\lambda(a) \cup \lambda(b)$
      by placing the origin at one of the nodes of $\lambda(b)$ that
      lies closest to the lattice edge.  Then the distance $\delta_{ab}$
      to the lattice edge is given by:
      \begin{equation} \label{eqn:delta}
        \delta_{ab}(m,n) = \begin{cases} m-n-1 & (m,n) \in \lambda(a) \\
          n-m & (m,n) \in \lambda(b).
        \end{cases}
      \end{equation}
      In every vertex flowline in Figure~\ref{fig:bad-vertex-ab}, the
      corners in $\lambda(a)$ can be paired as $a+/a-$ (if they exist).
      Thus, the $-1$ terms that appear in (\ref{eqn:delta}) cancel from
      $\sum_{\text{corners}} \delta_{ab}(i,j)$.  Every coordinate
      appears in the sum twice with opposite sign since consecutive
      corners are either labeled with the same letter that have opposite
      signs or are labeled with opposite letters that have the same
      sign.  The claim follows.
    \end{proof}
    
    As in the argument for the vertex flowlines in $\lambda(c)$, the
    flowline must have a positive corner at $e(1)$.  The claim above
    shows that the sum of the distances from the negative corners to
    the lattice edge is greater than that of the remaining positive
    corner.  This contradicts (\ref{eqn:energy}), and finishes the
    proof that $f$ and $g$ satisfy the third condition of
    Definition~\ref{defn:thick-flowlines}.

  \end{description}

  \subsubsection{The Energy Condition}
  \label{sec:energy-condition}
  
  The last thing to show is that $f$ and $g$ satisfy the energy
  condition.  By (\ref{eqn:h-energy}) and the assumption that $h$
  satisfies the energy condition, if one of the two satisfy the energy
  condition, then the other will as well.  The goal, then, is to show
  that the behavior of $e$ and $h$ can be chosen so that either $f$ or
  $g$ satisfies the energy condition.
  
  First, suppose that both $f$ and $g$ are external disks, and that
  $f$ contains the positive corner from $h$.  There is sufficient
  freedom to shift $e$ in $\lambda(p)$ so that $\delta(g(z')) =
  \delta(e(1)) = \sum_j \delta(g(w_j)) - n(g)$.  Thus, after the
  shift, $g$ satisfies the energy condition.
  
  From now on, assume that at least one of $f$ or $g$ is a flowline.
  As in Claim~\ref{clm:vertex-flowline}, the sum of the distances of
  the corners of a true flowline can be easily read off from its
  combinatorics.  If the sum of the distances, counted with sign, is
  equal to the negation of the defect, then the flowline satisfies the
  energy condition.  The following claim results from direct
  computations using the usual coordinates on full and half lattices:

  \begin{claim} \label{clm:flowline-dist}
    Suppose $f$ is a flowline that satisfies the first three
    conditions of Definition~\ref{defn:thick-flowlines}.
    \begin{enumerate}
    \item If $f$ is a vertex flowline, then $f$ satisfies the energy
      condition.
    \item If $f$ is a full flowline or a half-flowline to a $d$, then
      the contribution to the sum of the distances by the corners in
      an end in $\lambda(x)$ is equal to the width of the flowline in
      $\lambda(x)$.  The width contributes positively at the end with
      the positive corner, and negatively in $\lambda(d)$.
    \item If $f$ is a half-flowline from a $c$, then the contribution
      to the sum of the distances in $\lambda(c)$ is the width of the
      flowline and the contribution at the double point is the
      negation of the width of the flowline minus one.
    \end{enumerate}
  \end{claim}
  
  The claim shows that the energy condition holds if the flowline is a
  vertex flowline.  Further, if $e$ and $h$ can be arranged so that
  the widths at the ends of the flowline match up, then the flowline
  would satisfy the energy condition, and Lemma~\ref{lem:broken-disk}
  would follow.  Note that the extra $-1$ at the double-point end of
  the half-flowline from a $c$ corresponds with the defect of the
  half-flowline.  Thus, it remains to prove that $e$ and $h$ may be
  chosen so that the widths of the ends of the flowlines match up.

  There are two broad cases to consider.
  \begin{description}
  \item[The segment $e$ runs along the length of the flowline $f$]
    First, suppose that $e(1)$ lies in $\lambda(c)$ or $\lambda(d)$.
    If the other component of the purported broken disk is an external
    disk, then $e$ may be shifted to realize any width for the
    flowline.  For an example of this, see
    Figure~\ref{fig:parallel-e}.  In particular, $e$ can be shifted
    until the width of $f$ in $\lambda(c)$ or $\lambda(d)$ matches the
    width at the other end of $f$ (which is determined by the choice
    of $h$).
    
    So just consider the case where both $f$ and $g$ are flowlines.
    If $e(1)$ lies in $\lambda(d)$, then, since either $f$ or $g$ must
    have a positive corner there, one component is a vertex flowline
    and hence already satisfies the energy condition.  Otherwise, if
    $e(1)$ lies in $\lambda(c)$, then $h$ must have a positive corner
    in $\lambda(c)$.  The energy condition on $h$ implies that the
    width of $h$ in $\lambda(c)$ bounds the width of $h$ in any other
    lattice.  In particular, the width of $h$ bounds the width of $f$,
    so $e$ may be chosen in $\lambda(c)$ so as to realize the
    appropriate width.

    \begin{figure}[tbp]
      \centerline{\input{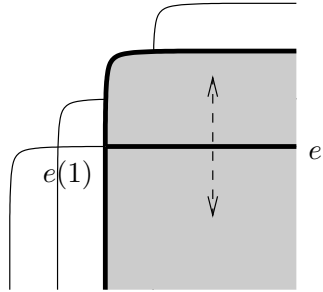}}
      \caption{If $e(1) \in \lambda(d)$, the obtuse disk has no
        corner in $\lambda(d)$, and the other disk is external, then
        the segment can be shifted to realize any width for the
        flowline.}
      \label{fig:parallel-e}
    \end{figure}

    \begin{figure}[tbp]
      \centerline{\input{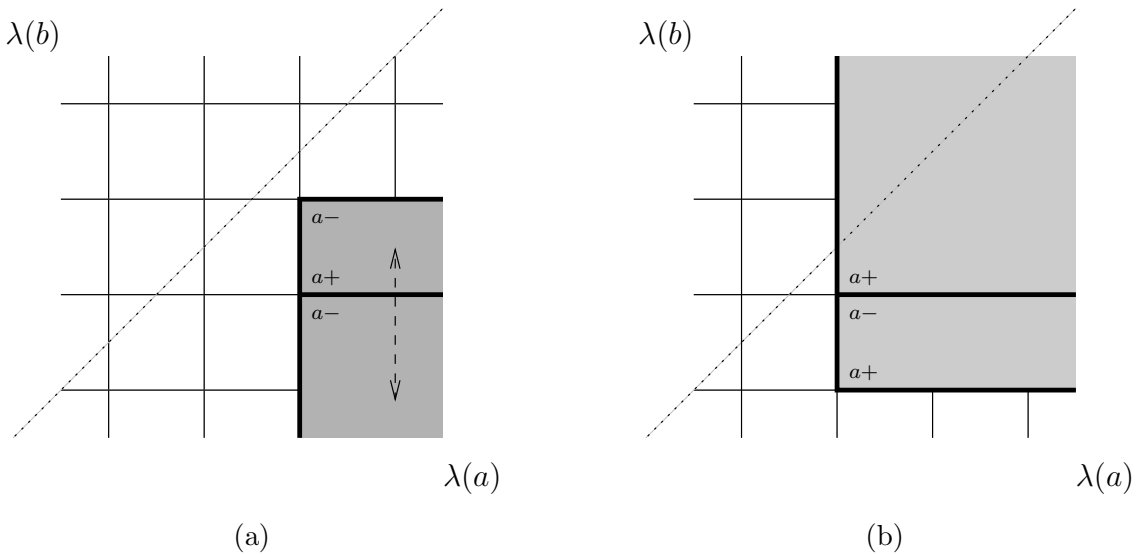}}
      \caption{If $e(1) \in \lambda(a) \cup \lambda(b)$ and the
        flowline is a half-flowline to a $d$, then either (a) $e$ may
        be shifted to realize any width, or (b) the width at $d$ is
        bounded by $\delta(a+)$ so any $e$ as shown works.}
      \label{fig:parallel-e-2}
    \end{figure}
    
    Finally, if $e(1)$ lies in a double-point lattice, similar
    arguments show that, again, either $e$ may be shifted to realize
    any width or that the widths in other lattices are bounded by the
    width in the double-point lattice.
    
  \item[The segment $e$ forms the end of the flowline] The arguments
    are entirely similar in this case.  Instead of shifting the
    segment $e$ in the lattice containing $e(1)$, it is necessary to
    shift the side of $h$ in which $e(1)$ lies. For example, see
    Figure~\ref{fig:half-to-d}(b).  As before, either the positive
    corner of $h$ lies in the flowline and gives a necessary bound on
    the width of the flowline elsewhere, or the side of $h$ in which
    $e(1)$ lies has complete freedom to shift to the required width.

  \end{description}

  This finishes the proof of Lemma~\ref{lem:broken-disk}, and hence
  completes the proof of Theorem~\ref{thm:dga}.
\end{proof}

\section{Proof of Invariance}
\label{sec:invariance}

The final loose end is the proof of Theorem~\ref{thm:invariance},
i.e.\ that the stable tame isomorphism type of $(\alg, \df)$ is an
invariant of Legendrian knots.  The proof consists of somewhat
technical combinatorial arguments generalizing those in \cite{chv,
  ens} to the contact circle bundle setting.

\subsection{Choices at $c$ and $d$ Points}
\label{sec:cd-invariance}

\subsubsection{Transverse Directions at $c$}
\label{sec:transverse-dir}

The first part of Theorem~\ref{thm:invariance} states that the stable
tame isomorphism type of the DGA of a knot diagram is independent of
the choices of transverse directions at the edge points labeled with
$c_i$.  Let $c$ label one of these points, and let $(\alg, \df)$ and
$(\alg', \df')$ be the DGAs given by the two choices of a transverse
direction.  Define a homomorphism $\phi: \alg \to \alg'$ to be the
identity on all generators besides $c^k$, $k=1,2,\ldots$, and to map
the generators $c^k$ according to the power series formula
\begin{equation}
  \boldsymbol{\phi}(1+\ps{c}) = (1+\ps{c}')^{-1}.
\end{equation}
Since the generating power series \ps{c} contains no constant term,
this indeed gives a tame isomorphism.  The first part of
Theorem~\ref{thm:invariance} follows from:


\begin{lem}
  $\phi$ is a chain map.
\end{lem}

\begin{proof}
  First, let $\ps{x} \neq \ps{c}$ be a generating power series of
  \alg.  The term $(1+\ps{c})^{\pm 1}$ appears in $\df \ps{x}$ if and
  only if $(1+\ps{c}')^{\mp 1}$ appears in $\df' \ps{x}$, so $\phi$ is
  a chain map on \ps{x}.  That $\phi$ is a chain map on \ps{c} is the
  result of a direct computation using
  Definition~\ref{defn:internal-diffl}.
\end{proof}

\subsubsection{Choice of $c_1$}
\label{sec:c-position}

The second part of Theorem~\ref{thm:invariance} asserts that the
stable tame isomorphism type of the DGA of a knot diagram does not
depend on whether a given edge point has a $c_i$ or $d_i$ label.  To
prove this, it suffices to show that the tame isomorphism type of the
DGA remains unchanged when the $c$ and $d$ labels are cyclically
shifted; i.e. the point labeled with $d_i$ becomes $c_i$ and the point
labeled $c_i$ becomes $d_{i-1}$.

\begin{figure}[tbp]
  \centerline{\input{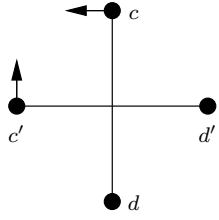}}
  \caption{Labels of special points around a double point in the  
    definition of $\phi$.}
  \label{fig:shift-diagram}
\end{figure}

The idea is to define a homomorphism $\phi: (\alg, \df) \to (\alg',
\df')$ via local contributions at each double point.  First, some
notation is necessary.  Let $\ps{y}_i$ be the contribution to $\dfi
\ps{c}$ of the half-flowline from $c_i$ to the double point between
$c_i$ and $d_i$.  Similarly, let $\ps{x}_i$ come from the
half-flowline from $c_i$ to the double point between $c_i$ and
$d_{i-1}$.  Suppose that the special points near a double point are
labeled as in Figure~\ref{fig:shift-diagram}.
\begin{itemize}
\item If $d$ passes through the double point,\footnote{In other words,
    if $d$ gets shifted to the position previously occupied by $c$.}
  then:
  \begin{equation}
    \boldsymbol{\phi} (\ps{d}) = \ps{d} + \ps{y}.
  \end{equation}
  Note that the \nov\ is included in $\ps{y}$. 
\item If $d'$ passes through the double point, then:
  \begin{equation}
    \boldsymbol{\phi} (\ps{d}) = \ps{d} + \ps{y'}.
  \end{equation}
\item If $c$ passes through the double point,
  then:
  \begin{equation} \label{eqn:local-c}
    \begin{split}
      \boldsymbol{\phi} (\ps{a}) &= \ps{a} (1+\ps{c})^{-1},\\
      \boldsymbol{\phi} (\ps{b}) &= (1+\ps{c}) \ps{b}.
    \end{split}
  \end{equation}
\item If $c'$ passes through the double point, then:
  \begin{equation} \label{eqn:local-c'}
    \begin{split}
      \boldsymbol{\phi} (\ps{a}) &= (1+\ps{c'})^{-1} \ps{a},\\
      \boldsymbol{\phi} (\ps{b}) &= \ps{b} (1+\ps{c'}).
    \end{split}
  \end{equation}
\end{itemize}
Only two of these possibilities can occur at any given double point.
Using the results of Section~\ref{sec:transverse-dir}, it is easy to
see that if the transverse direction at a $c$ changes in
Figure~\ref{fig:shift-diagram}, then $(1+\ps{c})^{\pm 1}$ changes to
$(1+\ps{c})^{\mp 1}$ in (\ref{eqn:local-c}) and (\ref{eqn:local-c'}).

The full map $\phi$ comes from taking the local contributions first
from the $d$ moves and then the $c$ moves into account, one at a time.
By Lemma~\ref{lem:tame-isom}, each of these steps is a tame
isomorphism of algebras.  It follows that $\phi$ is the composition of
finitely many tame isomorphisms, and hence is itself a tame
isomorphism of algebras.

It remains to show that $\phi$ is a chain map.  The proof may be
separated into two steps.

\begin{description}
\item[$\phi \dfe = \dfe' \phi$] First, suppose that a special point
  $w$ lies outside Figure~\ref{fig:shift-diagram}.  Further, suppose
  that a disk in $\dfe \ps{w}$ has a negative corner at the double
  point in Figure~\ref{fig:shift-diagram}.  If $d$ of $d'$ pass
  through the double point, then, locally, nothing changes in $\dfe
  \ps{w}$.  If $c$ passes through the double point and the boundary of
  the disk passes through $c$ before the shift, then:
  \begin{enumerate}
  \item The position of the $(1+\ps{c})^{\pm 1}$ factor with respect
    to the corner at the double point is the same in $\dfe \ps{w}$ and
    in the image of the double point corner under $\phi$.
  \item If $(1+\ps{c})^{\pm 1}$ appears in $\dfe \ps{w}$, then $\phi$
    contributes $(1+\ps{c})^{\mp 1}$.
  \end{enumerate}
  Thus, the $(1+\ps{c})^{\pm 1}$ factors cancel in $\phi \dfe \ps{w}$.
  On the other hand, after the shift of the $c$, the disk no longer
  passes through $c$, so there are no $(1+\ps{c})^{\pm 1}$ factors in
  $\dfe' \boldsymbol{\phi} \ps{w}$.  Hence, the local contributions of
  $\dfe' \phi \ps{w}$ and $\phi \dfe \ps{w}$ agree.
  
  If the boundary of the disk passes through $c$ after the shift, then
  the position and the exponent of $(1+\ps{c})^{\pm 1}$ are both the
  same for the external disk after the shift and for the image of the
  double point corner under $\phi$.  Thus, once again, the local
  contributions of $\dfe' \phi \ps{w}$ and $\phi \dfe \ps{w}$ agree.
  Since the argument for $c'$ is exactly the same, this finishes the
  case when $w$ lies outside of Figure~\ref{fig:shift-diagram} and
  $\dfe \ps{w}$ has a negative corner at the double point.

  Next, consider the special point $d$ in
  Figure~\ref{fig:shift-diagram}.  Suppose that $d$ passes through the
  double point and that a disk in $\dfe \ps{d}$ has a negative corner
  at the double point as in Figure~\ref{fig:shift-d}. Such disks
  represent words of the form $\ps{Va}$.  The terms in
  $\boldsymbol{\phi}(\ps{Va})$ differ from \ps{Va} at negative corners
  through which a $c$ passes.  This is exactly the situation
  considered above.  Thus, by abusing notation, write
  $\boldsymbol{\phi}(\ps{Va}) = \ps{Va}$.

  \begin{figure}[tbp]
    \centerline{\input{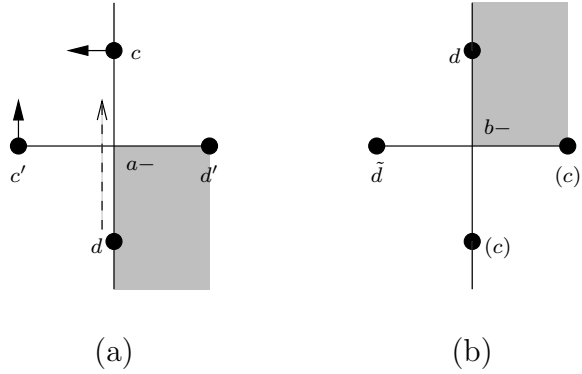}}
    \caption{(a) Disks in $\dfe \ps{d}$; (b) Disks in $\dfe' \ps{d}$.
      The points labeled $(c)$ are some unspecified $c_i$-labeled edge
      points.}
    \label{fig:shift-d}
  \end{figure}
  
  On the other hand, $\ps{y} = \ps{ba}\nov$ and so
  $\boldsymbol{\phi}(\ps{d}) = \ps{d} + \ps{b}\ps{a}\nov$.  The
  differential $\dfe'(\ps{d} + \ps{b}\ps{a}\nov)$ yields three types
  of disks on the right side of Figure~\ref{fig:shift-d}(b):
  \begin{enumerate}
  \item $\dfe' \ps{d}$ gives the disks pictured in
    Figure~\ref{fig:shift-d}(b).
  \item $\ps{b} \dfe'(\ps{a}) \nov$ matches the disks in $\dfe'
  \ps{d}$.
  \item $\dfe'(\ps{b})\ps{a}\nov$ matches $\ps{Va}$ and hence matches
    $\boldsymbol{\phi}\dfe \ps{d}$.
  \end{enumerate}
  The disks that lie on the other side of Figure~\ref{fig:shift-d}(a)
  match up via a parallel argument.  This finishes the proof that the
  local contribution near $d$ satisfies $\boldsymbol{\phi}\dfe \ps{d}
  = \dfe' \boldsymbol{\phi} \ps{d}$.  The proof for $d'$ is similar.

  A similar strategy works for the case of a disk in $\dfe \ps{b}$.
  For example, suppose that $c$ passes through the double point.
  Disks in the top left corner of Figure~\ref{fig:shift-diagram}
  contribute terms of the form $(1+\ps{c})\ps{W}$ to
  $\boldsymbol{\phi} \dfe \ps{b}$.  On the other hand, since $\dfe'
  \ps{c} = 0$,
  \begin{equation} \label{eqn:shift-b}
    \begin{split}
      \dfe' \boldsymbol{\phi}(\ps{b}) &= \dfe' \bigl( (1+\ps{c})\ps{b}
      \bigr) \\ &= (1+\ps{c}) \ps{W}.
    \end{split}
  \end{equation}
  Note that the two \ps{W} terms agree because, by the arguments
  above, $\phi$ is a chain map with respect to the local contributions
  at negative corners.  Further, if $c'$ passes through the double
  point as well, then the $(1+\ps{c}')$ terms in the analogue of
  (\ref{eqn:shift-b}) appear on the right side of \ps{b}.  Thus, the
  shifts of the $c$ and $c'$ do not interfere with each other
  algebraically.

  The proof for \ps{a} is the same.  This completes the proof that
  $\phi \dfe = \dfe' \phi$.
  
\item[ $\phi \dfi = \dfi' \phi$] First, consider the interaction
  between the internal differential and $\phi$ on \ps{d}.  Suppose
  that $d$ passes through the double point in
  Figure~\ref{fig:shift-diagram}.  Let $\tilde{d}$ be the $d$-labeled
  point that appears at the left of the diagram after the shift. See
  Figure~\ref{fig:shift-d}. Then $\ps{y} = \ps{b} \ps{a} \nov$ and:
  \begin{align*}
    \dfi' \boldsymbol{\phi} \ps{d} &= \dfi' (\ps{d} +
    \ps{b}\ps{a}\nov) \\ &= \ps{d}\ps{d} + \ps{d}\ps{b}\ps{a} \nov +
    \ps{b}\ps{a}\ps{d} \nov + 2 \cdot \ps{b}\tilde{\ps{d}}\ps{a} \nov
    + \ps{b}\ps{a}\ps{b}\ps{a} \nov^2 \\ &= \boldsymbol{\phi} \dfi
    \ps{d}.
  \end{align*}

  The next step is to prove that $\phi \dfi \ps{c} = \dfi' \phi
  \ps{c}$.  Consider the configuration in Figure~\ref{fig:shift-c},
  and, without loss of generality, suppose that $c$ shifts to the
  left.  Recall that $\ps{x}_i$ comes from the half-flowline from
  $c_i$ to the double point between $c_i$ and $d_{i-1}$.  The
  following straightforward lemma holds regardless of the labels of
  the points at the top and bottom of Figure~\ref{fig:shift-c}:

  \begin{lem} \label{lem:x-i}  
    \begin{align*}
      \boldsymbol{\phi}(\ps{x}_i) &= (1+\ps{c}_i)^{-1} \ps{x}_i
      (1+\ps{c}_i), \\ \boldsymbol{\phi}(\ps{y}_i) &= \ps{y}_i.
    \end{align*}
  \end{lem}

  \begin{figure}[tbp]
    \centerline{\input{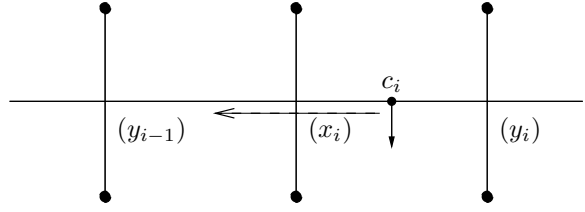}}
    \caption{The local configuration for $c$ shifting to the left.}
    \label{fig:shift-c}
  \end{figure}

  Now compute directly, using Lemma~\ref{lem:x-i}:
  \begin{equation*}
    \begin{split}
      \boldsymbol{\phi} \dfi \ps{c}_i &= \boldsymbol{\phi} \bigl(
      (1+\ps{c}_i)(\ps{d}_{i-1} + \ps{x}_i) + (\ps{d}_i + \ps{y}_i)
      (1+\ps{c}_i) \bigr) \\ &= (1+\ps{c}_i)(\ps{d}_{i-1} +
      \ps{y}_{i-1} + (1+\ps{c}_i)^{-1} \ps{x}_i (1+\ps{c}_i)) +
      \ps{d}_i (1+\ps{c}_i) \\ &= \dfi' \boldsymbol{\phi} \ps{c}_i.
    \end{split}
  \end{equation*}
  The last equality holds because after the shift, $c_i$ is on the
  other side of the double point that contributes $\ps{x}_i$, but the
  $c$ and $d$ labels on the crossing strand at that double point have
  switched places.  Thus, the double point still contributes
  $\ps{x}_i$.  This finishes the proof for \ps{c}.

  The internal differential on \ps{b} is more complicated.  For
  simplicity, only consider the case when $c$ and $d'$ pass through
  the double point; the other cases follow from similar direct
  computations.  Consider the configuration in
  Figure~\ref{fig:shift-b-int}.  In this case, $\ps{y} =
  \ps{b}_1\ps{a}_1 \nov$ and $\ps{y}' = \ps{ab}\nov$.  Computing
  directly,
  \begin{equation} \label{eqn:shift-b-1}
    \begin{split}
      \boldsymbol{\phi} \dfi \ps{b} &= \boldsymbol{\phi} \left(
        \ps{b}\ps{d}' + \ps{d}\ps{b} + \ps{b}\ps{a}\ps{b} \nov \right)
        \\ &= (1+\ps{c})\ps{b}(\ps{d}' + \ps{a}\ps{b}\nov) + (\ps{d} +
        \ps{b}_1 \ps{a}_1 \nov)(1+\ps{c})\ps{b} +
        (1+\ps{c})\ps{b}\ps{a}\ps{b} \nov \\ &=
        (1+\ps{c})\ps{b}\ps{d}' + (\ps{d} + \ps{b}_1 \ps{a}_1 \nov)
        (1+\ps{c})\ps{b}.  \\
    \end{split}
  \end{equation}
  On the other hand,
  \begin{equation} \label{eqn:shift-b-2}
    \begin{split}
      \dfi' \boldsymbol{\phi} (\ps{b}) &= \dfi' \bigl(
      (1+\ps{c})\ps{b} \bigr) \\ &= \bigl( (1+\ps{c})(\ps{d}_2 +
      \ps{b} \ps{a} \nov) + (\ps{d} + \ps{b}_1 \ps{a}_1 \nov)(1+
      \ps{c}) \bigr) \ps{b} \\ &\quad + (1+\ps{c})(\ps{b}\ps{d}' +
      \ps{d}_2\ps{b} + \ps{b}\ps{a}\ps{b} \nov) \\ &=
      (1+\ps{c})\ps{b}\ps{d}' + (\ps{d} + \ps{b}_1 \ps{a}_1 \nov)
      (1+\ps{c})\ps{b}.
    \end{split}
  \end{equation}
  If the local configuration near the double point labeled $(1)$
  changes so that the $(c)$ is on the left before the shift, then the
  order of the factors in \ps{y} switch places and
  (\ref{eqn:shift-b-1}) and (\ref{eqn:shift-b-2}) still agree.  Thus,
  $\dfi' \boldsymbol{\phi} (\ps{b}) = \boldsymbol{\phi} \dfi \ps{b}$.

  \begin{figure}[tbp]
    \centerline{\input{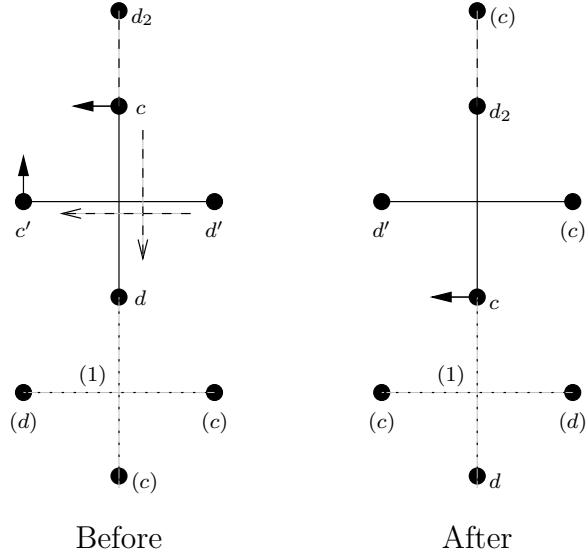}}
    \caption{The configuration for the internal differential of \ps{b} 
      when $c$ passes through the double point.  The parenthetical
      labels represent an arbitrary $c$ or $d$ label; exactly which
      one does not matter for this computation.}
    \label{fig:shift-b-int}
  \end{figure}

  The proof for \ps{a} is similar.  This completes the proof that
  $\phi \dfi = \dfi' \phi$.

\end{description}

\subsection{Legendrian Isotopy}
\label{sec:isotopy}

\subsubsection{Generic Isotopies}
\label{sec:generic-isotopy}

\begin{figure}[tbp]
  \centerline{\input{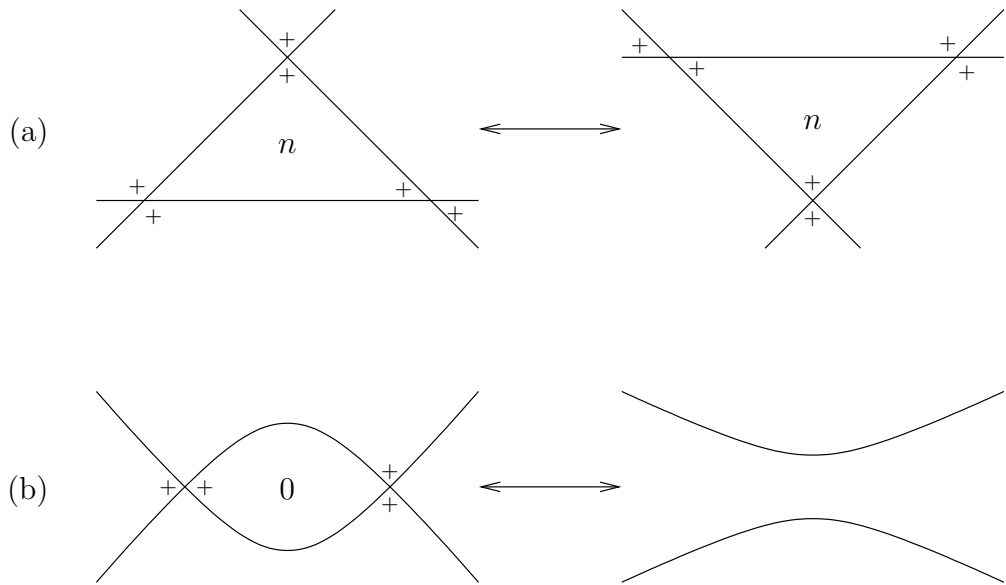}}
  \caption{(a) The triple-point move.  (b) The double-point move.}
  \label{fig:rm}
\end{figure}

The final part of Theorem~\ref{thm:invariance} asserts that the stable
tame isomorphism type of the DGA of a diagram of a Legendrian knot $L$
is invariant under Legendrian isotopy.  The following lemma translates
Legendrian isotopy of $L$ into combinatorial moves on the diagram
$(\Gamma^+_L, \vec{n})$:

\begin{lem}
  \label{lem:isotopy}
  If $L_0$ and $L_1$ are Legendrian isotopic Legendrian knots in $E$,
  then $(\Gamma^+_{L_0}, \vec{n}_0)$ and $(\Gamma^+_{L_1}, \vec{n}_1)$
  differ by a sequence of the local moves shown in
  Figure~\ref{fig:rm}.
\end{lem}

\begin{proof}
  Let $L_s$ be a Legendrian isotopy in $E$ between $L_0$ and $L_1$.
  Since the contact structure is transverse to the fibers, the isotopy
  projects to a homotopy of immersions $\gamma_s: S^1 \to F$, which
  may be assumed to be generic away from the transitions pictured in
  Figure~\ref{fig:rm}.
    
  The defects in Figure~\ref{fig:rm} need some justification.  Since
  the strands of $L$ are disjoint during both of these moves, the
  following inequality always holds for some fixed $\epsilon>0$:
  \begin{equation} \label{eqn:chord-lengths}
    \epsilon < l(a_i^0) < 2\pi -\epsilon,
  \end{equation}
  and similarly for $l(b_i^0)$.
  
  In the case of the triple-point move, in the limit as the area of
  the central triangle goes to zero, the equation for the defect
  reads:
  \begin{equation}
    2 \pi n = -l(a_1^0) + l(b_2^0) + l(b_3^0).
  \end{equation}
  Since the defect is an integer, the defect $n$ of the central
  triangle in Figure~\ref{fig:rm}(a) is either $0$ or $1$.  These two
  cases correspond to the two types of triple point moves for
  Legendrian knots in $\rr^3$; see \cite{chv,ens}.  A similar argument
  using equation (\ref{eqn:chord-lengths}) shows that the defect of
  the central lune in Figure~\ref{fig:rm}(b) must be $0$.
\end{proof}

The rest of this section is devoted to proving that the stable tame
isomorphism type of the DGA of a knot diagram is invariant under the
two moves in Figure~\ref{fig:rm}.

\subsubsection{Triple-Point Move}
\label{sec:triple-pt-invariance}

\begin{figure}[tbp]
  \centerline{\input{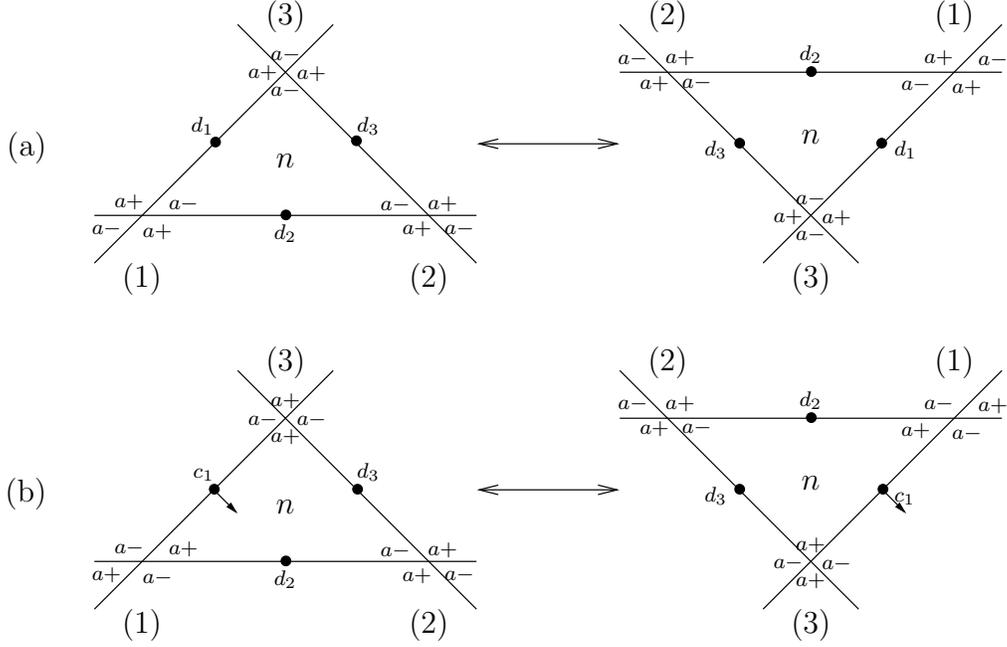}}
  \caption{The two possible configurations of the generators near a triple 
    point.  The parenthetical numbers are labels for the double
    points.  The DGAs for the left sides of (a) and (b) are denoted by
    $(\alg, \df)$, while the DGAs for the right sides are denoted by
    $(\alg', \df')$.}
  \label{fig:triple-pt}
\end{figure}

Up to choosing the position of $c_1$ and the transversal directions at
the $c$ points, the two possible configurations for the generators
near a triple-point move are pictured in Figure~\ref{fig:triple-pt}.
For the moment, consider only the diagrams in
Figure~\ref{fig:triple-pt}(a). Define a map $\phi: \alg \to \alg'$ by:
\begin{equation}
  \boldsymbol{\phi} (\ps{x}) = \begin{cases}
    \ps{a}_1 + \ps{b}_3 \ps{b}_2 \nov^n & \ps{x} = \ps{a}_1, \\
    \ps{a}_2 + \ps{b}_1 \ps{b}_3 \nov^n & \ps{x} = \ps{a}_2, \\
    \ps{a}_3 + \ps{b}_2 \ps{b}_1 \nov^n & \ps{x} = \ps{a}_3, \\
    \ps{d}_1 + \ps{a}_1 \ps{b}_1 \nov + \ps{b}_3 \ps{a}_3 \nov + \ps{b}_3
    \ps{b}_2 \ps{b}_1 \nov^{n+1} & \ps{x} = \ps{d}_1, \\
    \ps{d}_2 + \ps{a}_2 \ps{b}_2 \nov + \ps{b}_1 \ps{a}_1 \nov + \ps{b}_1
    \ps{b}_3 \ps{b}_2 \nov^{n+1} & \ps{x} = \ps{d}_2, \\
    \ps{d}_3 + \ps{a}_3 \ps{b}_3 \nov + \ps{b}_2 \ps{a}_2 \nov + \ps{b}_2
    \ps{b}_1 \ps{b}_3 \nov^{n+1} & \ps{x} = \ps{d}_3, \\
    \ps{x} & \text{otherwise.}
  \end{cases}
\end{equation}

\begin{lem} \label{lem:triple-pt}
  The map $\phi$ is a tame isomorphism of DGAs.
\end{lem}

\begin{proof}
  The map $\phi$ is the composition of three tame isomorphisms of
  algebras:
  \begin{enumerate}
  \item $\boldsymbol{\phi}_1 (\ps{d}_i) = \ps{d}_i + \ps{a}_i \ps{b}_i
    \nov$ for $i=1,2,3$.  It is the identity on all other generators.
  \item $\boldsymbol{\phi}_2 (\ps{a}_i) = \ps{a}_i + \ps{b}_{i-1}
    \ps{b}_{i-2}$, where $i-1$ and $i-2$ are interpreted cyclically.
    It is the identity elsewhere.
  \item $\boldsymbol{\phi}_3(\ps{d}_i) = \ps{d}_i + \ps{b}_{i-1}
    \ps{a}_{i-1} \nov$, where $i-1$ is interpreted cyclically.  It is
    the identity on all other generators.
  \end{enumerate}
  Each of these maps preserves the filtration by
  Lemma~\ref{lem:tame-isom}.  The first and last maps preserve grading
  since, by (\ref{eqn:shifted-grading}),
  \begin{equation*}
    \begin{split}
      |a^k| + |b^l| &= (k+l+1)\mu_E - 1 \\
      &= |d^{k+l+1}|.
    \end{split}
  \end{equation*}
  For $\phi_2$, the fact that $\df \ps{b}_1 = \ps{a}_2 \ps{a}_3
  T^{1-n} + \cdots$ implies
  \begin{equation*}
    |b_1^{k+l+1-n}| = |a_2^k| + |a_3^l| + 1.
  \end{equation*}
  Hence, again using (\ref{eqn:shifted-grading}),
  \begin{equation*}
    |a_1^{k+l+n}| = |b_2^k| + |b_3^l|.
  \end{equation*}
  It follows that $\phi$ itself is a tame isomorphism of filtered
  graded algebras.

  \begin{figure}[tbp]
    \centering{\input{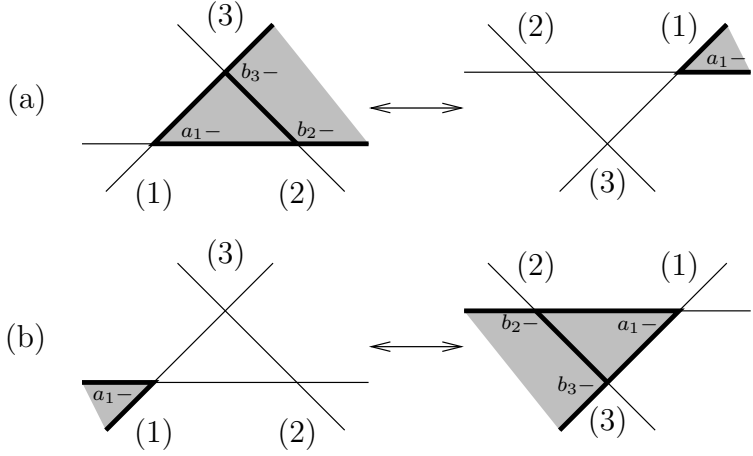}}
    \caption{Disks for differentials of generators that lie outside 
      the diagram.}
    \label{fig:triple-pt-external}
  \end{figure}
  
  Let \ps{x} be a generating power series for a special point that
  lies outside of Figure~\ref{fig:triple-pt}. To show that $\phi$ is a
  chain map on \ps{x} with respect to \dfe, consider the disks in
  Figure~\ref{fig:triple-pt-external}.  The disks in $\dfe \ps{x}$
  shown in Figure~\ref{fig:triple-pt-external}(a) give rise to terms
  containing $\ps{a}_1 + \ps{b}_3 \ps{b}_2 \nov^n$.  The $\ps{b}_3
  \ps{b}_2 \nov^n$ term in $\boldsymbol{\phi}(\ps{a}_1)$ cancels the
  one in $\dfe \ps{x}$, so the image under $\phi$ of the terms coming
  from the left side of Figure~\ref{fig:triple-pt-external}(a)
  correspond precisely to the terms in $\dfe' \ps{x}$ coming from the
  right side.  For the disks pictured in
  Figure~\ref{fig:triple-pt-external}(b), $\dfe \ps{x}$ has terms with
  $\ps{a}_1$ in them, while $\dfe' \ps{x}$ has terms with $\ps{a}_1 +
  \ps{b}_3 \ps{b}_2 \nov^n$.  This time, $\phi$ inserts $\ps{b}_3
  \ps{b}_2 \nov^n$ into $\dfe \ps{x}$ to give $\dfe' \ps{x}$.
  Symmetric arguments apply to disks in $\dfe \ps{x}$ that differ by a
  rotation from those in Figure~\ref{fig:triple-pt-external}.
  
  The map $\phi$ only affects $\dfi \ps{x}$ if $x$ is a $c$ generator
  situated immediately outside the diagram.  Suppose that $c$ lies on
  the upper right strand of Figure~\ref{fig:triple-pt}.  Calculate
  directly with the flowlines that appear in the figure:
  \begin{equation*}
    \begin{split}
      \boldsymbol{\phi} \dfi \ps{c} &= \boldsymbol{\phi} ((1+
      \ps{c})(\ps{d}_1 + \ps{b}_3 \ps{a}_3 \nov))
      \\
      &= (1+\ps{c})\bigl((\ps{d}_1 + \ps{a}_1 \ps{b}_1 \nov + \ps{b}_3
      \ps{a}_3
      \nov + \ps{b}_3 \ps{b}_2 \ps{b}_1 \nov^{n+1}) \\
      &\quad + (\ps{b}_3
      \ps{a}_3 \nov + \ps{b}_3 \ps{b}_2 \ps{b}_1 \nov^{n+1})\bigr) \\
      &= (1+\ps{c})(\ps{d}_1 + \ps{a}_1 \ps{b}_1 \nov) \\
      &= \dfi' \boldsymbol{\phi} \ps{c}.
    \end{split}
  \end{equation*}
  The other cases are symmetric, so $\phi$ is a chain map on \ps{x}.
  
  Next, consider the generating power series $\ps{a}_1$; as usual, the
  proofs for the other $\ps{a}_i$ are symmetric.  The left- and
  right-hand sides of the equation $\boldsymbol{\phi} \df \ps{a}_1 =
  \df' \boldsymbol{\phi} \ps{a}_1$ may be broken up as follows:
  \begin{equation} \label{eqn:df-a-splitting}
    \begin{split}
      \boldsymbol{\phi} \df \ps{a}_1 &= \ps{U}_{int} + \ps{U}_{ext},\\
      \df' \boldsymbol{\phi} \ps{a}_1 &= \ps{V}_{int} + \ps{V}_{ext} +
      \ps{W}_{int} + \ps{W}_{ext} + \ps{W}_{\Delta}.
    \end{split}
  \end{equation}
  The terms are:
  \begin{description}
  \item[$\ps{U}_{int}, \ps{U}_{ext}$]  the terms that make
    up $\boldsymbol{\phi} \dfi \ps{a}_1$ and $\boldsymbol{\phi} \dfe
    \ps{a}_1$, respectively.
  \item[$\ps{V}_{int}, \ps{V}_{ext}$]  the terms that make up
    $\dfi' \ps{a}_1$ and $\dfe' \ps{a}_1$, respectively.
  \item[$\ps{W}_{int}$]  the terms that make up $\dfi'
    (\ps{b}_3 \ps{b}_2 \nov^n)$.
  \item[$\ps{W}_{ext}$]  the terms that make up $\dfe'
    (\ps{b}_3 \ps{b}_2 \nov^n)$ except those that come from the
    central triangle in Figure~\ref{fig:triple-pt}.
  \item[$\ps{W}_{\Delta}$]  the remaining terms in $\dfe' (\ps{b}_3
    \ps{b}_2 \nov^n)$ that come from the central triangle.
  \end{description}

  \begin{claim}
    The terms in $\ps{U}_{ext}$ correspond to those in $\ps{V}_{ext} +
    \ps{W}_{ext}$.
  \end{claim}

  \begin{figure}[tbp]
    \centering{\input{figs/triple-pt-a1.pstex_t}}
    \caption{(a and b)  $\ps{U}_{ext}$; (c) $\ps{V}_{ext}$.}
    \label{fig:triple-pt-a1}
  \end{figure}
  
  Half of the disks in $\dfe \ps{a}_1$ are pictured in
  Figure~\ref{fig:triple-pt-a1}(a,b); the other half are reflections
  of those in the figure across the line through $a_1$ and $d_2$.  The
  terms represented by the disks in Figure~\ref{fig:triple-pt-a1}(a)
  are unchanged both by the triple-point move and by $\phi$, and hence
  appear in both $\ps{U}_{ext}$ and $\ps{V}_{ext}$.\footnote{For disks
    in $\dfe \ps{a}_1$ that come back to the triple point, use the
    arguments for \ps{x} above.}  The disks that appear in
  Figure~\ref{fig:triple-pt-a1}(b) give rise to terms of the form
  $\ps{Ub}_2$.  On the other hand, the disks in $\ps{V}_{ext}$ are of
  the form $\ps{b}_3 \ps{V}$, as shown in
  Figure~\ref{fig:triple-pt-a1}(c).  By the Leibniz rule for $\dfe'$,
  $\ps{W}_{ext}$ consists of terms of the form $\ps{Ub}_2 + \ps{b}_3
  \ps{V}$.  The defect for the \ps{U} term in $\ps{W}_{ext}$ differs
  from that of the $\ps{Ub}_2$ term in $\ps{U}_{ext}$ by $\nov^{-n}$,
  so the exponents of \nov\ agree.  This finishes the proof of the
  claim.

  \begin{claim}
    The terms in $\ps{U}_{int}$ correspond with those in $\ps{V}_{int}
    + \ps{W}_{int} + \ps{W}_{\Delta}$.
  \end{claim}
  
  The claim follows from direct computation:
  \begin{equation*}
    \begin{split}
      \ps{U}_{int} &= \boldsymbol{\phi}( \ps{d}_1 \ps{a}_1 +
      \ps{a}_1 \ps{d}_2) \\
      &= (\ps{d}_1 \ps{a}_1 + \ps{a}_1 \ps{d}_2) + (\ps{b}_3 \ps{a}_3
      \ps{a}_1 \nov + \ps{a}_1 \ps{a}_2 \ps{b}_2 \nov) \\
      &\quad + (\ps{d}_1 \ps{b}_3 \ps{b}_2 + \ps{b}_3 \ps{b}_2
      \ps{d}_2 + \ps{b}_3 \ps{a}_3 \ps{b}_3 \ps{b}_2\nov + \ps{b}_3
      \ps{b}_2 \ps{a}_2 \ps{b}_2\nov) \nov^n \\
      &= \ps{V}_{int} + \ps{W}_\Delta + \ps{W}_{int}.
    \end{split}
  \end{equation*}
  This finishes the proof that $\phi$ is a chain map on $\ps{a}_1$.

  \begin{figure}[tbp]
    \centering{\input{figs/triple-pt-b1.pstex_t}}
    \caption{The disks involved in $\dfe \ps{b}_1$.}
    \label{fig:triple-pt-b1}
  \end{figure}
  
  Next, consider the generating power series $\ps{b}_1$.  Aside from
  the terms that do not come from the inner triangle, the external
  differentials before and after the triple-point move come from the
  disks pictured in Figure~\ref{fig:triple-pt-b1}.  These disks are
  unchanged by the triple-point move and by $\phi$.
  
  For the internal differentials and the terms of the external
  differential that come from the inner triangle, compute:
  \begin{equation*}
    \begin{split}
      \boldsymbol{\phi} \df \ps{b}_1 &= \boldsymbol{\phi}( \ps{b}_1
      \ps{d}_1 + \ps{d}_2
      \ps{b}_1 + \ps{b}_1 \ps{a}_1 \ps{b}_1 \nov + \ps{a}_2 \ps{a}_3) \\
      &= \ps{b}_1 \ps{d}_1 + \ps{d}_2 \ps{b}_1 + \ps{b}_1 \ps{a}_1
      \ps{b}_1 \nov + \ps{a}_2 \ps{a}_3 \\
      &= \df' \boldsymbol{\phi} \ps{b}_1.
    \end{split}
  \end{equation*}
  Despite appearances, the equality between the first and second lines
  involves some cancellations.  This completes the proof that $\phi$
  is a chain map on $\ps{b}_1$.
  
  The calculation for the generating power series $\ps{d}_i$ is
  similar to that for the $\ps{a}_i$ and $\ps{b}_i$ generating power
  series. Lemma~\ref{lem:triple-pt} follows.
\end{proof}

For the second case in Figure~\ref{fig:triple-pt}, the tame
isomorphism is given by:
\begin{equation}
  \boldsymbol{\phi}(\ps{x}) = \begin{cases}
    \ps{a}_1 (1 + \ps{c}_1) & \ps{x} = \ps{a}_1 \\
    \ps{a}_2 + \ps{a}_1 (1 + \ps{c}_1) \ps{a}_3 & \ps{x} = \ps{a}_2 \\
    (1 + \ps{c}_1) \ps{a}_3 & \ps{x} = \ps{a}_3 \\
    (1+\ps{c}_1)^{-1}(\ps{b}_1 + (1+\ps{c}_1) \ps{a}_3 \ps{b}_2) &
    \ps{x} = \ps{b}_1 \\
    \ps{b}_3 (1+\ps{c}_1)^{-1}  + \ps{b}_2 \ps{a}_1 & \ps{x} = \ps{b}_3
    \\
    \ps{d}_2 + \ps{a}_1 \ps{b}_1 + (\ps{a}_2 + \ps{a}_1 (1 + \ps{c}_1)
    \ps{a}_3) \ps{b}_2 & \ps{x} = \ps{d}_2
    \\
    \ps{d}_3 + \ps{b}_2 \ps{a}_2 + (\ps{b}_3 (1+\ps{c})^{-1} +
    \ps{b}_2 \ps{a}_1) (1+\ps{c}) \ps{a}_3 & \ps{x} = \ps{d}_2 \\
    \ps{x} & \text{otherwise.}
  \end{cases}
\end{equation}
The proof that this is a tame isomorphism of DGAs is similar to the
proof of Lemma~\ref{lem:triple-pt}.

\subsubsection{Double-Point Move}
\label{sec:double-pt-invariance}

\begin{figure}[tbp]
  \centerline{\input{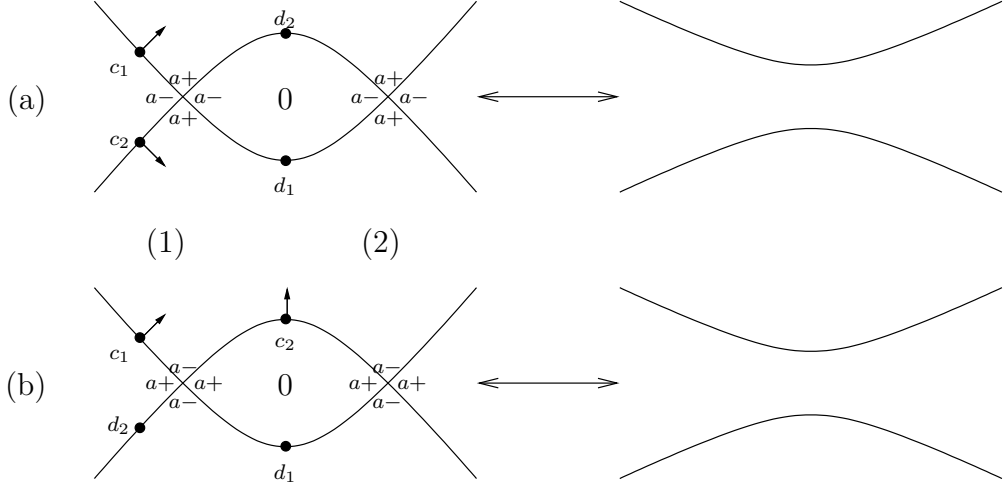}}
  \caption{Two types of possible diagrams for the double-point move.}
  \label{fig:double-pt}
\end{figure}

The last step in proving Theorem~\ref{thm:invariance} is to show that
the double point move does not change the stable tame isomorphism type
of a knot diagram's DGA.  Up to a choice of transverse directions at
the $c$ points and an overall shift in the $c$ and $d$ labels, there
are two possibilities for the labeling of the diagram; see
Figure~\ref{fig:double-pt}.  Only the first possibility is considered
in this section; the proof for the second is almost identical.

Let $(\alg, \df)$ be the DGA for the diagram on the left of
Figure~\ref{fig:double-pt}(a), and let $(\alg', \df')$ be the DGA for
the diagram on the right.  Define the following four stabilizations:
\begin{enumerate}
\item $\mathcal{E}_1 = \mathcal{E}(\boldsymbol{\alpha}_2,
  \boldsymbol{\beta}_1)$, with grading given by $|\beta^k_1| = |b^k_1|$.
\item $\mathcal{E}_2 = \mathcal{E}(\boldsymbol{\alpha}_1,
  \boldsymbol{\beta}_2)$, with grading given by $|\beta^k_2| = |b^k_2|$.
\item $\mathcal{E}_3 = \mathcal{E}(\boldsymbol{\delta}_1,
  \boldsymbol{\gamma}_1)$, with grading given by $|\gamma^k_1| = |c^k_1|$.
\item $\mathcal{E}_4 = \mathcal{E}(\boldsymbol{\delta}_2,
  \boldsymbol{\gamma}_2)$, with grading given by $|\gamma^k_2| = |c^k_2|$.
\end{enumerate}
Let $(S(\alg'), \df')$ be the result of stabilizing $\alg'$ with
$\mathcal{E}_1, \ldots, \mathcal{E}_4$.  The rest of this section will
carry out a five-step program to show that $(\alg, \df)$ is tame
isomorphic to $(S(\alg'), \df')$.  The first steps construct four
intermediate DGAs $(\alg_i, \df_i)$, $i=1,\ldots,4$, that are stable
tame isomorphic to $(\alg, \df)$; the last step shows that $(\alg_4,
\df_4) \simeq (S(\alg'), \df')$.

\emph{Step 1}

Let $\hat{\alg}$ be an algebra generated by the same power series
(except for $\ps{a}_2$ and $\ps{b}_1$) as \alg.

\begin{defn}
  \label{defn:alg_1}
  As an algebra, $\alg_1 = S_{\mathcal{E}_1}(\hat{\alg})$.
\end{defn}

To define a differential on $\alg_1$, let $W^k \in \alg$ be given by:
\begin{equation} \label{eqn:df-b_1}
  \df b_1^k = a_2^k + W^k
\end{equation}
Define $W_R^k \in \alg_1$ by the following inductive procedure: first,
let $W_R^0 = W^0$. Since the curvature over the lune in
Figure~\ref{fig:double-pt} is arbitrarily small, $l(a^k_2)$ is
arbitrarily close to $l(b^k_1)$.  Thus, any generator $w$ that appears
in $W^k$ must satisfy $l(w) < l(a_2^k)$.  In particular, any generator
$a_2^j$ appearing in $W^k$ must have $j < k$.  To define $W_R^k$,
replace any $a_2^j$ that appears in $W^k$ by $W_R^j$ and any $b_1^j$
by $0$.

\begin{defn}
  \label{defn:df_1}
  The differential $\df_1$ on $\alg_1$ is 
  \begin{equation}
    \df_1 \ps{x} = \df \ps{x} \Bigl\lvert_{\substack{a^k_2 = W^k_R 
        \\ b_1^k = 0}} \Bigr. \quad \text{for all } k \geq 0.
  \end{equation}
\end{defn}

That $\df_1 \circ \df_1 = 0$ will be proven by showing that there is a
tame isomorphism $\Phi_1$ of algebras between \alg\ and $\alg_1$ that
intertwines \df\ and $\df_1$.  In order to define $\Phi_1$, construct
a sequence of algebras $\alg^k$, $k=0, 1, \ldots$, where:
\begin{equation}
  F^j \alg^k = \begin{cases}
    F^j \alg_1 & 0 \leq j < k, \\
    F^j \alg & j \geq k.
  \end{cases}
\end{equation}
In other words, the generators $a_2^0, \ldots, a_2^{k-1}$ in \alg\ are
replaced by $\alpha_2^0, \ldots, \alpha_2^{k-1}$ and the generators
$b_1^0, \ldots, b_1^{k-1}$ are replaced by $\beta_1^0, \ldots,
\beta_1^{k-1}$.  Furthermore, $\alg^0 = \alg$ and $\alg_1$ is the
direct limit of the $\alg^k$. Each algebra in the sequence has a
differential defined by:
\begin{equation}
  \df^k \ps{x} = \df \ps{x} \Bigl\lvert_{\substack{a^j_2 = W^j_R 
      \\ b_1^j = 0}} \Bigr.
\end{equation}
for $j = 0, 1, \ldots, k-1$.  Note that the
projection operator $\tau: \alg_1 \to \hat{\alg}$ and the homotopy
operator $H: \alg_1 \to \alg_1$ defined in
Section~\ref{sec:isomorphisms} descend to the algebras $\alg^k$.  They
still satisfy:
\begin{equation} \label{eqn:new-H-chain}
      \tau \circ i + Id_{\hat{\alg}} = H \circ \df^k + \df^k \circ H.
\end{equation}

The next step in defining $\Phi_1$ is to construct tame isomorphisms
$\phi^k: \alg^k \to \alg^{k+1}$.  First, define an order $\prec$ on
the generators of $\alg^k$ by length, with $d^j \prec c^j$ and with
$\alpha_2^j$ and $\beta_1^j$ replacing $a_2^j$ and $b_1^j$,
respectively, in the ordering. Make the curvature of the lune small
enough so that there are no generators that lie between $a_2^k$ and
$b_1^k$.  Lemma~\ref{lem:bounded-curvature} and the definition of
\dfi\ then imply that every generator that appears in $\df^k x$
precedes $x$ in the ordering.

Define $\phi_0^k: \alg_1^k \to \alg_1^{k+1}$ by:
\begin{equation} \label{eqn:double-pt-psi}
  \phi_0^k(x^j) = \begin{cases}
    \beta_1^k & x^j = b_1^k, \\
    \alpha_2^k + W_R^k & x^j = a_2^k, \\
    x^j & \text{otherwise.}
  \end{cases}
\end{equation}
This is an elementary isomorphism.  Next, use this map to define
$\phi^k:\alg^k \to \alg^{k+1}$ by:
\begin{equation}
  \phi^k(x^j) = \begin{cases}
    \phi_0^k(x^j) & x^j = a_2^k, b_1^k, \\
    x^j & x^j = \alpha_2^j, \beta_1^j, \\
    x^j + H \phi_0^k \df^k x^j & \text{otherwise.}
  \end{cases}
\end{equation}

\begin{lem} \label{lem:inductive-chain}
  The map $\phi^k$ is a tame isomorphism of algebras that satisfies
  \begin{equation} \label{eqn:phi-chain}
    \df^{k+1} \phi^k = \phi^k \df^k.
  \end{equation}
  In particular, $\df^k$ is a differential on $\alg^k$.
\end{lem}

\begin{proof}[Proof]
  The map $\phi^k$ is a tame isomorphism of algebras by the remark
  after Lemma~\ref{lem:tame-isom}.
  
  To prove (\ref{eqn:phi-chain}), there are four types of generators
  $x^j$ to consider: $x^j \prec a_2^k$, $x^j = a_2^k$, $x^j = b_1^k$,
  and $x^j \succ b_1^k$.  The case when $x^j \prec a_2^k$ follows from
  two facts.  First, $a_2^l$ does not appear in $\df^k x^j$ for all $l
  \leq j$.  Thus, $H \phi^k_0 \df^k x^j = 0$, and hence $\phi^k(x^j) =
  x^j$.  Second, $\df^{k+1} x^j = \df^k x^j$ for $j \leq k$.  Putting
  these together yields:
  \begin{align*}
    \df^{k+1} \phi^k x^j &= \df^{k+1} x^j & \text{by the first fact} \\
    &= \df^k x^j & \text{by the second fact}\\
    &= \phi^k \df^k x^j & \text{by the first fact.}
  \end{align*}
  This proves the first case.
  
  For the second case, the following lemma is required:

  \begin{lem} \label{lem:df-b_1}
    $\df_1 a_2^k = \df_1 W_R^k$.
  \end{lem}
  
  \begin{proof}
    Since $\df \circ \df = 0$, (\ref{eqn:df-b_1}) implies, for all
    $k$, that:
    \begin{equation} \label{eqn:lemma-df-b1}
      \df a_2^k = \df W^k.
    \end{equation}
    The remainder of the proof is an induction on $k$.  For $k=0$,
    $W^0 = W^0_R$ by definition.  Thus, $\df_1 W^0 = \df_1 W^0_R$.  In
    general, a typical term in $W^k$ has the form $x_1 a_2^{j_1}x_2
    \cdots x_n a_2^{j_n}x_{n+1}$ with $j_i < k$.  In $W^k_R$, this
    term becomes $x_1 W_R^{j_1}x_2 \cdots x_n W_R^{j_n}x_{n+1}$.  The
    Leibniz rule and the inductive hypothesis imply:
    \begin{equation*}
      \begin{split}
        \df_1 (x_1 a_2^{j_1}x_2 \cdots x_n a_2^{j_n}x_{n+1}) &= \bigl(
        \df (x_1) a_2^{j_1} \cdots + x_1 \df (a_2^{j_1}) \cdots +
        \cdots \bigr) \Bigl\lvert_{\substack{a^j_2 = W^j_R \\ b_1^j =
            0}} \Bigr. \\
        &= \df (x_1) W_R^{j_1} \cdots + x_1 \df (W_R^{j_1}) \cdots +
        \cdots \\
        &= \df_1 (x_1 W_R^{j_1}x_2 \cdots x_n W_R^{j_n}x_{n+1}).
      \end{split}
    \end{equation*}
  \end{proof}
  To complete the second case, compute directly:
  \begin{align*}
    \df^{k+1} \phi^k a_2^k &= \df^{k+1} (\alpha_2^k + W_R^k) \\
    &= \df^{k+1} W_R^k && \text{since } \df^{k+1} \alpha_2^k = 0; \\
    &= \phi^k \df^k W_R^k && \text{by the proof of the first case}; \\
    &= \phi^k \df^k a_2^k && \text{by Lemma~\ref{lem:df-b_1}}.
  \end{align*}
  
  The third case also follows from a straightforward computation.
  
  For the final case, proceed by induction on the order $\prec$.  To
  begin, use (\ref{eqn:new-H-chain}) to get:
  \begin{align*}
    \phi^k \df^k x^j &= \tau \phi^k \df^k x^j + \df^{k+1} H \phi^k
    \df^k x^j + H \df^{k+1} \phi^k \df^k x^j. \\
    \intertext{Lemma~\ref{lem:bounded-curvature}, the definition of
      \dfi\, and the inductive hypothesis give:} &= \tau \phi^k \df^k
    x^j + \df^{k+1} H \phi^k
    \df^k x^j + H \phi^k \df^k \df^k x^j. \\
    \intertext{For every $a_2^k$ appearing in the image of $\df^k$,
      $\tau \phi^k$ substitutes in $W_R^k$. This is the same as $\tau
      \df^{k+1}$.  Thus:}
    &= \tau \df^{k+1} x^j + \df^{k+1} H \phi^k \df^k x^j. \\
    \intertext{Another application of (\ref{eqn:new-H-chain}) gives:}
    &= \df^{k+1} (x^j + H(\df^{k+1}x^j + \phi^k \df^k x^j)) + H
    \df^{k+1} \df^{k+1} x^j. \\
    \intertext{The second and last terms disappear, leaving:} &=
    \df^{k+1} \phi^k x^j.
  \end{align*}
  This completes the final case, and hence the proof of the lemma.
\end{proof}

With this machinery in hand, define $\Phi_1: \alg \to \alg_1$ by:
\begin{equation} \label{eqn:phi-1}
  \Phi_1 = \cdots \circ \phi^2 \circ \phi^1 \circ \phi^0.
\end{equation}
The map $\Phi_1$ satisfies the finiteness requirement for a tame
isomorphism of algebras since $\phi^k$ is the identity on
$F^{k-1}\alg$.  That $\Phi_1$ is a tame isomorphism of DGAs follows
from Lemma~\ref{lem:inductive-chain}.

\emph{Steps 2, 3, and 4}

Let $\hat{\alg_1}$ be generated by the same power series (except for
$\ps{a}_1$ and $\ps{b}_2$) as $\alg_1$. Suppose that:
\begin{equation} \label{eqn:df-b_2}
  \df_1 b_2^k = a_1^k + V^k.
\end{equation}
Define $V_R^k$ from $V^k$ by the same inductive procedure that
produced $W_R^k$.

\begin{defn}
  \label{defn:alg_2}
  As an algebra, $\alg_2 = S_{\mathcal{E}_2}(\hat{\alg_1})$.  The
  differential $\df_2$ on $\alg_2$ is defined by:
  \begin{equation}
    \df_2 \ps{x} = \df_1 \ps{x} \Bigl\lvert_{\substack{a^k_1 = V^k_R 
        \\ b_2^k = 0}} \Bigr.
  \end{equation}
  for all $k \geq 0$.
\end{defn}

The same procedure as in Step 1 provides a tame isomorphism $\Phi_2$
between $(\alg_1, \df_1)$ and $(\alg_2, \df_2)$.

Build $(\alg_3, \df_3)$ and $(\alg_4, \df_4)$ in much the same way:
first repeat the construction on $\alg_2$ with $c_1$ in the place of
$b_1$, $d_1$ in the place of $a_2$, and $X^k_R$ in the place of
$V^k_R$ to obtain $(\alg_3, \df_3)$.  Note that $X^k_R$ comes from
flowlines that start at $c_1$ and leave Figure~\ref{fig:double-pt}
along the upper left strand. Again, there exists a tame isomorphism
$\Phi_3: (\alg_2, \df_2) \to (\alg_3, \df_3)$.

Next, repeat the construction on $\alg_3$ with $c_2$ in the place of
$c_1$, $d_2$ in the place of $d_1$, and $Y^k_R$ in the place of
$X^k_R$ to obtain $(\alg_4, \df_4)$. As before, there exists a tame
isomorphism $\Phi_3: (\alg_3, \df_3) \to (\alg_4, \df_4)$.

\emph{Step 5}

The final step is to show that $\alg_4 \simeq S(\alg')$ as DGAs.  They
have the same generators by construction, so they are isomorphic as
algebras; it remains to show that $\df_4 = \df'$.

Let $x^k$ be a generator of $\alg$ that does not appear in
Figure~\ref{fig:double-pt}.  Write the differential of $x^k$ in $\alg$
as:
\begin{equation}
  \df x^k = W_0 + W_1 + W_{int} + W_{ext}.
\end{equation}
The terms on the right hand side are as follows:
\begin{itemize}
\item $W_0$ consists of words that do not contain any of the
  generators that appear in the left hand side of
  Figure~\ref{fig:double-pt}(a).
\item $W_1$ consists of words that contain either $b_i^j$ or $c_i^j$,
  $i=1,2$.
\item $W_{int}$ consists of words in $\dfi(x^k)$ that contain $d_i^j$,
  $i=1,2$.
\item $W_{ext}$ consists of words in $\dfe(x^k)$ that contain $a_i^j$
  but not $b_i^l$ or $c_i^l$, $i=1,2$.
\end{itemize}

Similarly, the differential of $x^k$ in $\alg'$ may be written as:
\begin{equation}
  \df' x^k = W_0 + W_{int}' + W_{ext}'.
\end{equation}
The two new terms on the right hand side are:
\begin{itemize}
\item $W_{int}'$, consisting of words coming from flowlines that pass
  through the right hand side of Figure~\ref{fig:double-pt}(a); and
\item $W_{ext}'$, consisting of words coming from disks that pass
  ``through the neck'' of the right hand side of
  Figure~\ref{fig:double-pt}.
\end{itemize}
The task is to prove that, under the substitutions in Steps 1 through
4 that transform \df\ into $\df_4$, $W_1 + W_{int} + W_{ext}$ becomes
$W_{int}' + W_{ext}'$.  

First of all, all words in $W_1$ contain generators that get mapped to
$0$ by one of $\df_1, \ldots, \df_4$.  Hence, $W_1$ does not appear in
$\df_4 x^k$ in any form.

\begin{figure}[tbp]
  \centering{\input{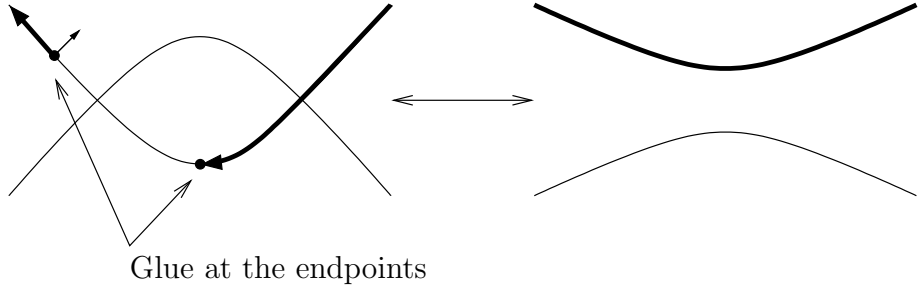}}
  \caption{Gluing two flowlines together to form a flowline that 
    passes through the diagram.}
  \label{fig:flowline-glue}
\end{figure}

Secondly, $W_{int}$ transforms into $W_{int}'$.  To see this, consider
the case of a generator that lies just off of the top right corner of
Figure~\ref{fig:double-pt}(a).\footnote{The case of a generator that
  lies off of the bottom right corner is symmetric.  The interior
  differentials of generators off of the left side do not involve any
  of the generators in Figure~\ref{fig:double-pt}(a).}  The terms with
$d_1^k$ in $W_{int}$ are untouched until Step 3, when $d^k_1$ is
replaced by $X^k_R$.  As noted above, $X^k_R$ comes from flowlines
that leave Figure~\ref{fig:double-pt}(a) via the top left strand.  The
result of the substitution of $X^k_R$ for $d^k_1$, then, is that the
flowlines coming into $d_1$ from the right are glued to flowlines
leaving $c_1$ to the left.  As shown in
Figure~\ref{fig:flowline-glue}, these form flowlines that pass through
the diagram in Figure~\ref{fig:double-pt}(a).  These are precisely the
flowlines that give $W_{int}'$.

\begin{figure}[tbp]
  \centering{\input{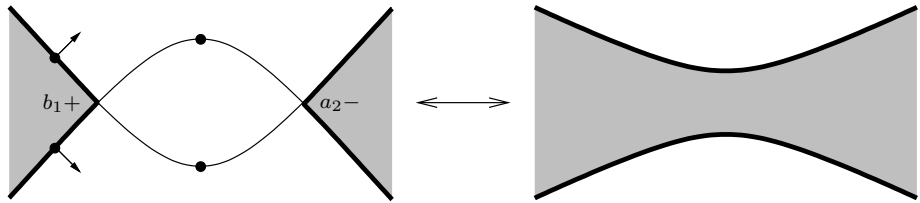}}
  \caption{Gluing two disks together to get a disk that passes through 
    the neck.}
  \label{fig:disk-glue}
\end{figure}

Thirdly, $W_{ext}$ transforms into $W_{ext}'$.  For conceptual
clarity, assume that $W_{ext}$ and $W^k_R$ are both monomials, i.e.
each comes from a single disk.  The general case is a sum over all the
disks involved in the constructions described below.  The first step
is to describe $W^k_R$ geometrically.  Recall that in the inductive
definition of $W^k_R$, every occurrence of $a^j_2$ in $W$ is replaced
by $W^j_R$. The substitution of $W^j_R$ can be represented
geometrically by the gluing pictured in Figure~\ref{fig:disk-glue}.
Thus, $W^k_R$ should be thought of as a disk that passes through the
neck on the right side of Figure~\ref{fig:disk-glue} (possibly
multiple times) that has a ``phantom'' positive corner at $b_1^j$ that
can be glued to a negative corner at $a_2^j$.  The construction of
$V^j_R$ is similar.

With the geometric description of $W^k_R$ in hand, suppose that
$W_{ext}$ has a negative corner at $a_2^j$.  In the definition of
$\df_1$, $a_2^j$ is replaced by $W^j_R$.  As before, this is the
algebraic realization of gluing $W^j_R$ to $W_{ext}$ as in
Figure~\ref{fig:disk-glue}. After replacing all generators $a_2^j$
with $W^j_R$ and $a_1^j$ with $V^j_R$, the result is a disk that
passes through the neck (possibly multiple times).  In other words,
the result is a summand of $W_{ext}'$.  Conversely, any disk in
$W_{ext}$ may be constructed in this manner by squeezing off the neck
of a disk in $W_{ext}'$.  Thus, $\df_4 x^k = \df' x^k$.

Finally, since $\df_4 = \df'$ on the $\alpha$, $\beta$, $\gamma$, and
$\delta$ generators, the arguments above show that $\df_4 = \df'$ on
all of $S(\alg')$.  This finishes the proof that $\alg_4$ is tame
isomorphic to $S(\alg')$ as DGAs, and hence the proof of
Theorem~\ref{thm:invariance}.


\section*{Acknowledgments}

This paper stems from my thesis research, and I would like to thank
Yasha Eliashberg for his insight and expert guidance.  Additionally, I
have benefited greatly from discussions and correspondence with Lenny
Ng, John Etnyre, and Frederic Bourgeois.

\bibliography{s1bundles}

\providecommand{\bysame}{\leavevmode\hbox to3em{\hrulefill}\thinspace}
\begin{thebibliography}{10}

\bibitem{aeb}
B.~Aebischer et~al., \emph{Symplectic geometry}, Prog. Math., vol. 124,
  Birkh\"auser, 1994.

\bibitem{austin-braam}
D.~M. Austin and P.~J. Braam, \emph{Morse-{B}ott theory and equivariant
  cohomology}, The Floer Memorial Volume, Birkh\"auser, Basel, 1995,
  pp.~123--183.

\bibitem{bennequin}
D.~Bennequin, \emph{Entrelacements et equations de {P}faff}, Asterisque
  \textbf{107--108} (1983), 87--161.

\bibitem{bourgeois:mb}
F.~Bourgeois, \emph{A {M}orse-{B}ott approach to contact homology}, Preprint,
  2002.

\bibitem{chv}
Yu. Chekanov, \emph{Differential algebras of {L}egendrian links}, Invent. Math.
  (2002), To Appear.

\bibitem{yasha:overtwisted}
Y.~Eliashberg, \emph{Classification of overtwisted contact structures on
  $3$-manifolds}, Invent. Math. \textbf{98} (1989), no.~3, 623--637.

\bibitem{yasha:filling}
\bysame, \emph{Filling by holomorphic discs and its applications}, Geometry of
  low-dimensional manifolds, 2 (Durham, 1989), Cambridge Univ. Press,
  Cambridge, 1990, pp.~45--67.

\bibitem{yasha:20yrs}
\bysame, \emph{Contact $3$-manifolds twenty years since {J}. {M}artinet's
  work}, Ann. Inst. Fourier (Grenoble) \textbf{42} (1992), no.~1-2, 165--192.

\bibitem{yasha:knots}
\bysame, \emph{{L}egendrian and transversal knots in tight contact
  $3$-manifolds}, Topological methods in modern mathematics (Stony Brook, NY,
  1991), Publish or Perish, Houston, TX, 1993, pp.~171--193.

\bibitem{yasha-fraser}
Y.~Eliashberg and M.~Fraser, \emph{Classification of topologically trivial
  {L}egendrian knots}, Geometry, topology, and dynamics (Montreal, PQ, 1995),
  Amer. Math. Soc., Providence, RI, 1998, pp.~17--51.

\bibitem{egh}
Y.~Eliashberg, A.~Givental, and H.~Hofer, \emph{Introduction to symplectic
  field theory}, Geom. Funct. Anal. (2000), no.~Special Volume, Part II,
  560--673, GAFA 2000 (Tel Aviv, 1999).

\bibitem{confol}
Y.~Eliashberg and W.~Thurston, \emph{Confoliations}, American Mathematical
  Society, Providence, RI, 1998.

\bibitem{efm}
Judith Epstein, Dmitry Fuchs, and Maike Meyer, \emph{Chekanov-{E}liashberg
  invariants and transverse approximations of {L}egendrian knots}, Pacific J.
  Math. \textbf{201} (2001), no.~1, 89--106.

\bibitem{etnyre-honda:non-existence}
J.~Etnyre and K.~Honda, \emph{On the non-existence of tight contact
  structures}, Ann. Math. (2) \textbf{153} (2001), 749--766.

\bibitem{etnyre-honda:knots}
\bysame, \emph{Knots and contact geometry}, J. Symplectic Geom. \textbf{1}
  (2002), X--Y.

\bibitem{ens}
J.~Etnyre, L.~Ng, and J.~Sabloff, \emph{Invariants of {L}egendrian knots and
  coherent orientations}, J. Symplectic Geom. \textbf{1} (2002), no.~2,
  321--368.

\bibitem{etnyre:intro}
J.~Etynre, \emph{Introductory lectures on contact geometry}, Preprint, 2001.

\bibitem{fuchs:augmentations}
D.~Fuchs, \emph{Chekanov-{E}liashberg invariants of {L}egendrian knots:
  Existence of augmentations}, Preprint, 2000.

\bibitem{giroux:bundles}
E.~Giroux, \emph{Structures de contact sur les vari\'et\'es fibr\'ees en
  cercles audessus d'une surface}, Comment. Math. Helv. \textbf{76} (2001),
  no.~2, 218--262.

\bibitem{honda:classification-2}
K.~Honda, \emph{On the classification of tight contact structures. {II}}, J.
  Differential Geom. \textbf{55} (2000), no.~1, 83--143.

\bibitem{lutz}
R.~Lutz, \emph{Structures de contact sur les fibres principaux en cercles de
  dimension trois}, Ann. Inst. Fourier, Grenoble \textbf{27} (1977), no.~3,
  1--15.

\bibitem{lenny:computable}
L.~Ng, \emph{Computable {L}egendrian invariants}, Topology (2002), To Appear.

\bibitem{sato-tsuboi}
A.~Sato and T.~Tsuboi, \emph{Contact structures of closed manifolds fibered by
  the circles}, Mem. Inst. Sci. Tech. Meiji Univ. \textbf{33} (1994), 41--46.

\bibitem{schwarz}
Matthias Schwarz, \emph{Morse homology}, Birkh\"auser Verlag, Basel, 1993.

\bibitem{turaev:shadow}
V.~Turaev, \emph{Shadow links and face models of statistical mechanics}, J.
  Differential Geom. \textbf{36} (1992), no.~1, 35--74.

\end{thebibliography}
\bibliographystyle{amsplain}

\end{document}